\documentclass[11pt]{article}
\usepackage{amsmath}
\usepackage{amsthm}

\usepackage{tikz}
\usetikzlibrary{positioning}

\usepackage{amsfonts}
\usepackage{amssymb}
\usepackage{mathtools}
\usepackage{soul}
\usepackage{bbm}
\usepackage{hyperref}
\usepackage[compress]{cite}
\usepackage[draft]{todonotes}   

\usepackage{cases}
\usepackage{xcolor}

\newcommand{\dt}{\partial_t}

\newcommand{\dvh}{\mathrm{div}_h\,}
\newcommand{\nablah}{\nabla_h}
\newcommand{\deltah}{\Delta_h}
\newcommand{\dz}{\partial_z}
\newcommand{\idx}{\,d\vec{x}}
\newcommand{\idxh}{\,dxdy}
\newcommand{\bblbrack}{\biggl\lbrack}
\newcommand{\bbrbrack}{\biggr\rbrack}
\newcommand{\bblparenthese}{\biggl(}
\newcommand{\bbrparenthese}{\biggr)}
\newcommand{\blparenthese}{\bigl(}
\newcommand{\brparenthese}{\bigr)}
\newcommand{\subeqref}[2]{$\eqref{#1}_{#2}$}
\newcommand{\abs}[2]{\bigl| #1 \bigr|^{#2}}
\newcommand{\norm}[2]{\bigl\Arrowvert #1 \bigr\Arrowvert_{#2}}
\newcommand{\hnorm}[2]{\bigl| #1 \bigr|_{#2}}
\newcommand{\Lnorm}[1]{L^{#1}}
\newcommand{\Hnorm}[1]{H^{#1}}

\theoremstyle{definition}

\theoremstyle{theorem}
\newtheorem{lemma}{Lemma}
\newtheorem{corollary}{Corollary}
\newtheorem{proposition}{Proposition}
\newtheorem{theorem}{Theorem}[section]
\theoremstyle{remark}
\newtheorem{remark}{Remark}

\numberwithin{equation}{section}
\allowdisplaybreaks
\title{Zero Mach Number Limit of the Compressible Primitive Equations Part I: Well-prepared Initial Data
}

\author{Xin Liu\footnote{Department of Mathematics, Texas A{\&}M University, College Station, TX 77843, USA. Email: stleonliu@gmail.com} \,\, and \,  Edriss S. Titi\footnote{Department of Mathematics, Texas A{\&}M University, College Station,  TX 77840, USA.  Department of Applied Mathematics and Theoretical Physics, University of Cambridge, Cambridge CB3 0WA, UK.
		Department of Computer Science and Applied Mathematics, Weizmann Institute of Science, Rehovot 76100, Israel. Email: {titi@math.tamu.edu}\, and \, {Edriss.Titi@damtp.cam.ac.uk}}
}
\date{May 22, 2019}
\newcommand\blfootnote[1]{%
	\begingroup
	\renewcommand\thefootnote{}\footnote{#1}%
	\addtocounter{footnote}{-1}%
	\endgroup
}

\begin{document}
\allowdisplaybreaks

\maketitle
\begin{abstract}
	This work concerns the zero Mach number limit of the compressible primitive equations. The primitive equations with the incompressibility condition are identified as the limiting equations. The convergence with well-prepared initial data (i.e., initial data without acoustic oscillations) is rigorously justified, and the convergence rate is shown to be of order $ \mathcal O(\varepsilon) $, as $ \varepsilon \rightarrow 0^+ $, where $ \varepsilon $ represents the Mach number. As a byproduct, we construct a class of global solutions to the compressible primitive equations, which are close to the incompressible flows.
\end{abstract}

\tableofcontents

\section{Introduction}

\subsection{The compressible primitive equations}\label{sec:intro-cpe}
The compressible primitive equations (see \eqref{CPE}, below) are used by meteorologists to perform theoretical investigations and practical weather predictions (see, e.g., \cite{Lions1992}). In comparison with the general hydrodynamic and thermodynamic equations, the vertical component of the momentum equations is missing in the compressible primitive equations. Instead, it is replaced by
 the hydrostatic balance equation (see \subeqref{CPE}{3}, below), which is also known as the quasi-static equilibrium equation. From the meteorologists' point of view, such an approximation is reliable and useful for two reasons: the balance of gravity and pressure dominates the dynamics in the vertical direction; and the vertical velocity is usually hard to observe in reality (see, e.g., \cite[Chapter 4]{Richardson1965}). On the other hand, by formally taking the zero limit of the aspect ratio between the vertical scale and the planetary horizontal scale, the authors in \cite{Ersoy2011a} derive the compressible primitive equations from the compressible hydrodynamic equations.
Such a deviation is very common in planetary scale geophysical models, which represents the fact that the vertical scale of the atmosphere (or ocean) is significantly smaller than the planetary horizontal scale. We refer, for more comprehensive meteorological studies, to \cite{Richardson1965,Washington2005}. See also,
\cite{MajdaAtmosphereOcean,FeireislSingularLimits,Klein2000,Klein2001,RKlein2010,Klein2005}
for more general discussions on multi-scale analysis.

As far as we know, there are very few mathematical studies concerning the compressible primitive equations (referred to as CPE hereafter). Lions, Temam, and Wang first introduced CPE into the mathematical community in \cite{Lions1992}. They formulated the commonly known primitive equations (referred to as PE hereafter) with the incompressibility condition as the representation of the compressible primitive equations in the pressure coordinates ($p$-coordinates) instead of the physical ones with the vertical spatial coordinate. On the other hand,
as mentioned before, the authors in \cite{Ersoy2011a} introduce these equations  with a formal deviation, and a rigorous justification is still an open question for now. In \cite{Ersoy2011a}, the stability of weak solutions is also investigated (see also \cite{Tang2015}). The stability is meant in the sense that a sequence of weak solutions, satisfying some entropy conditions, contains a subsequence converging to another weak solution, i.e., a very weak sense of stability. The existence of such weak solutions is recently constructed in \cite{LT2018b,wang2017global} (see also  \cite{Ersoy2012,Gatapov2005} for the existence of global weak solutions to some variant of compressible primitive equations in two spatial dimension). In \cite{LT2018a}, we also construct local strong solutions to CPE in two cases: with gravity but no vacuum; with vacuum but no gravity.

In analogy to the low Mach number limit in the study of compressible hydrodynamic equations,
this and our subsequent works are aiming at studying the low Mach number limit of the compressible primitive equations without gravity and Coriolis force. It is worth mentioning that  while taking into account the Coriolis force would not change  much our proof, considering gravity in our system causes challenging difficulties.
Let $ \varepsilon $ denote the Mach number, and let $ \rho^\varepsilon \in \mathbb R, v^\varepsilon \in \mathbb R^2, w^\varepsilon \in \mathbb R $ be the density, the horizontal and the vertical velocities, respectively.
System \eqref{CPE}, below, is obtained by rescaling the original CPE, which is similar to the rescaling of the compressible Navier--Stokes equations (see, e.g., \cite{FeireislSingularLimits}):
\begin{equation*}\label{CPE}\tag{CPE}
\begin{cases}
	\dt \rho^\varepsilon + \dvh (\rho^\varepsilon v^\varepsilon) + \dz (\rho^\varepsilon w^\varepsilon) = 0 & \text{in} ~ \Omega_h \times (0,1) , \\
	\dt (\rho^\varepsilon v^\varepsilon) + \dvh( \rho^\varepsilon v^\varepsilon\otimes v^\varepsilon) + \partial_z (\rho^\varepsilon w^\varepsilon v^\varepsilon)  \\
	~~~~ ~~ + \dfrac{1}{\varepsilon^2}\nablah P(\rho^\varepsilon)
	  = \dvh \mathbb S(v^\varepsilon) +  \partial_{zz} v^\varepsilon & \text{in} ~ \Omega_h \times (0,1), \\
	\partial_z P(\rho^\varepsilon) = 0 & \text{in} ~ \Omega_h \times (0,1),
\end{cases}
\end{equation*}
where $ P(\rho^\varepsilon) = (\rho^\varepsilon)^\gamma $ and $ \mathbb S (v^\varepsilon) = \mu (\nablah v^\varepsilon + \nablah^\top v^\varepsilon) + (\lambda - \mu) \dvh v^\varepsilon  \mathbb{I}_2 $ represent the pressure potential and the viscous stress tensor, respectively, with the shear and bulk viscosity coefficients $ \mu $ and $ \lambda-\mu + \frac{2}{3}\mu = \lambda - \frac{1}{3} \mu  $.
The physical requirements of $ \mu, \lambda, \gamma $ are
$ \lambda - \frac{1}{3} \mu  > 0$, $ \mu > 0 $ and $ \gamma > 1 $.
Moreover, we focus our study on the case when
$ \Omega_h := \mathbb T^2 \subset \mathbb R^2 $
.  \eqref{CPE} is complemented with the following stress-free and  impermeability physical boundary conditions:
\begin{equation*}\tag{BC-CPE} \label{bc-cpe}
	\dz v^\varepsilon\big|_{z=0,1} = 0,~ w^\varepsilon\big|_{z=0,1} = 0.
\end{equation*}
Hereafter, we have and will use $ \nabla_h, \dvh $ and $ \Delta_h $ to represent the horizontal gradient, the horizontal divergence, and the horizontal Laplace operator, respectively; that is,
\begin{gather*}
	\nabla_h := \biggl(\begin{array}{c}
		\partial_x \\ \partial_y
	\end{array}\biggr) , ~ \dvh := \nabla_h \cdot, ~
	\Delta_h := \dvh \nabla_h.
\end{gather*}
Notice that if we consider system \eqref{CPE} subject to periodic boundary conditions, with fundamental periodic domain $\Omega_h \times 2\mathbb T$, where $ 2 \mathbb T $ denotes the periodic domain in $ \mathbb R $ with period $ 2 $; and  subject to the following symmetry:
\begin{equation*}\tag{SYM-CPE}\label{SYM-CPE}
	\text{$ v^\varepsilon $ and $ w^\varepsilon $ are even and odd, respectively, in the $ z $-variable,}
\end{equation*}
then its solutions obey the above symmetry. Moreover, when such solutions are restricted to the physical domain $\Omega_h \times (0,1)$ they automatically satisfy the physical boundary conditions \eqref{bc-cpe}.

We recall the incompressible primitive equations:
\begin{equation*}\tag{PE}\label{PE0}
\begin{cases}
	\dvh v_p + \dz w_p = 0 & \text{in} ~ \Omega_h \times (0,1), \\
	\rho_0 (\dt v_p + v_p \cdot \nablah v_p + w_p \dz v_p ) + \nablah (c^2_s \rho_1) \\
	~~~~ ~~~~ = \mu \deltah v_p + \lambda \nablah \dvh v_p +  \partial_{zz}v_p & \text{in} ~ \Omega_h \times (0,1),\\
	\dz (c^2_s \rho_1) = 0& \text{in} ~ \Omega_h \times (0,1),
\end{cases}
\end{equation*}
complemented with stress-free and  impermeability physical boundary conditions:
\begin{equation*}\tag{BC-PE} \label{bc-pe}
	\dz v_p\big|_{z=0,1} = 0,~ w_p\big|_{z=0,1} = 0.
\end{equation*}

We observe again that if we consider system \eqref{PE0} subject to periodic boundary conditions, with fundamental periodic domain $\Omega_h \times 2\mathbb T$, and  subject to the following symmetry
\begin{equation*}\tag{SYM-PE}\label{SYM-PE}
	 \text{$ v_p $ and $ w_p $ are even and odd respectively in the $z$ variable,}
\end{equation*}
then its solutions obey the above symmetry. Moreover, when such solutions are restricted to the physical domain $\Omega_h \times (0,1)$ they automatically satisfy the physical boundary conditions \eqref{bc-pe}.

We aim at  investigating the asymptotic behavior, as $ \varepsilon \rightarrow 0^+ $, of the solutions to \eqref{CPE}. Owing to the symmetry property in \eqref{SYM-CPE} and the boundary conditions in \eqref{bc-cpe},
it suffices to study the following system:
\begin{equation}\label{CPE'}
\begin{cases}
\dt \rho^\varepsilon + \dvh (\rho^\varepsilon v^\varepsilon) + \dz (\rho^\varepsilon w^\varepsilon) = 0 & \text{in} ~ \Omega_h \times 2\mathbb T, \\
\dt (\rho^\varepsilon v^\varepsilon) + \dvh( \rho^\varepsilon v^\varepsilon\otimes v^\varepsilon) + \partial_z (\rho^\varepsilon w^\varepsilon v^\varepsilon) \\
~~~~ ~~ + \dfrac{1}{\varepsilon^2}\nablah P(\rho^\varepsilon)
 = \mu \deltah v^\varepsilon + \lambda \nablah \dvh v^\varepsilon + \partial_{zz} v^\varepsilon & \text{in} ~ \Omega_h \times 2\mathbb T, \\
\partial_z \rho^\varepsilon = 0 & \text{in} ~ \Omega_h \times 2\mathbb T,
\end{cases}
\end{equation}
subject to the periodic boundary condition and symmetry \eqref{SYM-CPE}.
We will show that the asymptotic system of \eqref{CPE'}, as $ \varepsilon \rightarrow 0^+ $, is the incompressible primitive equations, subject to the periodic boundary condition and symmetry \eqref{SYM-PE}:
\begin{equation}\label{PE}
\begin{cases}
	\dvh v_p + \dz w_p = 0 & \text{in} ~ \Omega_h \times 2\mathbb T, \\
	\rho_0 (\dt v_p + v_p \cdot \nablah v_p + w_p \dz v_p ) + \nablah (c^2_s \rho_1) \\
	~~~~ ~~~~ = \mu \deltah v_p + \lambda \nablah \dvh v_p +  \partial_{zz}v_p & \text{in} ~ \Omega_h \times 2\mathbb T,\\
	\dz (c^2_s \rho_1) = 0& \text{in} ~ \Omega_h \times 2\mathbb T,
\end{cases}
\end{equation}
with $ c^2_s = \gamma \rho_0^{\gamma-1} $ and $ \rho_0 = \text{constant} $. Here $ \rho_1 $ is the Lagrangian multiplier for the constraint \subeqref{PE}{1}  satisfying \begin{equation}\label{zero-average-pressure} \int_{\Omega_h \times 2\mathbb T} \rho_1 \idx  = 0. \end{equation}
In addition,
due to the conservation of 
linear momentum of
\eqref{PE}, we can impose the following condition for $
v_p $:
\begin{equation}\label{conservation}
	 \int_{\Omega_h\times 2\mathbb T} v_p \,dxdydz  = 0,
\end{equation}
for any time $ t \geq 0 $ as long as the solution exists.

We remark again that
the restrictions of the solutions of \eqref{CPE'} and \eqref{PE}
to the physical domain $ \Omega_h \times (0,1) $ solve the original system \eqref{CPE}, with physical boundary conditions \eqref{bc-cpe}, and system \eqref{PE0}, with the physical boundary conditions \eqref{bc-pe}, respectively, due to symmetry \eqref{SYM-CPE} and \eqref{SYM-PE}, provided the solutions exist and are regular enough.

\begin{figure}[h]
\centering
\begin {tikzpicture}[-latex, auto, node distance =3 cm and 6.5 cm, on grid,
state/.style ={rectangle,
draw, black , text=black, text width = 2.5 cm, align = center}]
\node[state] (A)  {
Compressible hydrodynamic equations
};
\node[state] (B) [below =of A] {{CPE}};
\node[state] (C) [right = of A] {{BE}};
\node[state] (D) [below =of C] {{PE}};
\path (A) edge [right] node[above] {Mach number $ \rightarrow 0^+ $} (C);
\path (A) edge [below]  node[left, text width = 2 cm, align = center] {aspect ratio $\rightarrow 0^+ $} (B);
\path (B) edge [right] node[below] {Mach number $ \rightarrow 0^+ $} (D);
\path (C) edge [below]  node[right, text width = 2 cm, align = center] {aspect ratio $\rightarrow 0^+ $} (D);
\path (A) edge [bend right = 25] node[above] {SARLM} (D);
\path (A) edge [bend right = -25] node[below] {LMSAR} (D);
\end{tikzpicture}
\caption{The PE diagram}
\label{pe-diagram}
\end{figure}
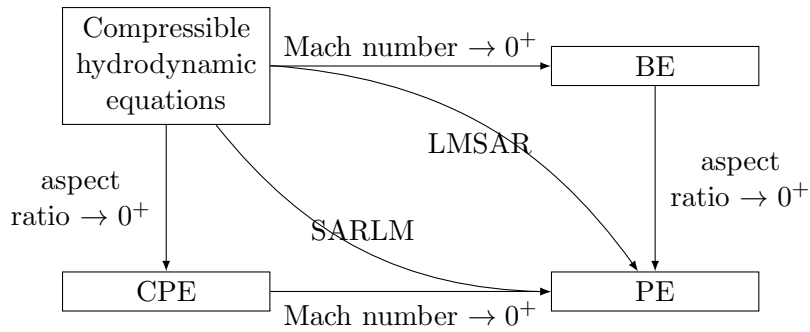

Historically, the limit system \eqref{PE}, besides being as the representation of the CPE in the $ p $-coordinates, is introduced as the limit system of Boussinesq equations (referred to as BE hereafter) when the aspect ratio between the vertical scale and the horizontal scale is very small, while the Boussinesq equations are the limit equations of the full compressible hydrodynamic equations with small Mach number and low stratification (see, e.g., \cite{JLLions1992}). That is to say, starting from the compressible hydrodynamic equations, by taking the low Mach number limit and then the small aspect ratio limit (referred to as LMSAR), one will arrive, formally, at the BE and then at the PE. On the other hand, by taking the small aspect ratio limit and then the low Mach number limit (referred to as SARLM), at least formally, the limit system of the compressible hydrodynamic equations is also the PE with the CPE as a middle state. Depending on the order of asymptotic limits, this gives us two directions from the hydrodynamic equations to the PE, which we will refer to as the PE diagram (see Figure \ref{pe-diagram}, below). The LMSAR part of the PE diagram has been shown to hold on solid ground in various settings (see, e.g., \cite{Rajagopal1996,Efeireisl2012,Novotny2011,wk,Azerad2001,Li2017}).
However, the validity of the SARLM part is relatively open. In order to fully justify the PE diagram, we investigate the low Mach number limit of the CPE in this work, which, as mentioned above, leads to the PE as the limit system. We remark that, the stratification effect of the gravity has been neglected in this work.

Each of the equations in the PE diagram have its own significances and have been studied separately in a large number of literatures. It will be certainly too ambitious to review all of those works. We refer  readers to the study of compressible hydrodynamic equations in, e.g., the books \cite{Lions1996,Lions1998,Majda2002vorticity,Feireisl2004}. As the limit system of the PE diagram, the primitive equations (PE) have been investigated intensively since they are introduced in \cite{Lions1992,JLLions1992,Lions1993,JLLions1994}. For instance, the global weak solutions are established in \cite{JLLions1992}. Local well-posedness with general data and global well-posedness with small data of strong solutions to the PE  in three spatial dimensions have been studied in
  \cite{GuillenGonzalez2001}
by Guill\'{e}n-Gonz\'{a}lez, Masmoudi and Rodr\'{i}guez-Bellido. Petcu and Wirosoetisno in \cite{Petcu2005} investigate the Sobolev and Gevrey regularity of the solutions to the PE.
  In \cite{HuTemamZiane2003}, in a domain with small depth, the authors address the global existence of strong solutions to PE. The well-posedness of unique global strong solutions is obtained by Cao and Titi in \cite{Cao2007} (see, also, \cite{Cao2003,Kobelkov2006,Kukavica2007a,Kukavica2007,Kukavica2014,hittmeir2017,Zelati2015,Hieber2016,Li2017a,Ignatova2012,Cao2012,Cao2014b,Cao2014,Cao2017,Cao2016,Cao2016a,Gerard-Varet2018} and the references therein for related studies). Considering the inviscid primitive equations, or hydrostatic incompressible Euler equations, the existence of solutions in the analytic function space and in the $H^s $ space are established in \cite{Brenier1999,Wong2012,Kukavica2011}.
  Renardy in \cite{Renardy2009} shows that the linearization of the equations at certain shear flows is ill-posed in the sense of Hadamard. Recently, the authors in \cite{Wong2014,Cao2015} construct finite-time blowup for the inviscid PE in the absence of rotation (i.e., without the Coriolis force).

In this work, we show that the PE can be viewed as the limit system of the CPE with the zero Mach number limit. The zero Mach number limit of the compressible hydrodynamic equations is a vast subject which has been studied for decades. Fruitful results have been obtained since the early works of Klainerman and Majda in \cite{Klainerman1981,Klainerman1982}, where the authors investigate the vanishing Mach number limit of compressible Euler equations with well-prepared initial data (see also \cite{SchochetCMP1986,Schochet1988}). Later by Ukai \cite{Ukai1986}, the theory of low Mach number limit of compressible Euler equations is extended to ill-prepared initial data (or called general data in some literatures). We remark here that the difference between the well-prepared and ill-prepared initial data, as in \cite{Metivier2001}, is that the well-prepared initial data have excluded the acoustic waves, while the ill-prepared initial data allow the interaction of the solutions with the high-frequency acoustic waves. In $ \mathbb R^n $, $ n = 2,3 $, such high-frequency acoustic waves disperse as shown in \cite{Ukai1986}, which implies strong convergences as the Mach number goes to zero. This can be also proved by applying the Strichartz's estimate (see, e.g., \cite{Ginibre1995,Lindblad1995,Keel1998}) for linear wave equations to the 
acoustic equations (see, e.g., \cite{DesjardinsGrenier1999}). In $ \mathbb T^n $, $ n = 2, 3 $, the high-frequency acoustic waves interact with each other and lead to fast oscillations and weak convergences when taking the low Mach number limit.
Such a fast oscillation phenomenon was first systematically studied in \cite{Schochet1994,Gallagher1998} for hyperbolic and parabolic systems, and by Lions and Masmoudi for  compressible Navier-Stokes equations in \cite{Lions1998a}. We refer, for the comparison of the whole space case, i.e., in $ \mathbb R^n $ and the periodic domain case, i.e., in $ \mathbb T^n $, to \cite{Masmoudi2001}. See also, \cite{Danchin2002,Danchin2002per,Danchin2005} for the studies in the Besov spaces. We acknowledge that the discussion here barely unveils the theory of the low Mach number limit, and we refer the reader to \cite{Alazard2006,Alazard2005,AlazardReview,Jiang2011,Metivier2001,FeireislSingularLimits,Rajagopal1996,Novotny2011,Masmoudi2007,Feireisl2008,Feireisl2011,Feireisl2015a,Efeireisl2012,wk} and the references therein for more comprehensive studies and recent progress.

In this work, we will focus on investigating the low Mach number limit of \eqref{CPE'} with well-prepared initial data. The convergence of the solutions of the CPE to the solution of the PE is in the strong sense. Furthermore, we are also able to obtain explicit convergence rate (see Theorem \ref{thm:low-mach}, below).
In particular, we obtain a class of global large solutions to \eqref{CPE'} with $ \varepsilon $ small enough.

We remark that in an upcoming paper \cite{CPE2PE2}, we will consider the low Mach number limit of \eqref{CPE'} with ill-prepared initial data, i.e., initial data with large, high-frequency acoustic waves.

\subsection{The low Mach number limit problem and main theorem}

In order to describe the aforementioned asymptotic limit, we study \eqref{CPE'} with $ (\rho^\varepsilon, v^\varepsilon, w^\varepsilon) $ close to an asymptotic state $ (\rho_0, v_p, w_p) $, where $ (\rho_0, v_p, w_p)  $ satisfies \eqref{PE}. For any $ \varepsilon > 0 $, the following ansatz is imposed:
\begin{equation}\label{ansatz}
\begin{cases}
	& \rho^\varepsilon := \rho_0 + \varepsilon^2 \rho_1 + \xi^\varepsilon ,\\ 
	& v^\varepsilon : = v_p + \psi^{\varepsilon,h}, \\
	& w^\varepsilon := w_p + \psi^{\varepsilon,z}.
\end{cases}
\end{equation}
Recall that $ \rho_1 $ has zero average in the domain (see \eqref{zero-average-pressure}).
This is motived by \cite{Hoff1998}. In addition, we shall employ the notation
\begin{equation*}
	\zeta^\varepsilon : = \varepsilon^2 \rho_1 + \xi^\varepsilon = \rho^\varepsilon - \rho_0.
\end{equation*}

For the sake of convenience, from time to time hereafter, we may drop the superscript $ \varepsilon $ from the functions. That is, we employ the notations $ \rho = \rho^\varepsilon, v = v^\varepsilon, w = w^\varepsilon, \xi = \xi^\varepsilon, \psi^h = \psi^{\varepsilon,h}, \psi^{z} = \psi^{\varepsilon,z}, \zeta = \zeta^\varepsilon $ throughout this work whenever there is no confusion.
Consequently, from \eqref{CPE'} and \eqref{PE}, the new unknown $(\xi, \psi^h, \psi^z) $ is governed by the following system:
\begin{equation}\label{eq:ori-perterb}
\begin{cases}
	\dt \xi + \rho_0(\dvh \psi^h + \dz \psi^z ) = - (\dvh (\xi v) + \dz (\xi w) ) \\
	~~ ~~ - \varepsilon^2 (\dt \rho_1 + \dvh (\rho_1 v) + \dz (\rho_1 w)) & \text{in} ~ \Omega_h \times 2\mathbb T,\\
	\rho \dt \psi^h + \rho v\cdot \nablah \psi^h + \rho w \dz \psi^h + \nablah (\varepsilon^{-2} (\rho^\gamma - \rho_0^\gamma - \varepsilon^2 c^2_s \rho_1)) \\
	~~ = \mu \deltah \psi^h + \lambda \nablah \dvh \psi^h + \partial_{zz} \psi^h + \rho_0^{-1}(\varepsilon^2 \rho_1 + \xi) \\
	~~ ~~ \times (\nablah (c^2_s \rho_1)
	 - \mu \deltah v_p - \lambda \nablah \dvh v_p - \partial_{zz} v_p ) \\
	~~ ~~ - \rho \psi^h \cdot\nablah v_p - \rho \psi^z \dz v_p & \text{in} ~ \Omega_h \times 2\mathbb T,\\
	\dz \xi = 0 & \text{in} ~ \Omega_h \times 2\mathbb T.
\end{cases}
\end{equation}
Observe that owing to the symmetry in \eqref{SYM-CPE} and \eqref{SYM-PE}, the following conditions hold automatically
\begin{equation}\label{bc-perturbed-eq}
	(\dz v, \dz v_p, \dz \psi^h)\bigr|_{z\in \mathbb Z} = 0 , ~(w, w_p, \psi^z)\bigr|_{z\in \mathbb Z}=0,
\end{equation}
for smooth enough functions, i.e., whenever one can make sense of these traces for $ z \in \mathbb Z $.
Recalling that $ c_s^2 = \gamma \rho_0^{\gamma - 1} $, we note that
\begin{equation*}
\begin{aligned}
	& \rho^\gamma - \rho_0^\gamma - \varepsilon^2 c^2_s \rho_1 =  \gamma \rho_0^{\gamma-1} (\rho-\rho_0) + \gamma (\gamma-1) \int_{\rho_0}^{\rho}(\rho-y)y^{\gamma-2} \,dy - \varepsilon^2 c^2_s \rho_1 \\
	& ~~~~ = c^2_s \xi + \mathcal R,
\end{aligned}
\end{equation*}
where
\begin{equation}\label{residue-est}
	\mathcal{R} = \mathcal{R}(\zeta) := \gamma (\gamma-1) \int_{\rho_0}^{\rho}(\rho-y)y^{\gamma-2} \,dy \leq C \zeta^2 \leq C (\varepsilon^4 \rho_1^2 + \xi^2 ),
\end{equation}
with $ C = C(\norm{\rho^{\gamma-2}}{\Lnorm{\infty}},\norm{\rho_0^{\gamma-2}}{\Lnorm{\infty}}) $.
Therefore, by denoting
\begin{equation}\label{nonlinearities}
\begin{aligned}
	& Q_p : = \rho_0^{-1} (\nablah (c^2_s \rho_1) - \mu \deltah v_p - \lambda \nablah \dvh v_p - \partial_{zz} v_p), \\
	& \mathcal F_1 : = \zeta Q_p,\\
	& \mathcal F_2 : = - \rho \psi^h \cdot\nablah v_p - \rho \psi^z \dz v_p,\\
	& \mathcal G_1 : = - (\dvh (\xi v) + \dz (\xi w) ), \\
	& \mathcal G_2 := - \varepsilon^2 (\dt \rho_1 + \dvh (\rho_1 v) + \dz (\rho_1 w)),
\end{aligned}
\end{equation}
we can write system \eqref{eq:ori-perterb} as,
\begin{equation}\label{eq:perturbation}
\begin{cases}
	 \dt \xi + \rho_0(\dvh \psi^h + \dz \psi^z ) = \mathcal G_1 + \mathcal G_2& \text{in} ~ \Omega_h \times 2\mathbb T,\\
	 \rho \dt \psi^h + \rho v\cdot \nablah \psi^h + \rho w \dz \psi^h + \nablah (\varepsilon^{-2} c^2_s \xi)  = \mu \deltah \psi^h\\
	 ~~~~ + \lambda \nablah \dvh \psi^h + \partial_{zz} \psi^h + \mathcal F_1 + \mathcal F_2 - \nablah (\varepsilon^{-2} \mathcal R) & \text{in} ~ \Omega_h \times 2\mathbb T, \\
	 \dz \xi = 0 & \text{in} ~ \Omega_h \times 2\mathbb T.
\end{cases}
\end{equation}

In order to recover the vertical velocity perturbation $ \psi^z $ from \eqref{CPE'}, we introduce the following notations, for any function $ f $ in $ \Omega_h \times 2\mathbb T $, which is also even in the $ z $-variable:
\begin{equation*}
\overline f(x,y,t) := \int_0^1 f(x,y,z',t) \,dz' ~~ \text{and} ~~ \widetilde f := f - \overline f.
\end{equation*}
The periodicity and symmetry of $ f $ imply that $ \overline f(x,y,t) = \int_{k}^{k+1} f(x,y,z',t)\,dz' $ for any $ k \in \mathbb{Z} $.
Notice that from \subeqref{CPE'}{3}, $ \rho $ is independent of the vertical variable $ z $.
Then by averaging \subeqref{CPE'}{1} over the vertical direction, thanks to \eqref{bc-perturbed-eq}, one will get
\begin{gather*}
	\dt \rho + \dvh (\rho \overline v) = 0, ~~ \text{and after comparing with \subeqref{CPE'}{1}},\\
	\dvh( \rho \widetilde v) + \dz (\rho w) = 0.
\end{gather*}
In particular, from the above, the vertical velocity $ w $ is determined through $ \rho, v $ by the formula, thanks to \eqref{bc-perturbed-eq}:
\begin{align}
	& \rho w = - \int_0^z \dvh (\rho \widetilde v) \,dz' = - \int_{0}^z \bigl( \rho \dvh v_p + \rho \dvh \psi^h - \rho \dvh \overline{\psi^h} {\nonumber} \\
	& ~~~~ ~~~~ + v \cdot\nablah \rho - \overline v \cdot \nablah \rho \bigr) \,dz', ~~~~ \text{and therefore} \label{id:vertical_velocity} \\
	& \rho \psi^z = \rho w - \rho w_p = - \int_0^z \bigl( \rho \dvh \widetilde{\psi^h} + \widetilde{v} \cdot\nablah \rho \bigr)  \,dz'   \nonumber \\
	& ~~~~ ~~~~  = - \int_0^z \bigl( \dvh(\rho \widetilde{\psi^h}) + \widetilde{v_p} \cdot\nablah \rho \bigr) \,dz',\label{id:vertical_perturbation}
\end{align}
where we have substituted the following identity thanks to \subeqref{PE}{1} and \eqref{bc-perturbed-eq},
\begin{equation}\label{id:vertical-vel-PE}
	w_p = - \int_0^z \dvh v_p \,dz'.
\end{equation}
Such facts imply that in \eqref{CPE'}, \eqref{PE} and \eqref{eq:perturbation}, the vertical velocities and the vertical perturbation, i.e., $ w^\varepsilon, w_p^\varepsilon, \psi^{\varepsilon,z}  $, are fully determined by $ v^\varepsilon, v_p^\varepsilon, \rho^\varepsilon $. Therefore, there is no need to impose initial data for $ w^\varepsilon, w_p^\varepsilon, \psi^{\varepsilon,z} $. 

System \eqref{eq:perturbation} (or equivalently \eqref{eq:ori-perterb}) is complemented with initial data,
\begin{equation}\label{initial-data-well-prepared}
	\begin{gathered}
		(\xi, \psi^h)\bigr|_{t=0} = (\xi_{in}, \psi^h_{in}) \in H^2(\Omega_h \times 2\mathbb T ;\mathbb R) \times H^2(\Omega_h\times 2\mathbb T;\mathbb R^2), ~ \text{where}\\
		\dz \xi_{in} = 0, ~ \text{and $\psi_{in}^h $ is even in the $ z $-variable},
	\end{gathered}
\end{equation}
with the compatibility conditions:
\begin{equation}\label{cmpt-cds-perturbed}
	\begin{aligned}
	& \xi_{in,1} = - \rho_0(\dvh \psi^h_{in} + \dz \psi^z_{in} )  - (\dvh (\xi_{in} v_{in}^{\varepsilon}) + \dz (\xi_{in} w_{in}^{\varepsilon}) ) \\
	 & ~~~~ - \varepsilon^2 \bigl(\rho_{1,in,1} + \dvh (\rho_{1,in} v_{in}^\varepsilon) + \dz (\rho_{1,in} w_{in}^\varepsilon)\bigr) ~~~~ \text{in} ~ \Omega_h \times 2\mathbb T,\\
	& \rho_{in}^\varepsilon \psi^h_{in,1} + \rho_{in}^\varepsilon v_{in}^\varepsilon \cdot \nablah \psi^h_{in} + \rho_{in}^\varepsilon w_{in}^\varepsilon \dz \psi^h_{in} \\
	& ~~~~ + \nablah \bigl(\varepsilon^{-2} ((\rho_{in}^{\varepsilon})^\gamma - \rho_0^\gamma - \varepsilon^2 c^2_s \rho_{1,in}) \bigr)
	 = \mu \deltah \psi^h_{in} + \lambda \nablah \dvh \psi^h_{in} \\
	& ~~~~ + \partial_{zz} \psi^h_{in} + \rho_0^{-1}(\varepsilon^2 \rho_{1,in} + \xi_{in})
	 \bigl(\nablah (c^2_s \rho_{1,in}) - \mu \deltah v_{p,in} \\
	& ~~~~ - \lambda \nablah \dvh v_{p,in} - \partial_{zz} v_{p,in} \bigr)
	- \rho_{in}^\varepsilon \psi^h_{in} \cdot\nablah v_{p,in} - \rho_{in}^\varepsilon \psi^z_{in} \dz v_{p,in} \\
	& ~~~~ \text{in} ~ \Omega_h \times 2\mathbb T, ~~~~
	 \dz \xi_{in} = 0 ~\text{and}~ \dz \xi_{in,1} = 0 ~~~~ \text{in} ~ \Omega_h \times 2\mathbb T, \\
	&  \text{with} ~ (\xi_{in,1}, \psi^h_{in,1}) \in L^2(\Omega_h \times 2\mathbb T) \times L^2(\Omega_h \times 2\mathbb T),
	\end{aligned}
\end{equation}
where $ \rho_{in}^\varepsilon = \rho_0 + \varepsilon^2 \rho_{1,in} + \xi_{in}, ~ v_{in}^{\varepsilon} = v_{p,in} + \psi^h_{in}, ~ w^\varepsilon_{in} = w_{p,in} + \psi^z_{in} $, and $ \psi^z_{in} $ is given by
\begin{equation*}
	\rho_{in}^\varepsilon \psi^z_{in} = - \int_0^z \bigl( \dvh (\rho_{in}^\varepsilon \widetilde{\psi^h_{in}} ) + \widetilde{v}_{p,in} \cdot \nablah \rho_{in}^\varepsilon \bigr) \,dz'.
\end{equation*}
Here $ v_{p,in}, \rho_{1,in}, \rho_{1,in,1}, w_{p,in} $ are initial values of $ v_{p}, \rho_1, \dt \rho_{1}, w_{p} $, respectively, while $ w_{p,in} $ is given by $ w_{p,in} = - \int_0^z \dvh v_{p,in} \,dz' $.

It is worth stressing
that we will chose the initial time derivatives of the perturbations, i.e., $ (\varepsilon^{-1}\xi_{in,1}, \psi^h_{in,1}) $ in \eqref{cmpt-cds-perturbed}, to be bounded, uniformly in $ \varepsilon $ (see \eqref{functional-initial-energy} and Theorem \ref{thm:low-mach}, below). The reason for such choices of initial data is to exclude the high-frequency acoustic waves which
corresponds to the fact that our initial data are well-prepared.

We denote the initial energy functional by
\begin{equation}\label{functional-initial-energy}
	\mathcal E_{in} := \norm{\psi^h_{in}}{H^2}^2
	+ \norm{\varepsilon\psi^h_{in,1}}{L^2}^2 + \norm{\varepsilon^{-1}\xi_{in}}{H^2}^2
	+ \norm{\xi_{in,1}}{L^2}^2.
\end{equation}

\begin{remark}
	$ (\rho_{in}^\varepsilon, v_{in}^\varepsilon ) $ is the corresponding initial datum of $ (\rho^\varepsilon, v^\varepsilon) $ for system \eqref{CPE'}. $ v_{p,in} $ is the initial datum of $ v_p $ for system \eqref{PE}. Accordingly, $ \rho_{1,in} = \rho_{1,in}(x,y), \rho_{1,in,1}= \rho_{1,in,1}(x,y) $ are determined by the following elliptic problems,
	\begin{align*}
	&  \begin{cases} - c^2_s \deltah \rho_{1,in} = \rho_0 \int_0^1 \dvh \blparenthese \dvh (v_{p,in}\otimes v_{p,in} ) \brparenthese \,dz ~~ \text{in} ~ \Omega_h,  \\
	\int_{\Omega_h} \rho_{1,in} \idxh = 0; \end{cases} \\
	&  \begin{cases} - c^2_s \deltah \rho_{1,in,1} = 2 \rho_0 \int_0^1 \dvh \blparenthese \dvh (v_{p,in}\otimes v_{p,in,1} ) \brparenthese \,dz \\ ~~ \text{in} ~ \Omega_h,  ~~
	\int_{\Omega_h} \rho_{1,in,1} \idxh = 0, \end{cases}
	\end{align*}
	where $ v_{p,in,1} $ is the initial value of $ \dt v_{p} $ determined by
	\begin{align*}
		& \rho_0  v_{p,in,1} = - \rho_0 (v_{p,in} \cdot \nablah v_{p,in} + w_{p,in} \dz v_{p,in} ) - \nablah (c^2_s \rho_{1,in}) \\
		& ~~~~  + \mu \deltah v_{p,in} + \lambda \nablah \dvh v_{p,in} +  \partial_{zz}v_{p,in} ~~~~ \text{in} ~ \Omega_h \times 2\mathbb T.
	\end{align*}
\end{remark}

As we stated before, we focus
in this work on the asymptotic limit as $ \varepsilon \rightarrow 0^+ $.
%
We have the following global regularity of the limit system \eqref{PE}:
\begin{theorem}[Global regularity of the PE]\label{thm:global-pe}
For $ \lambda < 4\mu < 12 \lambda  $, suppose that \eqref{PE} is complemented with initial data $ v_{p,in} \in H^1(\Omega_h \times 2\mathbb T) $, which is even in the $ z $-variable, and satisfies the compatibility conditions:
\begin{equation}\label{compatible-condition-pe}
	\int_{\Omega_h \times 2\mathbb T}  v_{p,in} \idx = 0, ~ \dvh  \overline{v}_{p,in} = 0. \end{equation}
Then there exists a solution $ (v_p, \rho_1) $, with $ \int_0^1 \dvh v_{p} \,dz = 0 $ and $ w_p $ given by \eqref{id:vertical-vel-PE}, to the primitive equations \eqref{PE}. Moreover,
there is a constant $ C_{p,in} $ depending only $ \norm{v_{p,in}}{\Hnorm{1}} $ such that
\begin{equation*}
	\sup_{0\leq t< \infty} \norm{v_p(t)}{\Hnorm{1}}^2 + \int_0^\infty \biggl( \norm{\nabla v_p(t)}{\Hnorm{1}}^2 + \norm{\dt v_p(t)}{\Lnorm 2 }^2 \biggr) \,dt  \leq C_{p,in}.
\end{equation*}
 Moreover,
\begin{equation*}
	\norm{v_p(t)}{\Lnorm2}^2 \leq C e^{-ct} \norm{v_{p,in} }{\Lnorm2}^2,
\end{equation*}
for some positive constants $ c, C $.
In addition, if $ v_{p,in} \in H^s(\Omega_h \times 2\mathbb T) $ for any integer $ s \geq 2 $, there is a positive constant $ C_{p,in,s} $, depending only on $ \norm{v_{p,in}}{\Hnorm{s}} $  such that
\begin{equation*}
\begin{aligned}
	& \sup_{0\leq t< \infty} \blparenthese \norm{v_p(t)}{H^s}^2 + \norm{\dt v_{p}(t)}{\Hnorm{s-2}}^2\brparenthese + \int_0^\infty \biggl( \norm{v_p(t)}{H^{s+1}}^2  \\
	& ~~~~ + \norm{\dt v_p(t)}{H^{s-1}}^2 \biggr) \,dt \leq C_{p,in,s}.
\end{aligned}
\end{equation*}	
\end{theorem}

\begin{proof}
	The local well-posedness of solutions to \eqref{PE} in the function space $ H^s $ has been established in \cite{Petcu2005}. What is left is the global regularity estimate, which is a direct consequence of Proposition \ref{prop:PE-Hs-data}, below.
\end{proof}

\begin{remark}
	We only have to be careful about the different estimates caused by the viscosity tensor. In particular, it is the $ L^q $ estimate of $ v_p $, below in section \ref{sec:global-pe}, that requires the constraints on the viscosity coefficients. For solutions with $ H^2 $ initial data, the result can be found in \cite{Li2017}. The new thing we treat here is the case when  $ v_{p,in} \in H^s $, for $ s > 2 $.
\end{remark}

Now we can state our main theorem in this work.
\begin{theorem}[Low Mach number limit of the CPE]\label{thm:low-mach} For $ \lambda < 4 \mu < 12 \lambda $, suppose $ v_{p,in} \in H^s(\Omega_h \times 2\mathbb T) $, with integer $ s \geq 3 $, and it satisfies the compatibility conditions \eqref{compatible-condition-pe}. Also, we complement \eqref{eq:perturbation} with initial data
$ (\xi^{\varepsilon}, \psi^{\varepsilon,h})\bigr|_{t=0} = (\xi_{in}, \psi^h_{in}) \in H^2(\Omega_h\times 2 \mathbb T) \times H^2(\Omega_h \times 2\mathbb T) $ as in \eqref{initial-data-well-prepared}, which satisfies the compatibility conditions \eqref{cmpt-cds-perturbed}. Recall that we also require $ v_{p,in}, \psi^h_{in} $ to be even in the $ z $-variable. Then there exists a positive constant $ \varepsilon_0 \in (0,1) $ small enough, such that if $ \varepsilon \in (0,\varepsilon_0) $ and
$ \mathcal E_{in} \leq \varepsilon^2 $, there exists a global unique strong solution $ (\xi^\varepsilon, \psi^{\varepsilon,h}) $ to system \eqref{eq:perturbation}. In particular, the following regularity is satisfied:
\begin{gather*}
	\xi^\varepsilon \in L^\infty(0,\infty;H^2(\Omega_h \times 2\mathbb T)), \\
	\dt \xi^\varepsilon \in L^\infty(0,\infty;L^2(\Omega_h \times 2\mathbb T)) \cap L^2(0,\infty;L^2(\Omega_h \times 2\mathbb T)), \\
	\nablah \xi^\varepsilon \in L^2(0,\infty;H^1(\Omega_h \times 2\mathbb T)), ~ \psi^{\varepsilon,h} \in L^\infty(0,\infty;H^2(\Omega_h \times 2\mathbb T)), \\
	\dt \psi^{\varepsilon,h} \in L^\infty(0,\infty;L^2(\Omega_h \times 2\mathbb T))\cap L^2(0,\infty;H^1(\Omega_h \times 2\mathbb T)), \\
	\nabla \psi^{\varepsilon,h} \in L^2(0,\infty;H^2(\Omega_h \times 2\mathbb T)).
\end{gather*}
In addition, we have the following estimate:
\begin{equation}\label{global-stability-thm}
		\begin{aligned}
			& \sup_{0 \leq t < \infty} \biggl\lbrace  \norm{\psi^{\varepsilon,h}(t)}{H^2}^2
				+ \norm{\varepsilon \dt \psi^{\varepsilon,h}(t)}{L^2}^2 + \norm{\varepsilon^{-1}\xi^\varepsilon(t)}{H^2}^2 \\
			& ~~~~
				+ \norm{\dt \xi^\varepsilon(t)}{L^2}^2  \biggr\rbrace
			+ \int_0^{\infty} \biggl\lbrace \norm{\nabla \psi^{\varepsilon,h}(t)}{\Hnorm{2}}^2 + \norm{\varepsilon \dt \psi^{\varepsilon,h}(t)}{\Hnorm{1}}^2\\
			& ~~~~ + \norm{\varepsilon^{-1} \nablah \xi^{\varepsilon}(t)}{\Hnorm{1}}^2
			 + \norm{\dt \xi^{\varepsilon}(t)}{\Lnorm{2}}^2 \biggr\rbrace \,dt  \leq C \varepsilon^2.
		\end{aligned}
	\end{equation}
	for some positive constant $ C $ depending only on $ \norm{v_{p,in}}{\Hnorm{3}} $, which is independent of $ \varepsilon $.
	In particular, $ (\rho^\varepsilon, v^\varepsilon, w^\varepsilon) $ as in \eqref{ansatz} is a globally defined strong solution to \eqref{CPE'} and the following asymptotic estimate holds:
	\begin{equation}\label{low-mach-limit-thm}
		\begin{aligned}
			& \sup_{0\leq t <\infty} \biggl\lbrace \norm{v^\varepsilon(t) - v_p(t)}{\Hnorm{2}}^2 + \norm{\rho^\varepsilon(t) - \rho_0}{\Hnorm{2}}^2 \\
			& ~~~~ + \norm{w^\varepsilon(t) - w_p(t)}{\Hnorm{1}}^2 \biggr\rbrace \leq C \varepsilon^2,
		\end{aligned}
	\end{equation}
	for some positive constant $ C $ depending only on $ \norm{v_{p,in}}{\Hnorm{3}} $, which is independent of $ \varepsilon $,
	where $ w^\varepsilon, w_p $ are given as in \eqref{id:vertical_velocity}, \eqref{id:vertical-vel-PE}, respectively.
\end{theorem}
\begin{remark}
	According to \eqref{global-stability-thm}, the time derivatives, in comparison to the spatial derivatives, have larger perturbations. However, thanks to the well-prepared data setting, they are bounded. 
\end{remark}
\begin{remark}
In addition to \eqref{conservation},
thanks to the conservation of  mass
and momentum in
 \eqref{CPE'}, we can impose the following conditions for $
\rho, v $:
\begin{equation}\label{conservation02}
	\int_{\Omega_h\times 2\mathbb T} \rho \,dxdydz = \rho_0\abs{\Omega_h\times 2\mathbb T}{}, ~
	\int_{\Omega_h\times 2\mathbb T} \rho v \,dxdydz  = 0.
\end{equation}
Then with a little more effort, under the same assumptions as in Theorem \ref{thm:low-mach}, one can further conclude that the perturbation energy $ \mathcal E (t) $ exponentially decays as time grows. This can be shown as follows. Consider any integer $ s \geq 1 $ and initial data $ v_{p,in} \in H^s(\Omega) $, satisfying the compatibility condition \eqref{compatible-condition-pe}, as in Theorem \ref{thm:global-pe}. Then the estimates in Theorem \ref{thm:global-pe} hold. Next, in the proof of Proposition \ref{prop:PE-Hs-data}, below, the differential inequalities (for example, \eqref{pe:Hs18Jun001} with $ s = s-1 $) yield,
	\begin{equation}\label{pe:Hs18Jun002}
		\dfrac{d}{dt} \norm{v_p}{\Hnorm{s}}^2 + c' \norm{v_p}{\Hnorm{s+1}}^2 \leq \mathfrak Y \norm{v_{p}}{\Hnorm{s}}^2,
	\end{equation}
	for any $ s \geq 1 $ and some positive constant $ c' $ with $ \int_0^\infty \mathfrak Y \,dt < \infty $. This implies,
	\begin{equation*}
		\dfrac{d}{dt} \biggl\lbrace e^{c't} \norm{v_p}{\Hnorm{s}}^2 \biggr\rbrace \leq \mathfrak Y e^{c't} \norm{v_{p}}{\Hnorm{s}}^2.
	\end{equation*}
	Thus after applying Gr\"onwall's inequality, one can obtain that
	\begin{equation*}
		\norm{v_p(t)}{\Hnorm{s}}^2 \lesssim e^{-c't}.
	\end{equation*}
	Furthermore, multiply \eqref{pe:Hs18Jun002} with $ e^{c''t} $, for some $ c'' \in (0,c') $, and integrate the resultant in the time variable; we deduce that, for arbitrary positive time $ T \in (0,\infty) $,
	\begin{align*}
		& \int_0^T e^{c''t} \norm{v_{p}(t)}{\Hnorm{s+1}}^2 \,dt \leq \int_0^T \mathfrak Y(t) e^{c''t} \norm{v_{p}(t)}{\Hnorm{s}}^2 \,dt + \norm{v_p(0)}{\Hnorm{s}}^2 \\
		& ~~~~ - e^{c''T} \norm{v_p(T)}{\Hnorm{s}}^2 + c'' \int_0^T e^{c''t} \norm{v_p(s)}{\Hnorm{s}}^2\,dt < \infty.
	\end{align*}
	Then, we have the following estimate of the $ H^s $ norm of $ v_p $:
	\begin{equation*}
		\sup_{0\leq t < \infty }e^{c_s t} \norm{v_{p}(t)}{\Hnorm{s}}^2 + \int_0^\infty e^{c_s t} \norm{v_p(t)}{\Hnorm{s+1}}^2 < \infty,
	\end{equation*}
	for $ s \geq 1 $ and some positive constant $ c_s $. With such estimates, similar arguments as in the proof of Proposition \ref{prop:PE-Hs-data} (see, for example, \eqref{pe:mi-006}), one can conclude that the time derivative of $ v_p $ will have similar estimates. That is:
	\begin{equation*}
	\begin{cases}
		\int_0^\infty e^{c_1 t} \norm{\dt v_p(t)}{\Lnorm{2}}^2 < \infty, & \text{for} ~ s = 1; \\
		\sup_{0\leq t < \infty }e^{c_s t} \norm{\dt v_{p}(t)}{\Hnorm{s-2}}^2 + \int_0^\infty e^{c_s t} \norm{\dt v_p(t)}{\Hnorm{s-1}}^2 < \infty, & \text{for} ~ s \geq 2.
	\end{cases}
	\end{equation*}
	Then, following the same lines as in Corollary \ref{cor:choosing-s}, below, will yield
	\begin{equation}\label{pe:Hs18Jun003}
		\int_0^\infty e^{c_pt} \mathfrak H_p(t) \,dt < \infty,
	\end{equation}
	for some positive constant $ c_p $,
	where $ \mathfrak H_p(t) $ is defined in
	\eqref{functional:dissipation-pe}, below.
	
	On the other hand, from the conservation of mass and momentum \eqref{conservation02}, one can derive the Poincar\'e type inequalities,
	\begin{equation*}
		\norm{\psi^{\varepsilon,h}}{\Lnorm{2}} \lesssim \norm{\nabla \psi^{\varepsilon,h}}{\Lnorm{2}} + \widetilde{Q}_2, 
		~~ \norm{\xi^\varepsilon}{\Lnorm{2}} \lesssim \norm{\nablah \xi^\varepsilon}{\Lnorm{2}},
	\end{equation*}
	where $ \widetilde{Q}_2 
	$ is at least quadratic of the perturbations. After
	applying these inequalities, the inequality in \eqref{estimates-total-02}, below, can be written as,
	\begin{equation*}
		\dfrac{d}{dt} \mathcal E_{LM} + C_1 \mathcal E_{LM} \lesssim \mathfrak H_p \varepsilon^2,
	\end{equation*}
	for some positive constant $ C_1 $, provided $ \mathcal E \lesssim \varepsilon^2 $ and $ \varepsilon $ small enough. Here, $ \mathcal E $ and $ \mathcal E_{LM} $ are defined in \eqref{functional:energy} and \eqref{functional:instant-energy}, respectively. Thus together with \eqref{pe:Hs18Jun003}, this inequality implies
	\begin{equation*}
		e^{C_2 t} \mathcal E_{LM}(t) < \infty,
	\end{equation*}
	for some positive constant $ C_2 $. Relation \eqref{equivalent-functionals} yields the decay of the perturbation energy $ \mathcal E $ as we claim. We leave the details to readers.
\end{remark}
	
%
%

This work will be organized as follows. In section \ref{sec:preliminaries}, we summarize the notations which will be commonly used in later paragraphs. Section \ref{sec:a-priori-est} focuses on the $ \varepsilon $-independent a priori estimates, which are the foundation of the low Mach number limit. In section \ref{sec:low-mach-number-limit}, we focus on the proof of Theorem \ref{thm:low-mach}. This will be shown through a continuity argument. Finally in section \ref{sec:global-pe}, we summarize the proof of Theorem \ref{thm:global-pe}.

\section{Preliminaries}\label{sec:preliminaries}
We use the notations
\begin{equation*}
	\int \cdot \idx = \int_{\Omega_h \times 2\mathbb T} \cdot \idx := \int_{\Omega_h \times 2\mathbb T} \cdot \,dxdydz, ~ \int_{\Omega_h} \cdot \idxh,
\end{equation*}
to represent the integrals in $ \Omega $ and $ \Omega_h $ respectively.
Hereafter, $ \partial_h \in \lbrace \partial_x, \partial_y \rbrace $ represents the horizontal derivatives, and $ \dz $ represents the vertical derivative.

	We will use $ \hnorm{\cdot}{}, \norm{\cdot}{} $ to denote norms in $ \Omega_h \subset \mathbb R^2 $ and $ \Omega_h \times 2\mathbb T \subset \mathbb R^3 $, respectively.
	After applying Ladyzhenskaya's and Agmon's inequalities
	in $ \Omega_h $ and $ \Omega $, directly we have
	\begin{equation}\label{ineq-supnorm}
	\begin{gathered}
		\hnorm{f}{\Lnorm{4}} \leq C \hnorm{f}{\Lnorm{2}}^{1/2} \hnorm{\nablah f}{\Lnorm{2}}^{1/2} + \hnorm{f}{\Lnorm{2}} , ~ 
		\hnorm{f}{\infty}
		\leq C \hnorm{f}{\Lnorm{2}}^{1/2} \hnorm{f}{\Hnorm{2}}^{1/2}, \\
		\norm{f}{\Lnorm{3}} \leq C \norm{f}{\Lnorm{2}}^{1/2} \norm{\nabla f}{\Lnorm{2}}^{1/2} + \norm{f}{\Lnorm{2}} ,
	\end{gathered}
	\end{equation}
	for the function $ f $ with bounded right-hand sides. Also, applying Minkowski's and H\"older's inequalities yields
	\begin{gather*}
	\hnorm{\overline{f}}{\Lnorm{q}} \leq \int_0^1 \hnorm{f(z)}{\Lnorm{q}} \,dz \leq C \norm{f}{\Lnorm{q}}, \\
	\text{and hence} ~ \norm{\widetilde{f}}{\Lnorm{q}} \leq C \norm{f}{\Lnorm{q}}, ~ q \in [1,\infty).
	\end{gather*}
	
We use $ \delta > 0 $ to denote a arbitrary constant which will be chosen later adaptively small. Correspondingly, $ C_\delta $ is some positive constant depending on $ \delta $.  In addition, for any quantities $ A $ and $ B $, $ A \lesssim B $ is used to denote that there exists a positive constant independent of the solutions such that $ A \leq C B $.

The following energy and dissipation functionals will be employed
\begin{align}
	& \mathcal E(t) := \norm{\psi^h(t)}{H^2}^2
	+ \norm{\varepsilon\psi^h_t(t)}{L^2}^2 + \norm{\varepsilon^{-1}\xi(t)}{H^2}^2
	+ \norm{\xi_t(t)}{L^2}^2, \label{functional:energy} \\
	& \mathcal D(t) := \norm{\nabla \psi^h(t)}{\Hnorm{2}}^2 + \norm{\varepsilon \psi^h_t(t)}{\Hnorm{1}}^2 + \norm{\varepsilon^{-1} \nablah \xi(t)}{\Hnorm{1}}^2 \nonumber \\
	& ~~~~ ~~~~ + \norm{\xi_t(t)}{\Lnorm{2}}^2.
	\label{functional:dissipation}
\end{align}
Then $ \mathcal E(0) = \mathcal E_{in} $, where $ \mathcal E_{in} $ is as in \eqref{functional-initial-energy}.
In this work,
we shall use $ Q(\mathcal E
) $ to denote a polynomial quantity, with positive coefficients, of $ \sqrt{\mathcal E}
$ and $ Q(0) = 0 $. In general, $ Q(\cdot) $ is a generic polynomial quantity, with positive coefficients, of the arguments and $ Q(0) = 0 $.

\section{$ \varepsilon $-independent \textbf{\textit{a priori}} estimate}\label{sec:a-priori-est}
This section is devoted to show the following:
\begin{proposition}\label{prop:a-prior-estimates}
	For any $ T > 0 $, $ t \in [0,T] $, suppose that the solution $ (v_p, \rho_1) $  with $ w_p $ given by \eqref{id:vertical-vel-PE}
	to \eqref{PE} satisfies
	\begin{equation}\label{boundness-and-decay}
	\begin{gathered}
		\norm{v_p(t)}{\Hnorm{3}}, \norm{v_{p,t}(t)}{\Lnorm{2}},
		\norm{\rho_1(t)}{\Hnorm{2}},  \norm{\rho_{1,t}(t) }{\Lnorm{2}},  \norm{w_p(t)}{\Hnorm{1}} \leq C,
		\\
		\int_0^T \blparenthese \norm{v_p(t)}{\Hnorm{3}}^2 + \norm{v_{p,t}(t)}{\Hnorm{1}}^2 + \norm{\rho_{1,tt}(t)}{\Lnorm{2}}^2 + \norm{\rho_{1,t}(t)}{\Hnorm{1}}^2 \\
		~~~~ + \norm{\rho_1(t)}{\Hnorm{2}}^2 + \norm{w_p(t)}{\Hnorm{2}}^2 + \norm{v_p(t)}{\Hnorm{2}} \brparenthese  \,dt \leq C,
	\end{gathered}
	\end{equation}
	for some positive constant $ C $, and
	\begin{equation}\label{boundness-density}
		\dfrac 1 2 \rho_0 < \rho < 2 \rho_0 ~~~~ \text{in} ~ (\Omega_h \times 2\mathbb T) \times [0,T).
	\end{equation}
	Then any solution $ (\psi^h, \xi) $ to \eqref{eq:perturbation}, with initial data as in \eqref{initial-data-well-prepared}, provided that it exists in the time interval $ \lbrack 0,T\rbrack $, with $ \psi^z $ given by \eqref{id:vertical_perturbation},  satisfies,
	\begin{equation}\label{prop:stability-CPE}
		\begin{aligned}
		& \sup_{0 \leq t \leq T} \mathcal E(t) + \int_0^T \mathcal D(t) \,dt  \leq C' e^{C' + Q(\sup_{0\leq t\leq T} \mathcal E(t))} \biggl\lbrace \varepsilon^2 + \mathcal E_{in} \\
		& ~~~~ + \bblparenthese \varepsilon^2 + (\varepsilon^2  + 1 ) Q(\sup_{0\leq t\leq T}  \mathcal E(t)) \bbrparenthese \int_0^T \mathcal D(t) \,dt \biggr\rbrace,
		\end{aligned}
	\end{equation}
	for some positive constant $ C' $ depending only on the bounds in \eqref{boundness-and-decay}. In particular, $ C'$ is  independent of $ \varepsilon $ and $ T $.
\end{proposition}
\begin{remark}
We remark here that, from the definition of  $ \mathcal E(t) $ in \eqref{functional:energy} and \eqref{ansatz},  \eqref{boundness-density} automatically holds for $ \varepsilon $ small enough if $ \sup_{0\leq t \leq T} \mathcal E(t) < \infty $ and \eqref{boundness-and-decay} holds.

Throughout the rest of this section, it is assumed that $ (\xi,\psi^h) $ with $ \psi^z $ given by \eqref{id:vertical_perturbation} is a solution to \eqref{eq:perturbation} which is smooth enough such that the estimates below can be established. To justify the arguments, one can employ the local well-posedness theory and the standard different quotient method to the corresponding lines below (replaced the differential operators by different quotients, for example). See, for instance, similar arguments in \cite{LT2018a,LT2018b}.
\end{remark}

We denote by $ \mathfrak G_p(t) $ a polynomial, with positive coefficients, quantity of the arguments $$ \norm{v_p(t)}{\Hnorm{3}}, \norm{v_{p,t}(t)}{\Lnorm{2}},
		\norm{\rho_1(t)}{\Hnorm{2}},  \norm{\rho_{1,t}(t) }{\Lnorm{2}},  \norm{w_p(t)}{\Hnorm{1}}, $$ and \begin{equation}\label{functional:dissipation-pe}\begin{gathered}
			\mathfrak H_p(t) : =  \norm{v_p(t)}{\Hnorm{3}}^2 + \norm{v_{p,t}(t)}{\Hnorm{1}}^2 + \norm{\rho_{1,tt}(t)}{\Lnorm{2}}^2 + \norm{\rho_{1,t}(t)}{\Hnorm{1}}^2 \\
		~~~~ + \norm{\rho_1(t)}{\Hnorm{2}}^2 + \norm{w_p(t)}{\Hnorm{2}}^2 + \norm{v_p(t)}{\Hnorm{2}}.
		\end{gathered}\end{equation}
In particular, \eqref{boundness-and-decay} of Proposition \ref{prop:a-prior-estimates}, is equivalent to $$ \sup_{0\leq t\leq T}\mathfrak G_p(t) + \int_0^T \mathfrak H_p(t)\,dt < C , $$
for some positive constant $ C $.
For the sake of convenience, we will shorten the notations $ \mathfrak G_p = \mathfrak G_p(t), \mathfrak H_p = \mathfrak H_p(t) $, below. We also remind the reader that we have assumed that \eqref{boundness-and-decay} and \eqref{boundness-density} hold throughout this section.

\subsection{Temporal derivatives}
We start by performing the time derivative estimate to the solutions to system \eqref{eq:perturbation}. Applying $ \dt $ to system \eqref{eq:perturbation} we will have the following system:
\begin{equation}\label{eq:perturbation-001}
	\begin{cases}
		\dt\xi_t + \rho_0 (\dvh \psi^h_t + \dz \psi^z_t) = \mathcal G_{1,t} + \mathcal G_{2,t} & \text{in} ~ \Omega_h \times 2\mathbb T, \\
		\rho \partial_{t} \psi^h_t + \rho v\cdot\nablah \psi^h_t + \rho w\dz  \psi^h_t + \nablah (\varepsilon^{-2} c^2_s \xi_t) = \mu  \deltah \psi^h_t \\
		~~~~ + \lambda \nablah \dvh \psi^h_t  + \partial_{zz} \psi^h_t + \mathcal F_{1,t} + \mathcal F_{2,t} - \nablah (\varepsilon^{-2} \mathcal R_t) + \mathcal H_t & \text{in} ~ \Omega_h \times 2\mathbb T,
	\end{cases}
\end{equation}
where
\begin{equation}\label{nonlinear-dt}
	\mathcal H_t := - \rho_t \psi^h_t - (\rho_t v + \rho v_t)\cdot\nablah \psi^h - ( \rho_t w + \rho w_t )\dz \psi^h.
\end{equation}
We will show the following:
\begin{lemma}\label{lm:temporal-derivative}
	In addition to the assumptions in Proposition \ref{prop:a-prior-estimates},
suppose that $ (\xi, \psi^h) $, with $ \psi^z $ given by \eqref{id:vertical_perturbation}, is a smooth solution to \eqref{eq:perturbation} in the time interval $ [0,T] $. We have
\begin{equation}\label{ene:temporal-derivative}
\begin{aligned}
	& \dfrac{d}{dt} \biggl\lbrace \dfrac 1 2 \norm{\rho^{1/2} \varepsilon \psi^h_t}{\Lnorm 2}^2 + \dfrac{c^2_s}{2\rho_0} \norm{\xi_t}{\Lnorm 2}^2 \biggr\rbrace + \mu \norm{\varepsilon \nablah \psi_t^h}{\Lnorm 2}^2 \\
	& ~~~~ + \lambda \norm{\varepsilon \dvh \psi^h_t}{\Lnorm{2}}^2
	 + \norm{\varepsilon \dz \psi^h_t}{\Lnorm{2}}^2   \leq \delta \bigl( \norm{\xi_t}{\Lnorm{2}}^2 + \norm{\varepsilon \nabla \psi^h_t}{\Lnorm{2}}^2 \\
	& ~~~~ + \norm{\varepsilon \psi^h_t}{\Lnorm{2}}^2 + \norm{\nabla \psi^h }{\Hnorm{2}}^2 \bigr)
	+ \varepsilon^2 C_\delta \bigl( Q(\mathcal E) + \mathfrak G_p \bigr)  \norm{\dz \psi^h}{\Hnorm{2}}^2 \\
	& ~~~~ + C_\delta Q(\mathcal E) \bigl( \norm{\varepsilon \psi^h_t}{\Lnorm{2}}^2 + \norm{\nabla \psi^h}{\Hnorm{1}}^2
	+ \norm{\xi_t}{\Lnorm{2}}^2 \bigr) \\
	& ~~~~
	+ C_\delta \bigl( Q(\mathcal E) + 1 + \mathfrak G_p \bigr) \mathfrak H_p \bigl( \norm{\psi^h}{\Hnorm{2}}^2 + \norm{\varepsilon \psi^h_t}{\Lnorm{2}}^2 + \varepsilon^2 \bigr) .
\end{aligned}
\end{equation}
\end{lemma}
\begin{proof}
Take the $ L^2 $-inner product of \subeqref{eq:perturbation-001}{2} with $ \varepsilon^2 \psi^h_t $. After applying integration by parts and substituting \subeqref{CPE'}{1}, we have the following:
\begin{equation*}
	\begin{aligned}
		& \dfrac{d}{dt} \biggl\lbrace \dfrac 1 2 \int \varepsilon^2 \rho \abs{\psi^h_t}{2} \idx \biggr\rbrace - \int c^2_s \xi_t \dvh \psi^h_t \idx + \int \biggl(  \mu \abs{\varepsilon\nablah\psi^h_t}{2} \\
		& ~~~~ ~~~~ + \lambda \abs{\varepsilon\dvh \psi^h_t}{2} + \abs{\varepsilon\dz \psi^h_t}{2} \biggr) \idx = \int \varepsilon^2 \mathcal F_{1,t} \cdot \psi^h_t \idx \\
		& ~~~~ + \int \varepsilon^2 \mathcal F_{2,t} \cdot \psi^h_t \idx - \int \nablah \mathcal R_t \cdot\psi^h_t \idx + \int \varepsilon^2 \mathcal H_t \cdot\psi^h_t \idx.
	\end{aligned}
\end{equation*}
We substitute \subeqref{eq:perturbation-001}{1} into the second term on the left-hand side of the above equation and compute the identity, thanks to \eqref{bc-perturbed-eq},
\begin{equation*}
	\begin{aligned}
		& - \int c^2_s \xi_t \dvh \psi^h_t \idx = \int  \rho_0^{-1} c^2_s \xi_t ( \dt \xi_t + \rho_0 \dz \psi^z_t - \mathcal G_{1,t} - \mathcal G_{2,t} ) \idx \\
		& ~~~~ = \dfrac{d}{dt} \biggl\lbrace \dfrac{c^2_s}{2\rho_0} \int \abs{\xi_t}{2} \idx \biggr\rbrace - \rho_0^{-1} c^2_s \int \xi_t \mathcal G_{1,t} \idx - \rho_0^{-1} c^2_s \int \xi_t \mathcal G_{2,t} \idx.
	\end{aligned}
\end{equation*}
Consequently, we have,
\begin{equation}\label{ee:002}
	\begin{aligned}
		& \dfrac{d}{dt} \biggl\lbrace \dfrac 1 2 \norm{\rho^{1/2} \varepsilon \psi^h_t}{\Lnorm 2}^2 + \dfrac{c^2_s}{2\rho_0} \norm{\xi_t}{\Lnorm 2}^2 \biggr\rbrace + \mu \norm{\varepsilon\nablah \psi_t^h}{\Lnorm2}^2 \\
		& ~~~~ + \lambda\norm{\varepsilon\dvh \psi^h_t}{\Lnorm 2}^2 + \norm{\varepsilon \dz \psi_t^h}{\Lnorm 2}^2  = \int \varepsilon^2 \mathcal F_{1,t} \cdot \psi^h_t \idx \\
		& ~~~~  + \int \varepsilon^2 \mathcal F_{2,t} \cdot \psi^h_t \idx + \rho_0^{-1} c^2_s \int \xi_t \mathcal G_{1,t} \idx + \rho_0^{-1} c^2_s \int \xi_t \mathcal G_{2,t} \idx \\
		& ~~~~  + \int \varepsilon^2 \mathcal H_t \cdot\psi^h_t \idx - \int \nablah \mathcal R_t \cdot\psi^h_t \idx =: \sum_{i=1}^{6} I_i. 
	\end{aligned}
\end{equation}
Next we estimate the right-hand side of \eqref{ee:002}.
After substituting \eqref{nonlinear-dt} into $ I_5 $, it can be written as
\begin{equation*}
	\begin{aligned}
		& I_{5} 
		= - \int \zeta_t \abs{\varepsilon \psi_t^h}{2}\idx - \int \bigr( \varepsilon^2 \zeta_t v \cdot\nablah \psi^h \cdot\psi^h_t + \varepsilon^2 \rho v_t \cdot \nablah \psi^h \cdot \psi^h_t  \bigr) \idx \\
		& ~~~~ - \int \bigl( \varepsilon^2 \zeta_t w \dz \psi^h \cdot\psi^h_t + \varepsilon^2 \rho w_t \dz \psi^h \cdot\psi^h_t  \bigr) \idx =: I_{5}' + I_{5}'' + I_{5}'''.
	\end{aligned}
\end{equation*}
Notice that $ \zeta = \varepsilon^2  \rho_1 + \xi $
is independent of the $ z $ variable. Then, for every $ \delta > 0 $ there exists a positive constant $ C_\delta $ such that
\begin{align*}
		& I_{5}' 
		\lesssim  \hnorm{\zeta_t}{\Lnorm{2}}  \int_0^1\hnorm{\abs{\varepsilon\psi_t^h}{2}}{\Lnorm{2}} \,dz = \norm{\zeta_t}{\Lnorm{2}} \int_0^1\hnorm{\varepsilon\psi_t^h}{\Lnorm{4}}^2 \,dz \\
		& ~~ \lesssim \norm{\zeta_t}{\Lnorm{2}} \int_0^1\biggl( \hnorm{\varepsilon\psi_t^h}{\Lnorm{2}}\hnorm{\varepsilon\nablah \psi_t^h}{\Lnorm{2}} + \hnorm{\varepsilon\psi_t^h}{\Lnorm{2}}^2 \biggr) \,dz\\
		& ~~ \lesssim \norm{\zeta_t}{\Lnorm{2}} (\norm{\varepsilon\psi_t^h}{\Lnorm{2}}\norm{\varepsilon\nablah \psi_t^h}{\Lnorm{2}} + \norm{\varepsilon\psi_t^h}{\Lnorm{2}}^2) \lesssim \delta \norm{\varepsilon\nablah \psi_t^h}{\Lnorm{2}}^2  \\
		& ~~~~ ~~~~+ \delta \norm{\varepsilon\psi^h_t}{\Lnorm{2}}^2 + C_\delta \bigl( \varepsilon^4 \norm{\rho_{1,t}}{\Lnorm{2}}^2 + \mathcal E \bigr) \norm{\varepsilon \psi^h_t}{\Lnorm{2}}^2,
\end{align*}
where we have applied the Minkowski, H\"older's, the Sobolev embedding and Young's inequalities.
On the other hand,
$ I_{5}'' $ can be estimated directly using H\"older's, the Sobolev embedding and Young's inequalities:
\begin{align*}
		& I_{5}'' \lesssim \varepsilon \norm{\zeta_t}{\Lnorm{2}} \norm{v}{\Lnorm\infty}  \norm{\nablah \psi^h}{\Lnorm6} (  \norm{\varepsilon\psi^h_t}{\Lnorm2}^{1/2} \norm{\varepsilon\nabla \psi^h_t}{\Lnorm2}^{1/2} +  \norm{\varepsilon\psi^h_t}{\Lnorm2}) \\
		& ~~~~ + \norm{\rho}{\Lnorm\infty}( \varepsilon \norm{v_{p,t}}{\Lnorm2} + \norm{\varepsilon\psi^h_t}{\Lnorm2}) \norm{\nablah \psi^h}{\Lnorm6}\\
		& ~~~~ ~~~~ \times ( \norm{\varepsilon\psi^h_t}{\Lnorm2}^{1/2} \norm{\varepsilon\nabla \psi^h_t}{\Lnorm2}^{1/2} +  \norm{\varepsilon\psi^h_t}{\Lnorm2}) \\
		& ~~ \lesssim \delta \norm{\varepsilon \nabla \psi^h_t}{\Lnorm{2}}^2 + \delta \norm{\varepsilon \psi^h_t}{\Lnorm{2}}^2 + C_\delta  Q(\mathcal E)  
		\bigl( \norm{\nablah \psi^h}{\Hnorm{1}}^2 + \norm{\varepsilon \psi^h_t}{\Lnorm{2}}^2 \bigr)\\
		& ~~~~ + C_\delta \bigl( Q(\mathcal E) + 1 + \mathfrak G_p \bigr) \mathfrak H_p \bigl( \norm{\psi^h}{\Hnorm{2}}^2 + \norm{\varepsilon \psi^h_t}{\Lnorm{2}}^2 \bigr).
\end{align*}
\noindent
On the other hand, from \eqref{id:vertical_velocity} (or \eqref{id:vertical_perturbation}), we have the identities:
\begin{equation}\label{id:vertical_velocity-001}
	\begin{aligned}
		w & = w_p  - \int_0^z \biggl( \dvh \widetilde{\psi^h} + \widetilde v \cdot\nablah \log \rho \biggr) \,dz' = w_p + \psi^z ,\\
		w_t & = w_{p,t} - \int_0^z  \biggl( \dvh \widetilde{\psi^h_t} + \widetilde v_t \cdot\nablah \log \rho + \widetilde v \cdot\nablah (\log \rho)_t \biggr) \,dz' \\
		& = w_{p,t} + \psi^z_t.
	\end{aligned}
\end{equation}
Therefore, $ I_{5}''' $ can be written as
\begin{equation*}
	\begin{aligned}
		& I_{5}''' = - \int  \varepsilon^2 \bigl(\zeta_t w_p + \rho w_{p,t} \bigr)\dz \bigl( \psi^h \cdot\psi^h_t \bigr) \idx
		 + \int \bblbrack \varepsilon^2  \int_0^z \bigl( \dvh \widetilde{\psi^h} \\
		& ~~~~  + \widetilde v \cdot\nablah \log \rho \bigr) \,dz' \times \bigl( \zeta_t \dz \psi^h \cdot\psi^h_t \bigr) \bbrbrack \idx + \int \bblbrack \varepsilon^2  \int_0^z  \bigl( \rho \dvh \widetilde{\psi^h_t} \\
		& ~~~~ + \widetilde v_t \cdot\nablah \rho + \rho \widetilde v \cdot\nablah (\log \rho)_t \bigr) \,dz' \times \bigl( \dz \psi^h \cdot\psi^h_t \bigr)  \bbrbrack \idx =: \sum_{i=1}^3 I_{5,i}'''. 
	\end{aligned}
\end{equation*}
Then, we plug in identity \eqref{id:vertical-vel-PE} and apply the H\"older, Minkowski and Young inequalities to infer,
\begin{align*}
	& I_{5,1}''' = \int \bblbrack \varepsilon^2 \bigl(\zeta_t \int_0^z \dvh v_p \,dz' + \rho \int_0^z \dvh v_{p,t} \,dz' \bigr) \times \bigl( \dz \psi^h \cdot \psi_t^h \bigr) \bbrbrack \idx \\
	& ~~~~ =  \int \bblbrack \varepsilon^2 \zeta_t \blparenthese \int_0^z \dvh v_p \,dz' \brparenthese \times \bigl( \dz \psi^h \cdot \psi_t^h \bigr) \bbrbrack \idx \\
	& ~~~~ - \int \bblbrack \varepsilon^2 \blparenthese \int_0^z v_{p,t} \,dz' \brparenthese \cdot \nablah (\rho \dz \psi^h \cdot \psi^h_t) \bbrbrack \idx\\
	& ~~ \lesssim \varepsilon \int_0^1 \bblbrack \biggl( \int_0^z \hnorm{\nablah v_p}{\Lnorm{8}} \,dz' \biggr) \times \blparenthese \hnorm{\zeta_t}{\Lnorm{2}} \hnorm{\dz \psi^h}{\Lnorm{8}} \hnorm{\varepsilon \psi_t^h}{\Lnorm{4}} \brparenthese \bbrbrack \,dz \\
	& ~~~~ + \varepsilon \int_0^1 \bblbrack \biggl(\int_0^z \hnorm{v_{p,t}}{\Lnorm{2}} \,dz'\biggr) \times \blparenthese \hnorm{\nablah \zeta}{\Lnorm{8}} \hnorm{\dz \psi^h}{\Lnorm{8}} \hnorm{\varepsilon \psi_t^h}{\Lnorm{4}} \brparenthese \bbrbrack \,dz \\
	& ~~~~  + \varepsilon \int_0^1 \bblbrack \biggl( \int_0^z \hnorm{v_{p,t}}{\Lnorm{4}} \,dz' \biggr) \times \blparenthese \hnorm{\rho}{\Lnorm{\infty}}  \hnorm{\dz \nablah \psi^h}{\Lnorm{2}} \hnorm{\varepsilon \psi_t^h}{\Lnorm{4}} \\
	& ~~~~ ~~~~ + \hnorm{\rho}{\Lnorm{\infty}} \hnorm{\dz \psi^h}{\Lnorm{4}} \hnorm{\varepsilon \nablah \psi_t^h}{\Lnorm{2}} \bigr)  \bbrbrack \,dz
	 \lesssim \varepsilon \bigl( \norm{\zeta_t}{\Lnorm{2}} \norm{v_p}{\Hnorm{2}} \\
	& ~~~~ ~~~~ +  \norm{\nablah \zeta}{\Hnorm{1}} \norm{v_{p,t}}{\Lnorm{2}} + \norm{\rho}{\Lnorm{\infty}} \norm{v_{p,t}}{\Hnorm{1}} \bigr) \norm{\dz \psi^h}{\Hnorm{1}} \\
	& ~~~~ ~~~~ \times \bigl(\norm{\varepsilon \nabla \psi_t^h}{\Lnorm{2}} + \norm{\varepsilon \psi_t^h}{\Lnorm{2}}\bigr)  \lesssim \delta \norm{\varepsilon \nabla \psi^h_t}{\Lnorm{2}}^2 + \delta \norm{\varepsilon \psi^h_t}{\Lnorm{2}}^2 \\
	& ~~~~ + C_\delta \bigl( Q(\mathcal E) + 1 + \mathfrak G_p \bigr) \mathfrak H_p \norm{\psi^h}{\Hnorm{2}}^2 .
	\end{align*}
\noindent On the other hand, a straight forward estimate shows that
\begin{align*}
	& I_{5,2}''' \lesssim \varepsilon \int_0^1 \bblparenthese \hnorm{\nablah \psi^h}{\Lnorm8} + \hnorm{v}{\Lnorm\infty}\hnorm{\nablah \zeta}{\Lnorm8} \bbrparenthese \,dz \times  \int_0^1 \hnorm{\zeta_t}{\Lnorm2} \hnorm{\dz \psi^h}{\Lnorm8} \hnorm{\varepsilon\psi^h_t}{\Lnorm4}  \,dz\\
	& ~~~~ \lesssim \varepsilon \norm{\zeta_t}{\Lnorm2} \bigl(\norm{\nablah \psi^h}{H^1}+ \norm{v}{H^2} \norm{\nablah \zeta}{H^1}\bigr)\norm{\dz \psi^h}{H^1} \\
	& ~~~~ ~~~~ \times \bigl( \norm{\varepsilon\nablah\psi^h_t}{\Lnorm2} + \norm{\varepsilon\psi^h_t}{\Lnorm2} \bigr) \lesssim \delta \norm{\varepsilon \nabla \psi^h_t}{\Lnorm{2}}^2 + \delta \norm{\varepsilon \psi^h_t}{\Lnorm{2}}^2 \\
	& ~~~~ + C_\delta Q(\mathcal E)
	 \norm{\dz \psi^h}{\Hnorm{1}}^2 + C_\delta \bigl( Q(\mathcal E) + 1 +\mathfrak G_p \bigr) \mathfrak H_p \norm{\psi^h}{\Hnorm{2}}^2 .
\end{align*}
\noindent
To estimate $ I_{5,3}''' $, we apply integration by parts as below,
\begin{align*}
	& I_{5,3}''' = \int \bblbrack \varepsilon^2 \int_0^z \bigl(  \rho \dvh \widetilde{\psi^h_t} + \widetilde{v_t} \cdot\nablah \rho -
	 \zeta_t \widetilde v \cdot \nablah \log\rho - \zeta_t \dvh \widetilde v \bigr)
	 \,dz \\
	 & ~~~~ ~~~~ \times \bigl( \dz \psi^h \cdot\psi^h_t \bigr) \bbrbrack \idx
	  - \int \bblbrack \varepsilon^2 \int_0^z \zeta_t \widetilde v \,dz \cdot \nablah (\dz \psi^h \cdot\psi^h_t ) \bbrbrack \idx,
\end{align*}
from which, we infer
\begin{align*}
	& I_{5,3}''' \lesssim \int_0^1 \biggl( \hnorm{\rho}{\Lnorm\infty} \hnorm{\varepsilon\nablah \psi^h_t}{\Lnorm2} + \hnorm{\varepsilon v_t}{\Lnorm 4} \hnorm{\nablah \zeta}{\Lnorm 4} \biggr) \,dz \\
	& ~~~~ ~~~~ \times \int_0^1 \hnorm{\dz \psi^h}{\Lnorm4} \hnorm{\varepsilon\psi^h_t}{\Lnorm 4} \,dz
	 + \int_0^1 \hnorm{\dz \psi^h}{\Lnorm8} \hnorm{\varepsilon\psi^h_t}{\Lnorm4} \,dz \\
	& ~~~~ ~~~~ \times \int_0^1 \biggl( \varepsilon\hnorm{\zeta_t}{\Lnorm 2} \hnorm{v}{\infty} \hnorm{\nablah \zeta}{\Lnorm 8} +\varepsilon \hnorm{\zeta_t}{\Lnorm 2} \hnorm{\nablah v}{\Lnorm 8} \biggr) \,dz
	  \\
	 & ~~~~ + \int_0^1 \biggl( \hnorm{\dz \psi^h}{\Lnorm \infty} \hnorm{\varepsilon\nablah \psi^h_t}{\Lnorm 2} + \hnorm{\nablah \dz \psi^h }{\Lnorm 4} \hnorm{\varepsilon \psi_t^h}{\Lnorm 4} \biggr) \,dz \\
	 & ~~~~ ~~~~ \times \int_0^1 \varepsilon \hnorm{\zeta_t}{\Lnorm 2} \hnorm{v}{\Lnorm \infty} \,dz
	\lesssim \bigl( \norm{\rho}{\Lnorm\infty}\norm{\varepsilon\nablah \psi^h_t}{\Lnorm2} \\
	& ~~~~ ~~~~ + \norm{\varepsilon v_{p,t}}{\Hnorm{1}} \norm{\nablah \zeta}{\Hnorm{1}} \bigr) \norm{\psi^h_z}{H^1}
	\bigl(\norm{\varepsilon\psi^h_t}{\Lnorm 2}^{1/2} \norm{\varepsilon\nablah \psi^h_t}{\Lnorm 2}^{1/2} \\
	& ~~~~ ~~~~ + \norm{\varepsilon\psi^h_t}{\Lnorm 2} \bigr) + \norm{\nablah \zeta}{H^1}\norm{\psi^h_z}{H^1}
	 \bigl(\norm{\varepsilon\psi^h_t}{\Lnorm2}\norm{\varepsilon\nablah\psi^h_t}{\Lnorm2} \\
	& ~~~~ ~~~~+\norm{\varepsilon\psi^h_t}{\Lnorm2}^2\bigr) + \varepsilon\norm{v}{H^2}\bigl(1+\norm{\nablah \zeta}{H^1}\bigr)\norm{\psi^h_z}{H^2}\norm{\zeta_t}{\Lnorm2} \\
	& ~~~~ ~~~~ \times \bigl(\norm{\varepsilon\nablah \psi^h_t}{\Lnorm2} + \norm{\varepsilon\psi^h_t}{\Lnorm 2}\bigr)\\
	& ~~ \lesssim \delta \norm{\varepsilon\nablah \psi_t^h}{\Lnorm 2}^2 + \delta \norm{\varepsilon \psi^h_t}{\Lnorm{2}}^2 + \varepsilon^2 C_\delta \bigl( Q(\mathcal E) + \mathfrak G_p\bigr) \norm{\psi^h_z}{\Hnorm{2}}^2 \\
	& ~~~~ + C_\delta Q(\mathcal E) \bigl( \norm{\varepsilon\psi_t^h}{\Lnorm2}^2 + \norm{\dz \psi^h}{\Hnorm{1}}^2 \bigr)
	 + C_\delta \bigl(Q(\mathcal E) + 1 + \mathfrak G_p \bigr) \mathfrak H_p \\
	 & ~~~~ ~~~~ \times \bigl( \norm{\psi^h}{H^2}^2 + \norm{\varepsilon \psi^h_t}{\Lnorm{2}}^2 \bigr).
\end{align*}
Therefore, we have shown
\begin{equation}\label{estimate:001}
	\begin{aligned}
		& I_{5} \lesssim \delta \norm{\varepsilon \nabla \psi^h_t}{\Lnorm 2}^2 + \delta \norm{\varepsilon \psi^h_t}{\Lnorm{2}}^2 +  \varepsilon^2 C_\delta \bigr(Q(\mathcal E) +  \mathfrak G_p \bigl)  \norm{\psi^h_z}{H^2}^2 \\
		& ~~~~ + C_\delta Q(\mathcal E) \bigl( \norm{\varepsilon \psi_t^h}{\Lnorm 2}^2 + \norm{\nabla \psi^h}{H^1}^2 \bigr) \\
		& ~~~~ + C_\delta \bigl(Q(\mathcal E) + 1 + \mathfrak G_p \bigr) \mathfrak H_p \bigl( \norm{\psi^h}{H^2}^2 + \norm{\varepsilon \psi^h_t}{\Lnorm{2}}^2 \bigr).
	\end{aligned}
\end{equation}
Next, after substituting \eqref{nonlinearities} into $ I_2 $, it follows that
\begin{align*}
	& I_2 = - \int \varepsilon^2 \bigl( \rho\psi^h\cdot \nablah v_{p,t} + \zeta_t \psi^h \cdot\nablah v_p + \rho \psi^h_t \cdot \nablah v_p \bigr) \cdot \psi^h_t \idx \\
	& ~~~~ - \int \varepsilon^2 \bigl(\rho \psi^z \dz v_{p,t} + \zeta_t \psi^z \dz v_p + \rho \psi^z_t \dz v_p \bigr) \cdot\psi_t^h \idx =: I_2' + I_2''.
\end{align*}
Similarly as before,
\begin{align*}
	& I_2' \lesssim \norm{\rho}{\Lnorm\infty}\norm{\nablah v_p}{\Lnorm6} \norm{\varepsilon \psi_t^h}{\Lnorm3}  \norm{\varepsilon \psi^h_t}{\Lnorm2} + \varepsilon \norm{\rho}{\Lnorm{\infty}} \norm{\nablah v_{p,t}}{\Lnorm{2}} \\
	& ~~~~ ~~~~ \times \norm{\psi^h}{\Lnorm{6}} \norm{\psi^h_t}{\Lnorm{3}} + \varepsilon \norm{\nablah v_p}{\Lnorm{\infty}} \norm{\psi^h}{\Lnorm6} \norm{\zeta_t}{\Lnorm2} \norm{\varepsilon \psi^h_t}{\Lnorm{3}} \\
	& ~~ \lesssim \norm{\rho}{\Hnorm{2}}\norm{\nablah v_p}{\Lnorm{2}} \bigl( \norm{\varepsilon\psi^h_t}{\Lnorm{2}}^{3/2} \norm{\varepsilon \nabla \psi^h_t}{\Lnorm{2}}^{1/2} + \norm{\varepsilon\psi^h_t}{\Lnorm{2}}^2 \bigr) \\
	& ~~~~ +\varepsilon \norm{\psi^h}{\Hnorm{1}}\norm{\nablah v_{p,t}}{\Lnorm{2}}\bigl( \norm{\varepsilon\psi^h_t}{\Lnorm{2}}^{1/2} \norm{\varepsilon \nabla \psi^h_t}{\Lnorm{2}}^{1/2} + \norm{\varepsilon\psi^h_t}{\Lnorm{2}} \bigr)  \\
	& ~~~~ + \varepsilon \norm{\nablah v_p}{\Hnorm{2}} \norm{\psi^h}{\Hnorm{1}}\norm{\zeta_t}{\Lnorm{2}} \bigl( \norm{\varepsilon \psi^h_t}{\Lnorm{2}}^{1/2}\norm{\varepsilon \nabla \psi^h_t}{\Lnorm{2}}^{1/2} + \norm{\varepsilon \psi^h_t}{\Lnorm{2}} \bigr) \\
	& ~~ \lesssim \delta \norm{\varepsilon \nabla \psi^h_t}{\Lnorm{2}} + C_\delta Q(\mathcal E)   \norm{\varepsilon \psi^h_t}{\Lnorm{2}}^2 
	+ C_\delta \bigl( Q(\mathcal E) + 1 + \mathfrak G_p \bigr) \mathfrak H_p \norm{\varepsilon \psi^h_t}{\Lnorm{2}}^2 \\
	& ~~~~ + \varepsilon^2 C_\delta \mathfrak H_p.
\end{align*}
On the other hand, after substituting \eqref{id:vertical_velocity-001}, $ I_2'' $ can be written as
\begin{align*}
	& I_2'' = \int \bblbrack \varepsilon^2 \int_0^z  \biggl( \dvh \widetilde{\psi^h} + \widetilde{v} \cdot\nablah\log \rho \biggr) \,dz \times \bblparenthese  \bigl( \rho \dz v_{p,t} + \zeta_t \dz v_p \bigr) \cdot\psi_t^h \bbrparenthese \bbrbrack \idx \\
	& ~~~~ + \int \bblbrack \varepsilon^2 \int_0^z \biggl( \dvh \widetilde{\psi^h_t} + \widetilde{v_t} \cdot\nablah \log \rho + \widetilde v \cdot\nablah (\log \rho)_t \biggr) \,dz \\
	& ~~~~ ~~~~ \times \biggl( \rho  \dz v_p \cdot\psi_t^h \biggr) \bbrbrack \idx
	 =: I_{2,1}'' + I_{2,2}''.
\end{align*}
Then we have the following estimate:
\begin{align*}
	& I_{2,1}'' \lesssim \int_0^1 \biggl( \hnorm{\nablah \psi^h}{\Lnorm4} + \hnorm{v}{\Lnorm\infty} \hnorm{\nablah \zeta}{\Lnorm 4} \biggr) \,dz \times \int_0^1 \bblparenthese \bigl(\varepsilon \hnorm{\dz v_{p,t}}{\Lnorm2} \\
	& ~~~~ ~~~~ +\varepsilon \hnorm{\zeta_t}{\Lnorm 2}\hnorm{\dz v_p}{\Lnorm\infty} \bigr) \hnorm{\varepsilon \psi_t^h}{\Lnorm 4}  \bbrparenthese \,dz \lesssim \bigl( \norm{v}{H^2} \norm{\nablah \zeta}{H^1} \\
	& ~~~~ ~~~~ + \norm{\nablah \psi^h}{H^1} \bigr)
	 \bigl( \varepsilon \norm{\dz v_{p,t}}{\Lnorm{2}} + \varepsilon \norm{\zeta_t}{\Lnorm{2}} \norm{\dz v_p}{\Hnorm{2}} \bigr) \\
	& ~~~~ ~~~~ \times \bigl( \norm{\varepsilon \psi^h_t}{\Lnorm{2}}^{1/2}\norm{\varepsilon \nabla \psi^h_t}{\Lnorm{2}}^{1/2} + \norm{\varepsilon \psi^h_t}{\Lnorm{2}} \bigr)
	\lesssim \delta \norm{\varepsilon \nabla \psi^h_t}{\Lnorm{2}}^2 \\
	& ~~~~ ~~~~ + C_\delta  Q(\mathcal E) \norm{\varepsilon \psi^h_t}{\Lnorm{2}}^2 + C_\delta \bigl( Q(\mathcal E) + 1 + \mathfrak G_p \bigr) \mathfrak H_p \norm{\varepsilon \psi^h_t}{\Lnorm{2}}^2
	 +\varepsilon^2 C_\delta \mathfrak H_p.
\end{align*}
To estimate $ I_{2,2}'' $, we first apply integration by parts as below,
\begin{align*}
	& I_{2,2}'' = \int \bblbrack \varepsilon^2 \int_0^z \biggl(  \rho \dvh \widetilde{\psi^h_t} + \widetilde{v_t} \cdot\nablah \rho -
	 \zeta_t \widetilde v \cdot \nablah \log\rho - \zeta_t \dvh \widetilde v \biggr)
	 \,dz \\
	 & ~~~~ ~~~~ \times  \bigl( \dz v_p \cdot\psi^h_t \bigr) \bbrbrack \idx
	  - \int \bblbrack \varepsilon^2 \bblparenthese \int_0^z \zeta_t \widetilde v \,dz \bbrparenthese \cdot \nablah (\dz v_p \cdot\psi^h_t ) \bbrbrack \idx,
\end{align*}
which yields, similarly to the estimate of $ I_{5,3}''' $,
\begin{align*}
	& I_{2,2}'' \lesssim \int_0^1 \biggl( \hnorm{\rho}{\Lnorm\infty} \hnorm{\varepsilon\nablah \psi^h_t}{\Lnorm2} + \hnorm{\varepsilon v_t}{\Lnorm4}\hnorm{\nablah \zeta}{\Lnorm4} \biggr) \,dz \\
	& ~~~~ ~~~~ \times \int_0^1 \hnorm{\dz v_p}{\Lnorm4} \hnorm{\varepsilon\psi^h_t}{\Lnorm4} \,dz + \int_0^1 \hnorm{\dz v_p}{\Lnorm8} \hnorm{\varepsilon \psi^h_t}{\Lnorm4} \,dz  \\
	& ~~~~ ~~~~ \times \int_0^1 \biggl( \varepsilon \hnorm{\zeta_t}{\Lnorm2} \hnorm{v}{\Lnorm \infty}\hnorm{\nablah \zeta}{\Lnorm8} + \varepsilon \hnorm{\zeta_t}{\Lnorm2} \hnorm{\nablah v}{\Lnorm8} \biggr) \,dz  \\
	& ~~~~ + \int_0^1 \biggl( \hnorm{\dz v_p}{\Lnorm\infty} \hnorm{\varepsilon \nablah \psi^h_t}{\Lnorm2} +\hnorm{\nablah\dz v_p}{\Lnorm4} \hnorm{\varepsilon\psi^h_t}{\Lnorm4} \biggr) \,dz\\
	& ~~~~ ~~~~ \times  \int_0^1 \varepsilon\hnorm{\zeta_t}{\Lnorm2}\hnorm{v}{\Lnorm\infty} \,dz
	\lesssim \bigl( \norm{\rho}{\Lnorm\infty}\norm{\varepsilon\nablah \psi^h_t}{\Lnorm2} \\
	& ~~~~ ~~~~ + \norm{\varepsilon v_{p,t}}{\Hnorm{1}} \norm{\nablah \zeta}{\Hnorm{1}} \bigr) \norm{\dz v_p}{H^1}
	 \bigl(\norm{\varepsilon\psi^h_t}{\Lnorm 2}^{1/2} \norm{\varepsilon\nablah \psi^h_t}{\Lnorm 2}^{1/2} \\
	 & ~~~~ ~~~~ + \norm{\varepsilon\psi^h_t}{\Lnorm 2} \bigr) +
	  \bigl(\norm{\varepsilon\psi^h_t}{\Lnorm2}\norm{\varepsilon\nablah\psi^h_t}{\Lnorm2}+\norm{\varepsilon\psi^h_t}{\Lnorm2}^2\bigr) \\
	 & ~~~~ ~~~~ \times \norm{\nablah \zeta}{H^1}\norm{\dz v_p}{H^1} + \varepsilon\norm{v}{H^2}\bigl(1+\norm{\nablah \zeta}{H^1}\bigr) \\
	& ~~~~ ~~~~ \times \norm{\dz v_p}{H^2}\norm{\zeta_t}{\Lnorm2}\bigl(\norm{\varepsilon\nablah \psi^h_t}{\Lnorm2} + \norm{\varepsilon\psi^h_t}{\Lnorm 2}\bigr)\\
	& ~~ \lesssim \delta \norm{\varepsilon\nablah \psi_t^h}{\Lnorm 2}^2 + C_\delta Q(\mathcal E) \norm{\varepsilon\psi_t^h}{\Lnorm2}^2 + C_\delta \bigl( Q(\mathcal E) + 1 + \mathfrak G_p \bigr) \mathfrak H_p \\
	& ~~~~ ~~~~ \times \bigl( \norm{\varepsilon \psi^h_t}{\Lnorm{2}}^2 + \varepsilon^2 \bigr).
\end{align*}
Therefore,
\begin{equation}\label{estimate:002}
	\begin{aligned}
		& I_{2} \lesssim \delta \norm{\varepsilon \nabla \psi^h_t}{2}^2 + C_\delta Q(\mathcal E) \norm{\varepsilon \psi_t^h}{2}^2\\
		& ~~~~ + C_\delta \bigl( Q(\mathcal E) + 1 + \mathfrak G_p \bigr) \mathfrak H_p \bigl( \norm{\varepsilon \psi^h_t}{\Lnorm{2}}^2 + \varepsilon^2 \bigr).
	\end{aligned}
\end{equation}
Now, we will estimate $ I_3 $, which reads
\begin{equation}\label{estimate:003}
\begin{aligned}
	& I_3 = - \rho_0^{-1} c^2_s \int \biggl( \xi_t \dvh (\xi_t v) + \xi_t \dvh (\xi v_t)   +  \xi_t \dz(\xi_t w) \\
	& ~~~~ ~~~~ + \xi_t \dz(\xi w_t)  \biggr) \idx
	 = - \dfrac{c^2_s}{2\rho_0} \int \abs{\xi_t}{2}\dvh \psi^h \idx \\
	& ~~~~ - \rho_0^{-1} c^2_s \int \biggl( \xi_t \xi \dvh \psi^h_t + \xi_t v_t \cdot\nablah \xi  \biggr) \idx\\
	& ~~ \lesssim \norm{\xi_t}{\Lnorm2}( \norm{\xi_t}{\Lnorm2} \norm{\nablah \psi^h}{H^2}+ \norm{\varepsilon^{-1} \xi}{H^2} \norm{\varepsilon \nablah \psi^h_t}{\Lnorm2}) \\
	& ~~~~ + (\varepsilon \norm{v_{p,t}}{\Lnorm6}+ \norm{\varepsilon\nabla \psi^h_t}{\Lnorm2} + \norm{\varepsilon\psi^h_t}{\Lnorm2})\norm{\varepsilon^{-1} \nablah \xi}{H^1}\norm{\xi_t}{\Lnorm2} \\
	& ~~ \lesssim \delta \norm{\varepsilon \nabla \psi^h_t}{\Lnorm2}^2 + \delta \norm{\varepsilon \psi^h_t}{\Lnorm{2}}^2 + \delta \norm{\nablah \psi^h}{H^2}^2 + C_\delta Q(\mathcal E)\norm{\xi_t}{2}^2\\
	& ~~~~ + \varepsilon^2 C_\delta \mathfrak H_p.
\end{aligned}
\end{equation}
Here we have employed the facts that $ \xi $ is independent of the $z$-variable and that $ \int_0^1 \dvh v_p \,dz = \int_0^1 \dvh v_{p,t} \,dz  = 0 $.
The rest is straightforward. For instance, substituting \eqref{nonlinearities} in $ I_1 $ yields
\begin{equation}\label{estimate:004}
\begin{aligned}
	& I_1 = \rho_0^{-1} \int \varepsilon \zeta_t \bigl( \nablah (c^2_s \rho_1) - \mu \deltah v_p -\lambda \nablah \dvh v_p - \partial_{zz} v_p \bigr) \cdot \varepsilon \psi^h_t \idx\\
	& ~~~~ + \rho_0^{-1} \int \varepsilon \zeta \bigl( \nablah (c^2_s \rho_{1,t}) - \mu \deltah v_{p,t} -\lambda \nablah \dvh v_{p,t} - \partial_{zz} v_{p,t} \bigr) \cdot \varepsilon \psi^h_t \idx\\
	& ~~ = \rho_0^{-1} \int \varepsilon \zeta_t \bigl( \nablah (c^2_s \rho_1) - \mu \deltah v_p -\lambda \nablah \dvh v_p - \partial_{zz} v_p \bigr) \cdot \varepsilon \psi^h_t \idx\\
	& ~~~~ - \rho_0^{-1} \int \biggl( \varepsilon c^2_s \rho_{1,t} \dvh (\zeta \varepsilon \psi_t^h) -  \varepsilon \mu \nablah v_{p,t} : \nablah (\zeta \varepsilon \psi^h_t)\\
	& ~~~~ ~~~~ -  \varepsilon \lambda \dvh v_{p,t} \dvh ( \zeta \varepsilon \psi^h_t) - \varepsilon \dz v_{p,t} \cdot \dz (\zeta \varepsilon \psi^h_t) \biggr) \idx \\
	& ~~ \lesssim \varepsilon \bigl( \norm{\nablah \rho_1}{\Hnorm{1}} + \norm{\nabla v_p}{\Hnorm{2}} \bigr) \norm{\zeta_t}{\Lnorm{2}} \bigl( \norm{\varepsilon \nabla \psi^h_t}{\Lnorm{2}} + \norm{\varepsilon \psi_t^h}{\Lnorm{2}} \bigr) \\
	& ~~~~ + \varepsilon \bigl( \norm{\rho_{1,t}}{\Lnorm{2}} + \norm{\nabla v_{p,t}}{\Lnorm{2}} \bigr)  \norm{\zeta}{\Hnorm{2}} \bigl( \norm{\varepsilon \nabla \psi^h_t}{\Lnorm{2}} + \norm{\varepsilon \psi_t^h}{\Lnorm{2}} \bigr) \\
	& ~~ \lesssim \delta \norm{\varepsilon \nabla \psi^h_t}{\Lnorm{2}}^2 + \delta \norm{\varepsilon \psi^h_t}{\Lnorm{2}}^2 + \varepsilon^2 C_\delta \bigl( Q(\mathcal E) + 1 + \mathfrak G_p \bigr) \mathfrak H_p.
\end{aligned}
\end{equation}
We list estimates for $ I_4,I_6 $ below:
\begin{align}
	& \nonumber I_4
	= - \rho_0^{-1} c^2_s \varepsilon^2 \int \biggl( \xi_t \bigl(\rho_{1,tt} + \rho_{1,t} \dvh \psi^h + v\cdot \nablah \rho_{1,t} \bigr) \\
	& \nonumber ~~~~ + \xi_t \bigl(\rho_1 \dvh \psi^h_t + v_t \cdot \nablah \rho_1 \bigr) \biggr) \idx \lesssim \varepsilon^2 \norm{\xi_t}{\Lnorm2} \biggl( \norm{\rho_{1,tt}}{\Lnorm2} \\
	& \nonumber ~~~~ + \norm{\rho_{1,t}}{\Lnorm2} \norm{\nablah \psi^h}{\Hnorm{2}} + \norm{\nablah \rho_{1,t}}{\Lnorm{2}} \norm{v}{\Hnorm{2}} \\
	& ~~~~ + \norm{\rho_1}{\Hnorm{2}} \norm{\nablah \psi^h_t}{\Lnorm2}
	\nonumber + \norm{\nablah \rho_1}{\Hnorm{1}} \norm{v_t}{\Hnorm{1}} \biggr) \lesssim \delta \norm{\xi_t}{\Lnorm{2}}^2 \\
	&  ~~~~ + \delta \norm{\varepsilon \nablah \psi^h_t}{\Lnorm2}^2 + \delta \norm{\nablah \psi^h}{\Hnorm{2}}
	 +\varepsilon^2  C_\delta \bigl( Q(\mathcal E) + 1 + \mathfrak G_p \bigr) \mathfrak H_p, \label{estimate:005} \\
	&  I_6  = \int \mathcal R_t \dvh \psi^h_t \idx
	\lesssim \norm{\varepsilon^{-1} \zeta}{\Lnorm\infty} \norm{\zeta_t}{\Lnorm2}\norm{\varepsilon \nablah\psi^h_t}{\Lnorm2} \lesssim \delta \norm{\varepsilon \nabla \psi^h_t}{\Lnorm2}^2 \nonumber \\
	& ~~~~ ~~ + C_\delta \bigl( \varepsilon^2 \norm{\rho_1}{\Hnorm{2}}^2 + \norm{\varepsilon^{-1} \xi}{\Hnorm{2}}^2 \bigr) \bigl( \varepsilon^4 \norm{\rho_{1,t}}{2}^2 + \norm{\xi_t}{2}^2 \bigr), \label{estimate:006}
\end{align}
where we have used the fact that from \eqref{residue-est}
$$ \abs{\mathcal R_t}{} = \abs{\gamma (\gamma-1) \int_{\rho_0}^{\rho} \rho_t y^{\gamma-2} \,dy}{} \leq \abs{\gamma \rho_t (\rho^{\gamma -1} - \rho_0^{\gamma-1})}{} \leq C \abs{\zeta_t}{} \abs{\zeta}{}. $$
Summing up inequalities \eqref{estimate:001}, \eqref{estimate:002}, \eqref{estimate:003}, \eqref{estimate:004}, \eqref{estimate:005}, \eqref{estimate:006} and \eqref{ee:002} completes the proof.
\end{proof}

The next lemma follows directly from system \eqref{eq:perturbation}, and it shows the estimates of the temporal derivatives of $ \xi, \psi^h $ in terms of the spatial derivatives.
\begin{lemma}\label{lm:temporal-diffusion} Under the same assumptions as in Lemma \ref{lm:temporal-derivative},
	\begin{align}
	\label{ee:006}
		& \norm{\xi_t}{\Lnorm2}^2 \leq C \norm{\nablah \psi^h}{\Lnorm2}^2 + \varepsilon^2 Q(\mathcal E) \bigl( \norm{\nablah \psi^h}{\Hnorm{1}}^2 + \norm{\varepsilon^{-1} \nablah \xi}{\Hnorm{1}}^2 \bigr) \nonumber \\
		& ~~~~ +  \varepsilon^2 C \bigl( Q(\mathcal E) + 1 + \mathfrak G_p \bigr) \mathfrak H_p, \\
\label{ee:008}
		& \norm{\varepsilon \rho\psi^h_t}{\Lnorm2}^2 \leq C \norm{\varepsilon^{-1} \nablah \xi}{\Lnorm2}^2 + \varepsilon^2 C \norm{\nabla \psi^h}{\Hnorm{1}}^2 + \varepsilon^2 Q(\mathcal E)\bigl(  \norm{\varepsilon^{-1}\nablah \xi}{\Hnorm{1}}^2 \nonumber\\
		& ~~~~ + \norm{\nabla \psi^h}{\Hnorm{1}}^2 \bigr) + \varepsilon^2 C \bigl( Q(\mathcal E) + 1 + \mathfrak G_p \bigr) \mathfrak H_p,
\end{align}
for some positive constant $ C $ independent of $ \varepsilon $.
\end{lemma}
\begin{proof}
Indeed, after integrating \subeqref{eq:perturbation}{1} in the $ z $ variable, we have, thanks to \eqref{bc-perturbed-eq},
\begin{equation*}
	\dt \xi + \rho_0 \dvh \overline{\psi^h} = - \dvh (\xi \overline v) - \varepsilon^2 (\dt \rho_1 + \dvh (\rho_1 \bar v) ).
\end{equation*}
Then directly, we have, after applying the Minkowski, H\"older and Sobolev embedding inequalities,
\begin{align*}
	& \norm{\xi_{t}}{\Lnorm2} = \hnorm{\xi_t}{\Lnorm2} \lesssim \int_0^1 \hnorm{\nablah \psi^h}{\Lnorm2} \,dz +  \int_0^1 \bigl( \hnorm{\xi}{\Lnorm4} \hnorm{\nablah v}{\Lnorm4} + \hnorm{\nablah \xi}{\Lnorm4} \hnorm{v}{\Lnorm4} \\
	& ~~~~ + \varepsilon^2 ( \hnorm{\rho_{1,t}}{\Lnorm2} + \hnorm{\rho_1}{\Lnorm4} \hnorm{\nablah v}{\Lnorm4} + \hnorm{\nablah \rho_1}{\Lnorm4} \hnorm{v}{\Lnorm4}  ) \bigr) \,dz \lesssim \norm{\nablah\psi^h}{\Lnorm2} \\
	& ~~~~  + \norm{\xi}{H^1} \norm{\nablah \psi^h}{H^1} + \norm{\nablah \xi}{\Hnorm{1}} \norm{\psi^h}{\Hnorm{1}} + \varepsilon \norm{v_p}{\Hnorm{2}} \norm{\varepsilon^{-1}\xi}{H^2} \\
	& ~~~~ + \varepsilon^2  \bigl( \norm{\rho_{1,t}}{\Lnorm{2}} + \norm{\rho_1}{\Hnorm{2}} \norm{v_p}{\Hnorm{2}} + \norm{\rho_1}{\Hnorm{2}} \norm{\psi^h}{H^2} \bigr).
\end{align*}
Taking square on both sides of the above inequality yields \eqref{ee:006}.


On the other hand, from \subeqref{eq:perturbation}{2}, we have
\begin{align*}
	& \rho \dt \psi^h = - \rho v\cdot \nablah \psi^h - \rho w \dz \psi^h - \nablah (\varepsilon^{-2} c^2_s \xi) + \mu \deltah \psi^h + \lambda \nablah \dvh \psi^h \\
	& ~~~~ + \partial_{zz} \psi^h + \mathcal F_1 + \mathcal F_2 - \nablah (\varepsilon^{-2} \mathcal R).
\end{align*}
Directly, we have
\begin{align*}
	& \norm{\rho \psi^h_t}{\Lnorm 2} \lesssim \norm{\rho}{\Lnorm\infty} \norm{v}{\Lnorm\infty} \norm{\nablah \psi^h}{\Lnorm2}  + \varepsilon^{-1}\norm{\varepsilon^{-1}\nablah\xi}{\Lnorm2} + \norm{\nabla^2 \psi^h}{\Lnorm2} \\
	& ~~~~ + \norm{\varepsilon^{-1} \zeta}{\Lnorm3} \norm{\varepsilon^{-1}\nablah\zeta}{\Lnorm6} + \norm{\mathcal F_1}{\Lnorm2} + \norm{\mathcal F_2}{\Lnorm2}+ \norm{\rho w \dz \psi^h}{\Lnorm2} \\
	& \lesssim \varepsilon^{-1} \norm{\varepsilon^{-1}\nablah \xi}{\Lnorm2} + \norm{\nabla \psi^h}{\Hnorm{1}} +  Q(\mathcal E) \bigl( \norm{\varepsilon^{-1}\nablah \xi}{\Hnorm{1}}
	 + \norm{\nabla \psi^h}{\Lnorm{2}} \bigr)\\
	& ~~~~ + \bigl( Q(\mathcal E) + \varepsilon \norm{\rho_1}{\Hnorm{1}} \bigr) \varepsilon \norm{\rho_1}{\Hnorm{2}} +  \norm{v_p}{\Hnorm{2}} \norm{\nablah \psi^h}{\Lnorm{2}} \\
	& ~~~~ + \norm{\mathcal F_1}{\Lnorm2} + \norm{\mathcal F_2}{\Lnorm2}+ \norm{\rho w \dz \psi^h}{\Lnorm2},
\end{align*}
where we have substituted the following inequality from \eqref{residue-est},
\begin{equation}\label{residue-est-h}
	\abs{\nablah \mathcal R}{} =\abs{\gamma(\gamma-1) \int_{\rho_0}^{\rho} \nablah \rho y^{\gamma-2} \,dy}{} \leq \abs{\gamma \nablah \rho (\rho^{\gamma-1} - \rho_0^{\gamma-1})}{} \lesssim \abs{\nablah \zeta}{} \abs{\zeta}{}.
\end{equation}
From \eqref{id:vertical_velocity}, \eqref{id:vertical_perturbation}, one has,
\begin{align*}
	& \rho w \dz \psi^h = \rho w_p \dz \psi^h - \int_0^z \rho \dvh \widetilde{\psi^h} + \widetilde v\cdot\nablah \rho \,dz  \dz \psi^h, \\
	& \mathcal F_1 = \zeta Q_p,\\
	& \mathcal F_2 = - \rho \psi^h \cdot\nablah v_p + \bblparenthese \int_0^z \bigl( \rho \dvh \widetilde{\psi^h} + \widetilde{v}\cdot \nablah \rho \bigr) \,dz \bbrparenthese \dz v_p.
\end{align*}
Therefore,
\begin{align*}
	& \norm{\mathcal F_1}{\Lnorm2} + \norm{\mathcal F_2}{\Lnorm2}+ \norm{\rho w \dz \psi^h}{\Lnorm2} \lesssim \norm{\zeta}{\infty} \norm{Q_p}{\Lnorm2} \\
	& ~~~~ + \norm{\rho}{\Lnorm\infty}\bigl( \norm{w_p}{\Lnorm3} \norm{\dz \psi^h}{\Lnorm6} + \norm{\psi^h}{\Lnorm\infty} \norm{\nablah v_p}{\Lnorm2} \\
	& ~~~~+ \norm{\nablah\psi^h}{\Lnorm3} \norm{\dz v_p}{\Lnorm6} + \norm{\nablah\psi^h}{\Lnorm3} \norm{\dz \psi^h}{\Lnorm6} \bigr) \\
	& ~~~~ + \norm{v}{\Lnorm\infty} \norm{\nablah \zeta}{\Lnorm3}\norm{\dz \psi^h}{\Lnorm6} + \norm{v}{\Lnorm\infty} \norm{\nablah \zeta}{\Lnorm3}\norm{\dz v_p}{\Lnorm6}\\
	& \lesssim  \bigl( \varepsilon^2 \norm{\rho_1}{\Hnorm{2}} + \varepsilon \norm{\varepsilon^{-1} \xi}{\Hnorm{2}}\bigr)\bigl( \norm{\rho_1}{\Hnorm{1}} + \norm{v_p}{\Hnorm{2}} \bigr) \\
	& ~~~~ + \bigl( \norm{w_p}{\Hnorm{1}} + \norm{v_p}{\Hnorm{2}} + \norm{\nabla \psi^h}{\Hnorm{1}} \bigr) \norm{\psi^h}{\Hnorm{2}}\\
	& ~~~~ +\bigl( \norm{v_p}{\Hnorm{2}} + \norm{\psi^h}{\Hnorm{2}} \bigr)^2 \bigl( \varepsilon^2 \norm{\rho_1}{\Hnorm{2}} + \varepsilon \norm{\varepsilon^{-1} \nablah \xi}{\Hnorm{1}} \bigr).
\end{align*}
Summing the above inequalities yields \eqref{ee:008}.
\end{proof}

%
%

\subsection{Horizontal derivatives}
We derive the required estimates for the horizontal derivatives in this subsection. 
After applying $ \partial_{hh} = \partial_h^2 $ to system \eqref{eq:perturbation},  we obtain the following system:
\begin{equation}\label{eq:perturbation-020}
\begin{cases}
\dt \xi_{hh} + \rho_0 ( \dvh \psi^h_{hh} + \dz \psi^z_{hh} ) = \mathcal G_{1,hh} + \mathcal G_{2,hh} & \text{in} ~ \Omega, \\
\rho \dt \psi^h_{hh} + \rho v\cdot \nablah \psi^h_{hh} + \rho w \partial_z \psi^h_{hh} + \nablah ( \varepsilon^{-2} c^2_s \xi_{hh}) \\
~~~~ = \mu \deltah \psi^h_{hh}
 + \lambda\nablah \dvh \psi^h_{hh} + \partial_{zz} \psi^h_{hh}  + \mathcal F_{1,hh} \\
 ~~~~~ ~~~~~ + \mathcal F_{2,hh} - \nablah (\varepsilon^{-2} \mathcal R_{hh})
 + \mathcal H_{hh} & \text{in} ~ \Omega ,
\end{cases}
\end{equation}
where
\begin{equation}\label{nonlinear-dhh}
	\begin{aligned}
	& \mathcal{H}_{hh} : = - \rho_{hh}\psi^h_t - (\rho v)_{hh} \cdot\nablah \psi^h - ( \rho w)_{hh} \dz \psi^h - 2 \rho_h \psi^h_{ht} \\
	& ~~~~ - 2 (\rho v)_h \cdot \nablah \psi^h_h - 2 (\rho w)_h\dz \psi^h_h.
	\end{aligned}
\end{equation}

\begin{lemma}\label{lm:horizontal-derivative}
Under the same assumptions as in Lemma \ref{lm:temporal-derivative}, we have
	\begin{equation}\label{ee:horizontal-derivative}
	\begin{aligned}
		& \dfrac{d}{dt} \biggl\lbrace \dfrac{1}{2} \norm{\rho^{1/2}\psi^h_{hh}}{\Lnorm 2}^2 + \dfrac{c^2_s}{2\rho_0} \norm{\varepsilon^{-1}\xi_{hh}}{\Lnorm 2}^2  \biggr\rbrace + C_{\mu,\lambda} \norm{\nabla \psi^h_{hh} }{\Lnorm 2}^2 \\
		& ~~~~ \leq \delta \bigl( \norm{\nabla \psi^h}{\Hnorm{2}}^2  +  \norm{\varepsilon^{-1} \nablah \xi}{\Hnorm{1}}^2 +  \norm{\varepsilon \nabla \psi^h_t}{\Lnorm{2}}^2 \\
		& ~~~~ +  \norm{\varepsilon\psi^h_t}{\Lnorm{2}}^2 \bigr) + C_\delta Q(\mathcal E) \bigl( \norm{\nabla\psi^h}{\Hnorm{1}}^2 + \norm{\varepsilon^{-1} \nablah \xi}{\Hnorm{1}}^2 \bigr) \\
		& ~~~~ + C_\delta \bigl( Q(\mathcal E) + 1 + \mathfrak G_p \bigr) \mathfrak H_p \bigl( \norm{\psi^h}{\Hnorm{2}}^2 + \norm{\varepsilon^{-1} \xi}{\Hnorm{2}}^2 + \varepsilon^2 \bigr),
	\end{aligned}
\end{equation}
for some positive constant $ C_{\mu,\lambda} $, which is independent of $ \varepsilon $.
\end{lemma}

\begin{proof}
Take the inner product of \subeqref{eq:perturbation-020}{2} with $ \psi^h_{hh} $ and integrate the resultant over $ \Omega $. Similarly as before, we will have the following:
\begin{align*}
	& \dfrac{d}{dt} \biggl\lbrace \dfrac 1 2 \norm{\rho^{1/2}\psi^h_{hh}}{\Lnorm 2}^2 \biggr\rbrace - \int \varepsilon^{-2} c^2_s  \xi_{hh} \dvh \psi^h_{hh} \idx + \mu \norm{\nablah \psi^h_{hh}}{\Lnorm2}^2 \\
	& ~~~~ + \lambda \norm{\dvh\psi^h_{hh}}{\Lnorm 2}^2 + \norm{\dz \psi^h_{hh}}{\Lnorm2}^2 = \int \mathcal F_{1,hh} \cdot \psi^h_{hh} \idx \\
	& ~~~~ + \int \mathcal F_{2,hh} \cdot \psi^h_{hh} \idx - \int \nablah(\varepsilon^{-2}\mathcal R_{hh}) \cdot \psi^h_{hh}\idx + \int \mathcal H_{hh} \cdot \psi^h_{hh} \idx.
\end{align*}
Again, we employ equation \subeqref{eq:perturbation-020}{1} to derive that,
\begin{align*}
	& - \int \varepsilon^{-2} c^2_s  \xi_{hh} \dvh \psi^h_{hh} \idx = \int \bblparenthese  \varepsilon^{-2} \rho_0^{-1} c^2_s \xi_{hh} \bigl( \dt \xi_{hh} + \rho_0 \dz \psi^z_{hh} - \mathcal G_{1,hh} \\
	& ~~~~ - \mathcal G_{2,hh}  \bigr) \bbrparenthese \idx = \dfrac{d}{dt} \biggl\lbrace \dfrac{c^2_s}{2\rho_0} \int \varepsilon^{-2} \abs{\xi_{hh}}{2} \idx \biggr\rbrace - \rho_0^{-1} c^2_s \int \varepsilon^{-2} \xi_{hh} \mathcal G_{1,hh} \idx \\
	& ~~~~ - \rho_0^{-1} c^2_s \int \varepsilon^{-2} \xi_{hh} \mathcal G_{2,hh} \idx.
\end{align*}
Therefore, the following equality holds:
\begin{equation}\label{ee:003}
\begin{aligned}
	& \dfrac{d}{dt} \biggl\lbrace \dfrac{1}{2} \norm{\rho^{1/2}\psi^h_{hh}}{\Lnorm2}^2 + \dfrac{c^2_s}{2\rho_0} \norm{\varepsilon^{-1}\xi_{hh}}{\Lnorm2}^2  \biggr\rbrace + \mu \norm{\nablah \psi^h_{hh}}{\Lnorm2}^2 \\
	& ~~~~ + \lambda \norm{\dvh\psi_{hh}^h}{\Lnorm2}^2 + \norm{\dz \psi^h_{hh}}{\Lnorm2}^2 = \int \mathcal F_{1,hh} \cdot \psi^h_{hh}\idx \\
	& ~~~~ + \int \mathcal F_{2,hh} \cdot \psi^h_{hh}\idx + \rho_0^{-1} c^2_s \int \varepsilon^{-2} \xi_{hh} \mathcal G_{1,hh}\idx \\
	& ~~~~ + \rho_0^{-1} c^2_s \int \varepsilon^{-2} \xi_{hh} \mathcal G_{2,hh} \idx + \int \mathcal H_{hh} \cdot \psi_{hh}^h \idx \\
	& ~~~~ - \int \nablah ( \varepsilon^{-2} \mathcal R_{hh}) \cdot \psi_{hh}^h \idx =: \sum_{i=7}^{12} I_i.
\end{aligned}
\end{equation}
In the following, we shall estimate the terms on the right-hand side of \eqref{ee:003}.
$ I_{11} $ can be written, after substituting \eqref{nonlinear-dhh}, as
\begin{align*}
	& I_{11} = - \int \biggl( \zeta_{hh} \psi_t^h \cdot \psi_{hh}^h + 2 \zeta_h \psi_{ht}^h \cdot \psi_{hh}^h \biggr) \idx - \int \biggl( (\rho v)_{hh} \cdot \nablah \psi^h \cdot \psi_{hh}^h \\
	& ~~~~ + 2 (\rho v)_h \cdot \nablah \psi_h^h \cdot \psi_{hh}^h \biggr) \idx
	- \int \biggl( ( \rho w)_{hh}  \dz \psi^h \cdot \psi_{hh}^h \\
	& ~~~~ + 2(\rho w)_h \dz \psi_h^h \cdot \psi_{hh}^h \biggr) \idx =: I_{11}' + I_{11}'' + I_{11}'''.
\end{align*}
Straightforward estimates show the following inequalities:
\begin{align*}
	& I_{11}' \lesssim \hnorm{\zeta_{hh}}{\Lnorm2}\int_0^1 \hnorm{\psi^h_t}{\Lnorm4} \hnorm{\psi_{hh}^h}{\Lnorm4} \,dz + \hnorm{\zeta_h}{\Lnorm4} \int_0^1 \hnorm{\psi^h_{ht}}{\Lnorm2} \hnorm{\psi^h_{hh}}{\Lnorm4} \,dz \\
	& ~~ \lesssim \norm{\varepsilon^{-1}\nabla \zeta}{\Hnorm{1}} \bigl(  \norm{\varepsilon\nablah \psi_t^h}{\Lnorm2} + \norm{\varepsilon\psi_t^h}{\Lnorm2}\bigr) \\
	& ~~~~ \times \bigl( \norm{\psi^h_{hh}}{\Lnorm2}^{1/2}\norm{\nablah \psi^h_{hh}}{\Lnorm2}^{1/2} + \norm{\psi^h_{hh}}{\Lnorm2} \bigr)
	\lesssim \delta \norm{\nablah \psi_{hh}^h}{\Lnorm2}^2 \\
	& ~~~~ + \delta \norm{\varepsilon \nablah\psi^h_t}{\Lnorm{2}}^2 + \delta \norm{\varepsilon\psi^h_t}{\Lnorm{2}}^2 + C_\delta Q(\mathcal E) \norm{\nabla \psi^h}{\Hnorm{1}}^2 \\
	& ~~~~ + C_\delta \bigl( \mathfrak G_p + 1 \bigr) \mathfrak H_p \norm{\psi^h}{\Hnorm{2}}^2, \\
	& I_{11}'' \lesssim \bigl( \norm{\zeta_{hh}}{\Lnorm2} \norm{v}{\Lnorm\infty} + \norm{\rho}{\Lnorm\infty} \norm{v_{hh}}{\Lnorm2} + \norm{\zeta_h}{\Lnorm3} \norm{v_h}{\Lnorm6} \bigr)\\
	& ~~~~ \times \norm{\nablah \psi^h}{\Lnorm6} \norm{\psi_{hh}^h}{\Lnorm3} + \bigl( \norm{\zeta_h}{\Lnorm6} \norm{v}{\Lnorm\infty} + \norm{\rho}{\Lnorm\infty} \norm{v_h}{\Lnorm6} \bigr) \\
	& ~~~~\times \norm{\nablah^2 \psi^h}{\Lnorm2} \norm{\psi_{hh}^h}{\Lnorm3} \lesssim \bigl( \norm{\nablah \zeta}{\Hnorm{1}} + 1 \bigr) \norm{v}{\Hnorm{2}} \norm{\nabla \psi^h}{\Hnorm{1}} \\
	& ~~~~ \times \bigl( \norm{\psi^h_{hh}}{\Lnorm2}^{1/2} \norm{\nabla \psi^h_{hh}}{\Lnorm2}^{1/2} + \norm{\psi_{hh}^h}{\Lnorm2}\bigr) \lesssim \delta \norm{\nabla \psi_{hh}^h}{\Lnorm2}^2 \\
	& ~~~~ + C_\delta Q(\mathcal E) \norm{\nabla \psi^h}{\Hnorm{1}}^2 + C_\delta \bigl( Q(\mathcal E) + 1 + \mathfrak G_p \bigr) \mathfrak H_p \norm{\psi^h}{\Hnorm{2}}^2.
\end{align*}
In order to estimate $ I_{11}''' $, we first take $ \partial_h, \partial_{hh} $ in \eqref{id:vertical_velocity} to obtain the identities:
\begin{equation}\label{id:vertical-velocity-020}
	\begin{aligned}
		w_{h} & = w_{p,h} - \int_0^z \biggl( \dvh \widetilde{\psi^h_h} + \widetilde v_h \cdot \nablah \log \rho + \widetilde v \cdot \nablah (\log \rho)_h \biggr) \,dz \\
		& = w_{p,h} + \psi^z_h, \\
		w_{hh} & = w_{p,hh} - \int_0^z \biggl( \dvh \widetilde{\psi^h_{hh}} + \widetilde v_{hh} \cdot\nablah \log \rho + 2 \widetilde v_h \cdot \nablah (\log \rho)_h \\
		& ~~~~ + \widetilde v\cdot\nablah (\log \rho)_{hh} \biggr) \,dz = w_{p,hh} + \psi^z_{hh}.
	\end{aligned}
\end{equation}
Then $ I_{11}''' $ can be written as
\begin{align*}
	& I_{11}''' = - \int (\zeta_{hh} w_p + 2 \zeta_h w_{p,h} + \rho w_{p,hh} )\dz\psi^h \cdot\psi^h_{hh} \idx \\
	& ~~~~ - 2 \int ( \zeta_h w_p + \rho w_{p,h}) \dz \psi^h_h \cdot\psi_{hh}^h \idx + \int \biggl( \int_0^z \rho \widetilde v \cdot \nablah (\log \rho)_{hh} \,dz \\
	& ~~~~ \times \bigl( \dz \psi^h \cdot\psi^h_{hh} \bigr) \biggr) \idx + \int \biggl\lbrack \int_0^z \biggl( \zeta_{hh}  \dvh \widetilde{\psi^h} + \zeta_{hh} \widetilde v \cdot \nablah \log \rho \\
	& ~~~~  + 2 \zeta_h \dvh \widetilde{\psi^h_{h}}
	 + 2 \zeta_h \widetilde v_h \cdot\nablah \log \rho  + 2 \zeta_h \widetilde v \cdot \nablah (\log \rho)_h + \rho \dvh \widetilde{\psi^h_{hh}} \\
	& ~~~~ + \rho \widetilde{v}_{hh} \cdot\nablah\log \rho
	 + 2 \rho \widetilde{v}_{h} \cdot \nablah (\log \rho)_h  \biggr) \,dz \times  \bigl( \dz \psi^h \cdot \psi^h_{hh} \bigr) \biggr\rbrack \idx \\
	& ~~~~ + 2 \int \biggl\lbrack \int_0^z \biggl( \zeta_h \dvh \widetilde{\psi^h}
	 + \zeta_h \widetilde v \cdot \nablah \log \rho + \rho\dvh\widetilde{\psi^h_h} + \rho \widetilde v_h \cdot\nablah \log \rho \\
	& ~~~~  + \rho \widetilde v \cdot\nablah (\log \rho)_h  \biggr) \,dz \dz \times  \bigl( \psi^h_h \cdot \psi_{hh}^h \bigr) \biggr\rbrack \idx =: \sum_{i=1}^5 I_{11,i}'''.
\end{align*}
Similarly, after applying the H\"older, Sobolev embedding and Young inequalities, one has
\begin{align*}
& I_{11,1}''' + I_{11,2}''' \lesssim \bigl(\norm{\zeta_{hh}}{\Lnorm2}\norm{w_p}{\Lnorm\infty} + \norm{\zeta_h}{\Lnorm6} \norm{w_{p,h}}{\Lnorm3} \\
& ~~~~ + \norm{\rho}{\Lnorm\infty}\norm{w_{p,hh}}{\Lnorm2} \bigr) \norm{\dz \psi^h}{\Lnorm6} \norm{\psi^h_{hh}}{\Lnorm3} \\
& ~~~~ + \bigl(\norm{\zeta_h}{\Lnorm6}\norm{w_p}{\Lnorm6} + \norm{\rho}{\Lnorm\infty}\norm{w_{p,h}}{\Lnorm3}\bigr) \norm{\dz \psi^h_{h}}{\Lnorm2} \norm{\psi^h_{hh}}{\Lnorm6}\\
& \lesssim \delta \norm{\nabla\psi^h_{hh}}{\Lnorm2}^2 + \delta \norm{\psi^h_{hh}}{\Lnorm{2}}^2 + C_\delta \bigl( Q(\mathcal E) + 1 + \mathfrak G_p \bigr) \mathfrak H_p \norm{\psi^h}{\Hnorm{2}}^2.
\end{align*}
On the other hand, after applying integration by parts, $ I_{11,3}''' $ is written as
\begin{align*}
	& I_{11,3}''' = - \int  \biggl( \int_0^z (\log \rho)_{hh} (\rho \dvh \widetilde{ v} + \widetilde v \cdot \nablah \zeta ) \,dz \biggr) \times \bigl( \dz \psi^h \cdot\psi^h_{hh} \bigr) \idx \\
	& ~~~~ - \int \biggl( \int_0^z (\log \rho)_{hh} \rho \widetilde{v} \,dz \biggr) \cdot \nablah (\dz \psi^h \cdot \psi^h_{hh}) \idx,
\end{align*}
which can be estimated as below:
\begin{align*}
	& I_{11,3}''' \lesssim \int_0^1 ( \hnorm{\zeta_{hh}}{\Lnorm2} + \hnorm{\zeta_h}{\Lnorm4}^2   )(\hnorm{\rho}{\Lnorm\infty}\hnorm{\nablah v}{\Lnorm8} + \hnorm{v}{\Lnorm\infty} \hnorm{\nablah \zeta}{\Lnorm8}) \,dz \\ & ~~~~ \times \int_0^1 \hnorm{\dz \psi^h}{\Lnorm8} \hnorm{\psi_{hh}^h}{\Lnorm4} \,dz + \int_0^1 ( \hnorm{\zeta_{hh}}{\Lnorm2} + \hnorm{\zeta_h}{\Lnorm4}^2  ) \hnorm{\rho}{\Lnorm\infty} \hnorm{v}{\Lnorm\infty} \,dz \\ & ~~~~ \times \int_0^1 (\hnorm{\dz \psi^h}{\Lnorm\infty} \hnorm{\nablah \psi_{hh}^h}{\Lnorm2} + \hnorm{\nablah \dz \psi^h}{\Lnorm4} \hnorm{\psi_{hh}^h}{\Lnorm4}) \,dz \\
	& ~~ \lesssim \bigl( \norm{ \nablah \zeta}{\Hnorm{1}} + \norm{\nablah \zeta}{\Hnorm{1}}^2 \bigr)\bigl(\norm{\rho}{\Lnorm\infty} \norm{v}{\Hnorm{2}} + \norm{v}{\Hnorm{2}} \norm{\zeta}{\Hnorm{2}}\bigr)\\
	& ~~~~ \times \norm{\dz\psi^h}{\Hnorm{1}} \bigl( \norm{\psi^h_{hh}}{\Lnorm2}^{1/2} \norm{\nablah \psi^h_{hh}}{\Lnorm2}^{1/2} + \norm{\psi^h_{hh}}{\Lnorm2} \bigr) \\
	& ~~ + \bigl( \norm{\nablah \zeta}{\Hnorm{1}} + \norm{\nablah \zeta}{\Hnorm{1}}^2\bigr) \norm{\rho}{\Lnorm\infty} \norm{v}{\Hnorm{2}}  \bigl( \norm{\dz \psi^h}{\Lnorm{2}}^{1/2} \\
	& ~~~~ \times \norm{\dz \psi^h}{\Hnorm{2}}^{1/2} \norm{\nablah \psi^h_{hh}}{\Lnorm2} + ( \norm{\nablah \dz \psi^h}{\Lnorm{2}}^{1/2} \norm{\nablah^2 \dz \psi^h}{\Lnorm{2}}^{1/2} \\
	& ~~~~ + \norm{\nablah \dz \psi^h}{\Lnorm{2}}) ( \norm{\psi^h_{hh}}{\Lnorm2}^{1/2} \norm{\nablah \psi^h_{hh}}{\Lnorm2}^{1/2} \\
	& ~~~~ + \norm{\psi^h_{hh}}{\Lnorm2} )  \bigr)
	\lesssim \delta \norm{\nabla^3 \psi^h}{\Lnorm2}^2 + C_\delta Q(\mathcal E) \norm{\nabla \psi^h}{\Hnorm{1}}^2 \\
	& ~~~~ + C_\delta \bigl( Q(\mathcal E) + 1 + \mathfrak G_p \bigr) \mathfrak H_p \norm{\psi^h}{\Hnorm{2}}^2 .
\end{align*}
Applying the H\"older, Minkowski, Sobolev embedding and Young inequalities implies,
\begin{align*}
	& I_{11,4}''' \lesssim \int_0^1 \biggl( \hnorm{\zeta_{hh}}{\Lnorm2}\hnorm{\nablah \psi^h}{\Lnorm8} + \hnorm{\zeta_{hh}}{\Lnorm2} \hnorm{v}{\Lnorm\infty}\hnorm{\nablah \zeta}{\Lnorm8} + \hnorm{\zeta_h}{\Lnorm8} \hnorm{\nablah \psi^h_h}{\Lnorm2} \\
	& ~~~~ + \hnorm{\zeta_h}{\Lnorm8}\hnorm{v_h}{\Lnorm4}\hnorm{\nablah \zeta}{\Lnorm4} + \hnorm{\zeta_h}{\Lnorm8} \hnorm{v}{\Lnorm\infty} \bigl( \hnorm{\nablah^2\zeta}{\Lnorm2} + \hnorm{\nablah \zeta}{\Lnorm4}^2\bigr)  \\
	& ~~~~ + \hnorm{\rho}{\Lnorm 8} \hnorm{\nablah \psi^h_{hh}}{\Lnorm2} + \hnorm{\rho}{\Lnorm\infty}\hnorm{v_{hh}}{2} \hnorm{\nablah\zeta}{\Lnorm8} + \hnorm{\rho}{\Lnorm\infty} \hnorm{v_h}{\Lnorm8} \\
	& ~~~~ \times \bigl( \hnorm{\nablah^2 \zeta}{\Lnorm2} + \hnorm{\nablah \zeta}{\Lnorm4}^2\bigr) \biggr) \,dz  \times
	\int_0^1 \hnorm{\dz \psi^h}{\Lnorm8} \hnorm{\psi^h_{hh}}{\Lnorm4} \,dz \\
	& \lesssim \biggl( \norm{\nablah \zeta}{\Hnorm{1}} \norm{\nabla \psi^h}{\Hnorm{1}} + ( \norm{\nablah \zeta}{\Hnorm{1}}^2+ \norm{\nablah \zeta}{\Hnorm{1}}^3) \norm{v}{\Hnorm{2}} \\
	& ~~~~ + \norm{\rho}{\Hnorm{1}} \norm{\nablah \psi^h_{hh}}{\Lnorm2} + \norm{\rho}{\Hnorm{2}} \norm{v}{\Hnorm{2}} ( \norm{\nablah \zeta}{\Hnorm{1}} + \norm{\nablah \zeta}{\Hnorm{1}}^2 )  \biggr) \\
	& ~~~~ \times \norm{\dz \psi^h}{\Hnorm{1}} \times \bigl(\norm{\psi^h_{hh}}{2}^{1/2} \norm{\nablah \psi^h_{hh}}{2}^{1/2} + \norm{\psi^h_{hh}}{2}\bigr)   \\
	& \lesssim \delta \norm{\nabla\psi^h_{hh}}{2}^2 + C_\delta  Q(\mathcal E) \norm{\nabla \psi^h}{\Hnorm{1}}^2 + C_\delta \bigl( Q(\mathcal E) + 1 + \mathfrak G_p \bigr) \mathfrak H_p \norm{\psi^h}{\Hnorm{2}}^2,\\
	& I_{11,5}''' \lesssim  \int_0^1 \biggl( \hnorm{\zeta_h}{\Lnorm4} \hnorm{\nablah \psi^h}{\Lnorm4} + \hnorm{\zeta_h}{\Lnorm4} \hnorm{v}{\Lnorm\infty} \hnorm{\nablah \zeta}{\Lnorm4} + \hnorm{\rho}{\Lnorm\infty}\hnorm{\nablah \psi_h^h}{\Lnorm2}\\
	& ~~~~ + \hnorm{\rho}{\Lnorm\infty} \hnorm{v_h}{\Lnorm4}\hnorm{\nablah \zeta}{\Lnorm4} + \hnorm{\rho}{\Lnorm\infty}\hnorm{v}{\Lnorm\infty} (\hnorm{\nablah^2\zeta}{\Lnorm2} + \hnorm{\nablah \zeta}{\Lnorm4}^2) \biggr)  \,dz \\
	& ~~~~ \times \int_0^1 \hnorm{\dz \psi^h_h}{\Lnorm4} \hnorm{\psi^h_{hh}}{\Lnorm4}  \,dz  \lesssim \biggl( \norm{\zeta}{\Hnorm{2}} \norm{\psi^h}{\Hnorm{2}} + \norm{v}{\Hnorm{2}} \norm{\zeta}{\Hnorm{2}}^2 \\
	& ~~~~ + \norm{\rho}{\Lnorm\infty} \norm{\psi^h}{\Hnorm{2}} + \norm{\rho}{\Lnorm\infty} \norm{v}{\Hnorm{2}} \norm{\zeta}{\Hnorm{2}} + \norm{\rho}{\Lnorm\infty}\norm{v}{\Hnorm{2}}\\
	& ~~~~ \times ( \norm{\zeta}{\Hnorm{2}}+ \norm{\zeta}{\Hnorm{2}}^2 )   \biggr)  \bigl( \norm{\dz \psi^h_h}{\Lnorm2}^{1/2} \norm{\nablah \dz \psi^h_h}{\Lnorm2}^{1/2} + \norm{\dz \psi_h^h}{\Lnorm2} \bigr) \\
	& ~~~~ \times \bigl( \norm{\psi^h_{hh}}{\Lnorm2}^{1/2} \norm{\nablah \psi^h_{hh}}{\Lnorm2}^{1/2} + \norm{\psi^h_{hh}}{\Lnorm2}\bigr) \\
	& \lesssim \delta \norm{\nabla\psi^h_{hh}}{\Lnorm2}^2 + C_\delta  Q(\mathcal E) \norm{\nabla \psi^h}{\Hnorm{1}}^2 + C_\delta \bigl( Q(\mathcal E) + 1 + \mathfrak G_p\bigr) \mathfrak H_p \norm{\psi^h}{\Hnorm{2}}^2.
\end{align*}
Therefore, we have established the following estimate
\begin{equation}\label{estimate:007}
	\begin{aligned}
		& I_{11} \lesssim \delta \norm{\nabla^3 \psi^h}{\Lnorm2}^2 + \delta \norm{\psi^h_{hh}}{\Lnorm{2}}^2 + \delta \norm{\varepsilon \nabla \psi^h_t}{\Lnorm{2}}^2 + \delta \norm{\varepsilon \psi^h_t}{\Lnorm{2}}^2 \\
		& ~~~~ + C_\delta Q(\mathcal E)  \norm{\nabla \psi^h}{\Hnorm{1}}^2 + C_\delta \bigl( Q(\mathcal E) + 1 + \mathfrak G_p  \bigr) \mathfrak H_p \norm{\psi^h}{\Hnorm{2}}^2.
	\end{aligned}
\end{equation}
Other terms on the right-hand side of \eqref{ee:003} can be estimated in the same manner. We will list them below.
In order to estimate $ I_8 $, notice
\begin{align*}
	& \mathcal F_{2,hh} = - ( \rho \psi^h)_{hh} \cdot \nablah v_p - (\rho \psi^z)_{hh} \dz v_p - 2 (\rho \psi^h)_h \cdot\nablah v_{p,h} \\
	& ~~~~ - 2 (\rho \psi^z)_h \dz v_{p,h} - \rho \psi^h \cdot\nablah v_{p,hh} - \rho \psi^z \dz v_{p,hh}.
\end{align*}
Then $ I_8 $ can be written as
\begin{align*}
	& I_8 = - \int \bigl( (\rho\psi^h)_{hh}\cdot \nablah v_p + 2 (\rho\psi^h)_h \cdot\nablah v_{p,h} + \rho \psi^h \cdot\nablah v_{p,hh}  \bigr) \cdot\psi^h_{hh} \idx \\
	& ~~~~ - \int \bigl( (\rho\psi^z)_{hh} \dz v_p + 2 ( \rho\psi^z )_h \dz v_{p,h} + \rho \psi^z \dz v_{p,hh} \bigr) \cdot \psi_{hh}^h \idx : = I_8' + I_8'',
\end{align*}
with the estimate
\begin{align*}
	& I_8'
	\lesssim \biggl( \norm{\zeta_{hh}}{\Lnorm2} \norm{\psi^h}{\Lnorm\infty} \norm{\nablah v_p}{\Lnorm6} + \norm{\zeta_h}{\Lnorm6}\norm{\psi^h_h}{\Lnorm3} \norm{\nablah v_p}{\Lnorm6} \\
	& ~~~~ + \norm{\rho}{\Lnorm\infty} \norm{\psi^h_{hh}}{\Lnorm2}\norm{\nablah v_p}{\Lnorm6} + \norm{\zeta_h}{\Lnorm6} \norm{\psi^h}{\Lnorm\infty} \norm{\nablah v_{p,h}}{\Lnorm2} \\
	& ~~~~ + \norm{\rho}{\Lnorm\infty} \norm{\psi_h^h}{\Lnorm6} \norm{\nablah v_{p,h}}{\Lnorm2} + \norm{\rho}{\Lnorm\infty} \norm{\psi^h}{\Lnorm6} \norm{\nablah v_{p,hh}}{\Lnorm2} \biggr) \\
	& ~~~~ \times \norm{\psi^h_{hh}}{3} \lesssim \bigl( \norm{\nablah\zeta}{\Hnorm{1}} \norm{\psi^h}{\Hnorm{2}}\norm{v_p}{\Hnorm{2}}
	 + \norm{\rho}{\Lnorm{\infty}} \norm{\psi^h}{\Hnorm{2}} \norm{v_p}{\Hnorm{3}} \bigr) \\
	 & ~~~~ \times  \bigl( \norm{\psi^h_{hh}}{\Lnorm2}^{1/2} \norm{\nabla \psi^h_{hh}}{\Lnorm2}^{1/2} + \norm{\psi^h_{hh}}{\Lnorm2} \bigr) \\
	& \lesssim \delta \norm{\nabla \psi^h_{hh}}{\Lnorm2}^2 + \delta \norm{\psi^h_{hh}}{\Lnorm{2}}^2 + C_\delta \bigl( Q(\mathcal E) + 1 + \mathfrak G_p \bigr) \mathfrak H_p \norm{\psi^h}{\Hnorm{2}}^2.
\end{align*}
Next, after substituting \eqref{id:vertical_velocity-001} and \eqref{id:vertical-velocity-020} into the term $ I_{8}'' $, we have the following:
\begin{align*}
	& I_8'' = \int \biggl\lbrack \int_0^z \biggl( \zeta_{hh} \bigl(\dvh \widetilde{\psi^h} + \widetilde v\cdot\nablah \log \rho \bigr) + 2 \zeta_h \bigl( \dvh \widetilde{\psi^h_h} \\
	& ~~~~ + \widetilde v_h \cdot\nablah \log \rho
	 + \widetilde v \cdot\nablah (\log \rho )_h \bigr) + \rho \bigl( \dvh \widetilde{\psi^h_{hh}} + \widetilde v_{hh} \cdot\nablah \log \rho \\
	& ~~~~ + 2 \widetilde v_h \cdot\nablah (\log \rho)_h \bigr) \biggr) \,dz \times \bigl( \dz v_p \cdot \psi^h_{hh} \bigr) \biggr\rbrack \idx
	 + \int \biggl\lbrack \int_0^z \rho \widetilde v \cdot \nablah(\log \rho)_{hh} \,dz \\
	& ~~~~ \times \bigl( \dz v_p \cdot \psi^h_{hh} \bigr) \biggr\rbrack \idx
	 + 2 \int \bblbrack \int_0^z \bblparenthese \zeta_h \bigl( \dvh \widetilde{\psi^h} + \widetilde v\cdot \nablah \log \rho \bigr)
	+ \rho \bigl( \dvh \widetilde{\psi^h_h} \\
	& ~~~~ + \widetilde v_h \cdot\nablah \log \rho + \widetilde v \cdot \nablah (\log \rho)_h\bigr) \bbrparenthese \,dz \times \blparenthese  \dz v_{p,h} \cdot\psi^h_{hh} \brparenthese \bbrbrack \idx \\
	& ~~~~ + \int \bblbrack \int_0^z \rho \bigl(\dvh \widetilde{\psi^h} + \widetilde v\cdot\nablah \log \rho \bigr) \,dz \times \blparenthese \dz v_{p,hh} \cdot \psi^h_{hh} \brparenthese \bbrbrack \idx \\
	& ~~~~ =: I_{8,1}'' + I_{8,2}'' + I_{8,3}'' + I_{8,4}''.
\end{align*}
Applying integration by parts in $ I_{8,2}'' $ yields the following:
\begin{align*}
	& I_{8,2}'' = - \int \bblbrack \int_0^z \bblparenthese (\log \rho)_{hh} \rho \dvh{\widetilde v} + (\log \rho)_{hh} \widetilde v \cdot
	\nablah \rho \bbrparenthese \,dz \\
	& ~~~~ \times \blparenthese \dz v_p \cdot\psi^h_{hh}\brparenthese \bbrbrack \idx
	- \int\int_0^z (\log\rho)_{hh} \rho \widetilde v \,dz \cdot \nablah ( \dz v_p \cdot\psi^h_{hh}) \idx \\
	& \lesssim \int_0^1 \bigl( \hnorm{\zeta_{hh}}{\Lnorm2} + \hnorm{\zeta_h}{\Lnorm4}^2 \bigr) \bigl(\hnorm{\rho}{\Lnorm\infty} \hnorm{\nablah v}{\Lnorm4} + \hnorm{v}{\Lnorm\infty} \hnorm{\nablah \zeta}{\Lnorm4} \bigr) \,dz \\
	& ~~~~ \times \int_0^1 \hnorm{\dz v_p}{\Lnorm\infty} \hnorm{\psi^h_{hh}}{\Lnorm4}  \,dz
	+ \int_0^1 \bigl(\hnorm{\zeta_{hh}}{\Lnorm2} + \hnorm{\zeta_h}{\Lnorm4}^2\bigr) \hnorm{\rho}{\Lnorm\infty} \hnorm{v}{\Lnorm\infty} \,dz \\
	& ~~~~ \times \int_0^1 \bblparenthese \hnorm{\nablah \dz v_p}{\Lnorm4} \hnorm{\psi^h_{hh}}{\Lnorm4} + \hnorm{\dz v_p}{\Lnorm\infty}\hnorm{\nablah \psi^h_{hh}}{\Lnorm2} \bbrparenthese  \,dz \\
	& \lesssim \bigl(\norm{\nablah \zeta}{\Hnorm{1}} + \norm{\nablah \zeta}{\Hnorm{1}}^2 \bigr) \biggl(  (\norm{\rho}{\Lnorm\infty} \norm{\nabla v}{\Hnorm{1}} + \norm{v}{\Hnorm{2}} \norm{\nablah \zeta}{\Hnorm{1}} ) \\
	& ~~~~ \times \norm{\dz v_p}{\Lnorm\infty} + \norm{\rho}{\Lnorm\infty}\norm{v}{\Hnorm{2}} \norm{\nablah \dz v_p}{\Lnorm4} \biggr) \\
	& ~~~~ \times \bigl( \norm{\psi^h_{hh}}{\Lnorm2}^{1/2} \norm{\nablah \psi^h_{hh}}{\Lnorm2}^{1/2} + \norm{\psi^h_{hh}}{\Lnorm2} \bigr)
	 + \bigl(\norm{\nablah \zeta}{\Hnorm{1}} + \norm{\nablah \zeta}{\Hnorm{1}}^2 \bigr) \\
	 & ~~~~ \times \norm{\rho}{\Lnorm\infty}\norm{v}{\Hnorm{2}} \norm{\dz v_p}{\Lnorm\infty} \norm{\nablah \psi^h_{hh}}{\Lnorm2}
	\lesssim \delta \norm{\nabla\psi^h_{hh}}{\Lnorm2}^2 \\
	& ~~~~ + \delta \norm{\psi^h_{hh}}{\Lnorm{2}}^2 + C_\delta \bigl( Q(\mathcal E) + 1 + \mathfrak G_p\bigr) \mathfrak H_p \bigl( \varepsilon^4 + \varepsilon^2 \norm{\varepsilon^{-1} \xi}{\Hnorm{2}}^2 \bigr) .
\end{align*}
The estimates of the rest of $ I_8'' $ are listed below:
\begin{align*}
	& I_{8,1}'' \lesssim \int_0^1 \bblparenthese \hnorm{\zeta_{hh}}{\Lnorm2}\bigl( \hnorm{\nablah \psi^h}{\Lnorm4}+ \hnorm{v}{\Lnorm\infty} \hnorm{\nablah \zeta}{\Lnorm4} \bigr) + \hnorm{\zeta_h}{\Lnorm4}\bigl(\hnorm{\nablah \psi^h_h}{\Lnorm2} \\
	& ~~~~ + \hnorm{v_h}{\Lnorm4} \hnorm{\nablah\zeta}{\Lnorm4}
	 + \hnorm{v}{\Lnorm\infty}(\hnorm{\nablah^2\zeta}{\Lnorm2} + \hnorm{\nablah \zeta}{\Lnorm4}^2) \bigr)  + \hnorm{\rho}{\Lnorm4} \hnorm{\nablah \psi^h_{hh}}{\Lnorm2}\\
	 & ~~~~ + \hnorm{\rho}{\Lnorm\infty}\bigl(\hnorm{v_{hh}}{\Lnorm2} \hnorm{\nablah \zeta}{
	 \Lnorm4}
	+ \hnorm{ v_h}{\Lnorm4} ( \hnorm{\nablah^2\zeta}{\Lnorm2} + \hnorm{\nablah \zeta}{\Lnorm4}^2) \bigr)  \bbrparenthese \,dz \\
	& ~~~~ \times \int_0^1 \hnorm{\dz v_p}{\Lnorm\infty} \hnorm{\psi^h_{hh}}{\Lnorm4} \,dz  \lesssim \biggl( \norm{\nablah \zeta}{\Hnorm{1}}(\norm{\nabla \psi^h}{\Hnorm{1}} \\
	& ~~~~ + \norm{v}{\Hnorm2}(\norm{\nablah \zeta}{\Hnorm{1}}+ \norm{\nablah \zeta}{\Hnorm1}^2 )) + \norm{\rho}{\Hnorm{1}} \norm{\nablah \psi^h_{hh}}{\Lnorm2} \\
	& ~~~~ + \norm{\rho}{\Lnorm\infty} ( \norm{\nabla v}{\Hnorm1} \norm{\nablah \zeta}{\Hnorm1} + \norm{\nabla v}{\Hnorm1}(\norm{\nablah \zeta}{\Hnorm1} + \norm{\nablah \zeta}{\Hnorm1}^2))   \biggr) \\
	& ~~~~ \times \norm{\dz v_p}{\Lnorm\infty} \bigl( \norm{\psi_{hh}^h}{\Lnorm2}^{1/2} \norm{\nablah \psi^h_{hh}}{\Lnorm2}^{1/2} + \norm{\psi^h_{hh}}{\Lnorm2} \bigr) \lesssim \delta \norm{\nabla\psi^h_{hh}}{\Lnorm2}^2 \\
	& ~~~~ + \delta \norm{\psi^h_{hh}}{\Lnorm{2}}^2 + C_\delta \bigl( Q(\mathcal E) + 1 + \mathfrak G_p\bigr) \mathfrak H_p \bigl( \norm{\psi^h}{\Hnorm{2}}^2 + \varepsilon^4 + \varepsilon^2 \norm{\varepsilon^{-1} \xi}{\Hnorm{2}}^2 \bigr),\\
	& I_{8,3}''\lesssim \int_0^1 \bblparenthese \hnorm{\zeta_h}{\Lnorm4}\bigl(\hnorm{\nablah \psi^h}{\Lnorm4} + \hnorm{v}{\Lnorm\infty}\hnorm{\nablah \zeta}{\Lnorm4}\bigr) + \hnorm{\rho}{\Lnorm\infty} \bigl( \hnorm{\nablah^2 \psi^h}{\Lnorm2} \\
	& ~~~~ + \hnorm{v_h}{\Lnorm4} \hnorm{\nablah \zeta}{\Lnorm4}  + \hnorm{v}{\Lnorm\infty}(\hnorm{\nablah^2 \zeta}{\Lnorm2} + \hnorm{\nablah \zeta}{\Lnorm4}^2) \bigr) \bbrparenthese \,dz \\
	& ~~~~ \times \int_0^1 \hnorm{\dz v_{p,h}}{\Lnorm4} \hnorm{\psi^h_{hh}}{\Lnorm4} \,dz
	\lesssim \biggl( \norm{\nablah \zeta}{\Hnorm1} ( \norm{\nabla \psi^h}{\Hnorm1} \\
	& ~~~~ + \norm{v}{\Hnorm2} \norm{\nablah \zeta}{\Hnorm1})
	 + \norm{\rho}{\Lnorm\infty} \bigl\lbrack \norm{\nabla \psi^h}{\Hnorm1}
	+ \norm{v}{\Hnorm2} (\norm{\nablah \zeta}{\Hnorm1} \\
	& ~~~~ + \norm{\nablah \zeta}{\Hnorm1}^2 )\bigr\rbrack \biggr) \norm{v_{p,hz}}{\Lnorm4}
	 \bigl( \norm{\psi^h_{hh}}{\Lnorm2}^{1/2} \norm{\nablah \psi^h_{hh}}{\Lnorm2}^{1/2} + \norm{\psi^h_{hh}}{\Lnorm2} \bigr) \\
	& ~~~~ \lesssim \delta \norm{\nabla\psi^h_{hh}}{\Lnorm2}^2 + \delta \norm{\psi^h_{hh}}{\Lnorm{2}}^2 + C_\delta \blparenthese Q(\mathcal E) + 1 + \mathfrak G_p \brparenthese \mathfrak H_p \blparenthese \norm{\psi^h}{\Hnorm{2}}^2 \\
	& ~~~~ ~~~~ + \varepsilon^4 + \varepsilon^2 \norm{\varepsilon^{-1}\xi}{\Hnorm{2}}^2 \brparenthese,\\
	& I_{8,4}''\lesssim \int_0^1 \hnorm{\rho}{\Lnorm\infty} \bigl( \hnorm{\nablah \psi^h}{\Lnorm4} + \hnorm{v}{\Lnorm\infty} \hnorm{\nablah \zeta}{\Lnorm4} \bigr) \,dz \times \int_0^1 \hnorm{\dz v_{p,hh}}{\Lnorm{2}} \hnorm{\psi^h_{hh}}{\Lnorm{4}} \,dz\\
	& ~~ \lesssim \norm{\rho}{\Lnorm\infty} \bigl( \norm{\nabla \psi^h}{\Hnorm{1}} + \norm{v}{\Hnorm{2}}\norm{\nablah \zeta}{\Hnorm{1}} \bigr) \norm{v_{p,hhz}}{\Lnorm{2}} \bigl( \norm{\psi_{hh}^h}{\Lnorm{2}}^{1/2} \\
	& ~~~~ \times \norm{\nabla\psi_{hh}^h}{\Lnorm{2}}^{1/2} + \norm{\psi_{hh}^h}{\Lnorm{2}} \bigr) \lesssim \delta \norm{\nabla \psi^h_{hh}}{\Lnorm{2}}^2 + \delta \norm{\psi^h_{hh}}{\Lnorm{2}}^2 \\
	& ~~~~ + C_\delta \blparenthese Q(\mathcal E) + 1 + \mathfrak G_p \brparenthese \mathfrak H_p \blparenthese \norm{\psi^h}{\Hnorm{2}}^2 + \varepsilon^4 + \varepsilon^2 \norm{\varepsilon^{-1} \xi}{\Hnorm{2}}^2 \brparenthese.
\end{align*}
Hence, we have got
\begin{equation}\label{estimate:008}
	\begin{aligned}
		& I_8 \lesssim \delta \norm{\nabla\psi^h_{hh}}{2}^2 + \delta \norm{\psi^h_{hh}}{\Lnorm{2}}^2
		 + C_\delta \blparenthese Q(\mathcal E) + 1 + \mathfrak G_p \brparenthese \mathfrak H_p \\
		 & ~~~~ \times \blparenthese \norm{\psi^h}{\Hnorm{2}}^2 + \varepsilon^4 + \varepsilon^2 \norm{\varepsilon^{-1}\xi}{\Hnorm{2}}^2 \brparenthese.
	\end{aligned}
\end{equation}
Next, we establish the estimate of $ I_9 $. Notice, by making use of the facts that $ \xi_{hh}, \xi_{h}, \xi $ are independent of the $ z $-variable, and that $ \int_0^1 \dvh v_{p,hh} \,dz = \int_0^1 \dvh v_{p,h} \,dz = \int_0^1 \dvh v_p \,dz = 0 $, thanks to \eqref{bc-perturbed-eq}, one can derive
\begin{equation}\label{estimate:009}
\begin{aligned}
	& I_9 = - \rho_0^{-1} c^2_s \int \varepsilon^{-2} \xi_{hh} \dvh (\xi_{hh} v + 2 \xi_h v_h + \xi v_{hh}) \idx \\
	& ~~ = - \rho_0^{-1} c^2_s \int \bblparenthese \varepsilon^{-2} \xi_{hh} ( \dfrac{1}{2}\xi_{hh} \dvh \psi^h + 2 \xi_h \dvh \psi^h_h + 2 v_h \cdot\nablah \xi_h \\
	& ~~~~ + \xi \dvh \psi^h_{hh} + v_{hh} \cdot \nablah \xi ) \bbrparenthese \idx \lesssim \hnorm{\varepsilon^{-1}\nablah \xi}{\Hnorm{1}}^2 \int_0^1 \bblparenthese \hnorm{\nablah \psi^h}{\Lnorm\infty} \\
	& ~~~~ + \hnorm{v_{p,h}}{\Lnorm\infty} \bbrparenthese \,dz
	 + \hnorm{\varepsilon^{-1}\xi_{hh}}{\Lnorm2} \bblparenthese \hnorm{\varepsilon^{-1}\xi_h}{\Lnorm4} \int_0^1 \blparenthese \hnorm{\nablah^2 \psi^h}{\Lnorm4} \\
	 & ~~~~ + \hnorm{v_{p,hh}}{\Lnorm4} \brparenthese \,dz + \hnorm{\varepsilon^{-1}\xi}{\Lnorm\infty} \int_0^1 \hnorm{\nablah \psi^h_{hh}}{\Lnorm2} \,dz  \bbrparenthese
	\lesssim  \norm{\varepsilon^{-1}\nablah \xi}{\Hnorm{1}}^2 \\
	& ~~~~ \times \bigl( \norm{\nabla \psi^h}{\Hnorm{2}} + \norm{v_{p}}{\Hnorm{3}}  \bigr) + \norm{\varepsilon^{-1}\nablah \xi}{\Hnorm{1}} \norm{\varepsilon^{-1}\xi}{\Hnorm{2}} \\
	& ~~~~ \times \norm{\nablah \psi^h_{hh}}{\Lnorm{2}}
	 \lesssim \delta \norm{\nabla^3 \psi^h}{\Lnorm{2}}^2 +  C_\delta Q(\mathcal E) \norm{\varepsilon^{-1} \nablah \xi}{\Hnorm{1}}^2 \\
	 & ~~~~ + C_\delta Q(\mathcal E) \mathfrak H_p \norm{\varepsilon^{-1} \xi}{\Hnorm{2}}^2.
\end{aligned}
\end{equation}
Similarly, $ I_7, I_{10}, I_{12} $ can be estimated as follows:
\begin{align}
	& I_7 = \rho_0^{-1} \int \zeta_{hh} \bigl( \nablah (c^2_s \rho_1) - \mu \deltah v_p - \lambda \nablah \dvh v_p - \partial_{zz} v_p  \bigr) \cdot \psi^h_{hh}\idx \nonumber \\
	& ~~~~+ 2  \rho_0^{-1} \int \zeta_{h} \bigl( \nablah (c^2_s \rho_{1,h}) - \mu \deltah v_{p,h} - \lambda \nablah \dvh v_{p,h} - \partial_{zz} v_{p,h}  \bigr) \cdot \psi^h_{hh}\idx \nonumber \\
	& ~~~~ -  \rho_0^{-1} \int \bigl( \nablah (c^2_s \rho_{1,h}) - \mu \deltah v_{p,h} - \lambda \nablah \dvh v_{p,h} - \partial_{zz} v_{p,h}  \bigr) \cdot ( \zeta \psi^h_{hh})_h \idx  \nonumber \\
	& \lesssim \norm{\zeta_{hh}}{\Lnorm{2}} \norm{ \nablah (c^2_s \rho_1) - \mu \deltah v_p - \lambda \nablah \dvh v_p - \partial_{zz} v_p  }{\Lnorm{6}} \norm{\psi^h_{hh}}{\Lnorm{3}} \nonumber \\
	& ~~~~ + \norm{\nablah (c^2_s \rho_{1,h}) - \mu \deltah v_{p,h} - \lambda \nablah \dvh v_{p,h} - \partial_{zz} v_{p,h}}{\Lnorm{2}}  \nonumber \\
	& ~~~~ \times \bigl( \norm{\zeta_h}{\Lnorm{6}} \norm{\psi^h_{hh}}{\Lnorm{3}} + \norm{\zeta}{\Lnorm{\infty}}\norm{\psi^h_{hhh}}{\Lnorm{2}} \bigr)  \lesssim \delta \norm{\nabla \psi^h_{hh}}{\Lnorm{2}}^2 \nonumber  \\
	& ~~~~ + \delta \norm{\psi^h_{hh}}{\Lnorm{2}}^2
	 + C_\delta \blparenthese Q(\mathcal E) + 1  + \mathfrak G_p \brparenthese \mathfrak H_p \blparenthese \varepsilon^4 + \varepsilon^2 \norm{\varepsilon^{-1} \xi}{\Hnorm{2}}^2 \brparenthese, \label{estimate:010} \\	
	& I_{10} \nonumber
	= - \rho_0^{-1} c^2_s \int \bblparenthese \xi_{hh} \bigl( \dt \rho_{1,hh} + \rho_{1,hh} \dvh \psi^h + 2 \rho_{1,h} \dvh \psi^h_h + \rho_1 \dvh \psi^h_{hh} \\
	& \nonumber ~~~~ + v \cdot\nablah \rho_{1,hh} + 2 v_h \cdot\nablah \rho_{1,h} + v_{hh} \cdot\nablah \rho_1  \bigr) \bbrparenthese \idx \lesssim \norm{\nablah \xi}{\Hnorm{1}}\nonumber  \\
	& \nonumber ~~~~ \times \biggl( \norm{\dt \rho_{1,hh}}{\Lnorm2}
	  + \norm{\rho_{1,hh}}{\Lnorm3} \norm{\nablah\psi^h}{\Lnorm6} + \norm{\rho_{1,h}}{\Lnorm\infty} \norm{\nablah \psi^h_h}{\Lnorm2} \\
	 & \nonumber ~~~~ + \norm{\rho_1}{\Lnorm\infty} \norm{\nablah \psi^h_{hh}}{\Lnorm2}
	 + \norm{v}{H^2} \norm{\rho_1}{H^3} + \norm{v_h}{\Lnorm6} \norm{\nablah \rho_{1,h}}{\Lnorm3} \nonumber \\
	& ~~~~ + \norm{v_{hh}}{\Lnorm2} \norm{\nablah\rho_1}{\Lnorm\infty} \biggr) \nonumber
	\lesssim \delta \norm{\nablah \psi^h_{hh}}{\Lnorm2}^2 + \delta \norm{\varepsilon^{-1} \nablah \xi}{\Hnorm{1}}^2 \\
	& ~~~~ + C_\delta  \blparenthese 1 + \mathfrak G_p \brparenthese \mathfrak H_p \blparenthese \varepsilon^2 + \varepsilon^2 \norm{\varepsilon^{-1} \xi}{\Hnorm{2}}^2 + \varepsilon^2 \norm{\psi^h}{\Hnorm{2}}^2 \brparenthese
	, \label{estimate:011} \\
	& I_{12} = \int \varepsilon^{-2} \mathcal R_{hh} \dvh \psi^h_{hh} \idx \lesssim \varepsilon^{-2} \bigl( \norm{\zeta_{hh}}{\Lnorm2} \norm{\zeta}{\Lnorm\infty} + \norm{\zeta_h}{\Lnorm4}^2 \bigr) \nonumber \\
	& ~~~~ \times \norm{\nablah\psi_{hh}^h}{\Lnorm2} \nonumber
	\lesssim \delta \norm{\nablah \psi^h_{hh}}{\Lnorm2}^2 + C_\delta Q(\mathcal E)  \norm{\varepsilon^{-1} \nablah \xi}{\Hnorm{1}}^2 \nonumber \\
	& ~~~~ + C_\delta \blparenthese Q(\mathcal E) + \mathfrak G_p \brparenthese \mathfrak H_p \blparenthese  \varepsilon^4 + \varepsilon^2 \norm{\varepsilon^{-1}\xi}{\Hnorm{2}}^2 \brparenthese. \label{estimate:012} 
\end{align}
Here we have applied \eqref{residue-est-h} and that
\begin{equation}\label{residue-est-hh}
	\abs{\mathcal R_{hh}}{} = \abs{\gamma(\gamma-1) \int_{\rho_0}^\rho \rho_{hh} y^{\gamma-2} \,dy + \gamma(\gamma-1) \rho_h^2 \rho^{\gamma-2}}{} \lesssim \abs{\zeta_{hh}}{} \abs{\zeta}{} + \abs{\zeta_h}{2}.
\end{equation}
Therefore, \eqref{ee:horizontal-derivative} follows from \eqref{ee:003} and inequalities \eqref{estimate:007}, \eqref{estimate:008}, \eqref{estimate:009}, \eqref{estimate:010}, \eqref{estimate:011}, \eqref{estimate:012}, above.
\end{proof}

Next, we will derive the required estimate of $ \xi_{hh} $.
After integrating \subeqref{eq:perturbation}{2} over $ z \in (0,1) $, we have the following equation, thanks to \eqref{bc-perturbed-eq}:
\begin{equation}\label{eq:average-001}
	\begin{aligned}
		& \rho \dt  \overline{\psi^h} +  \int_0^1 \bigl( \rho v \cdot\nablah \psi^h - \rho w_z \psi^h \bigr) \,dz + \nablah (\varepsilon^{-2} c^2_s \xi) \\
		& ~~~~ ~~~~ = \mu \deltah \overline{\psi^h} + \lambda \nablah \dvh \overline{\psi^h} - \nablah (\varepsilon^{-2}\mathcal R) + \int_0^1 \blparenthese \mathcal F_1 + \mathcal F_2 \brparenthese \,dz.
	\end{aligned}
\end{equation}
After applying $ \partial_h $ to \eqref{eq:average-001}, one has
\begin{equation}\label{eq:average-002}
	\begin{aligned}
		& \varepsilon^{-2} c^2_s \nablah \xi_h = \underbrace{ - \zeta_h \dt \overline{\psi^h} - \rho \dt \overline{\psi^h_h}}_{R_1} \\
		& ~~~~ \underbrace{- \int_0^1 \bblparenthese \zeta_h v \cdot \nablah \psi^h + \rho v_h \cdot\nablah \psi^h + \rho v \cdot\nablah \psi^h_h \bbrparenthese \,dz}_{R_2}\\
		& ~~~~ + \underbrace{\int_0^1 \bblparenthese \zeta_h w_z \psi^h + \rho w_{hz} \psi^h + \rho w_z \psi^h_h \bbrparenthese \,dz}_{R_3} \\
		& ~~~~ + \underbrace{ \mu \deltah \overline{\psi^h_h} + \lambda \nablah \dvh \overline{\psi^h_h}  - \nablah (\varepsilon^{-2} \mathcal R_{h})}_{R_4}
		 + \underbrace{\int_0^1 \bblparenthese \mathcal F_{1,h} + \mathcal F_{2,h} \bbrparenthese \,dz}_{R_{5}}.
	\end{aligned}
\end{equation}
What we need is to estimate the $ L^2 $-norm of the terms on the right-hand side of \eqref{eq:average-002}. In fact, after applying the Minkowski, H\"older and Sobolev embedding inequalities, one has
\begin{align*}
	& \hnorm{R_1}{\Lnorm2} \lesssim \int_0^1 \hnorm{\zeta_h}{\Lnorm4} \hnorm{\dt\psi^h}{\Lnorm4} \,dz + \int_0^1 \hnorm{\rho}{\Lnorm\infty}\hnorm{\dt \psi^h_h}{\Lnorm2} \,dz \lesssim \norm{\zeta_h}{\Hnorm{1}} \\
	& ~~~~ \times \bigl( \norm{\dt \psi^h}{\Lnorm{2}}^{1/2}
	 \norm{\nabla \dt \psi^h}{\Lnorm{2}}^{1/2} + \norm{\dt \psi^h}{\Lnorm{2}} \bigr) + \norm{\rho}{\Lnorm{\infty}} \norm{\dt \psi^h_h}{\Lnorm{2}} \\
	 & ~~~~ \lesssim \norm{\nabla \dt \psi^h}{\Lnorm{2}}
	+  \bigl( \norm{\varepsilon^{-1}\xi}{\Hnorm{2}}+ \mathfrak G_p \norm{\varepsilon \rho_1}{\Hnorm{2}} \bigr) \norm{\varepsilon \dt \psi^h}{\Lnorm{2}}, \\
	& \hnorm{R_2}{\Lnorm2} \lesssim \int_0^1 \biggl( \hnorm{\zeta_h}{\Lnorm4} \hnorm{v}{\Lnorm\infty} \hnorm{\nablah \psi^h}{\Lnorm4} + \hnorm{\rho}{\Lnorm\infty}\hnorm{v_h}{\Lnorm4}\hnorm{\nablah \psi^h}{\Lnorm4} + \hnorm{\rho}{\Lnorm\infty}\\
	& ~~~~ \times \hnorm{v}{\Lnorm\infty}
	 \hnorm{\nablah \psi^h_h}{\Lnorm2} \biggr)  \,dz
	\lesssim \bigl( \norm{\zeta}{\Hnorm{2}} + \norm{\rho}{\Lnorm\infty} \bigr)\norm{v}{\Hnorm{2}}\norm{\nablah \psi^h}{\Hnorm{1}}, \\
	& \hnorm{R_4}{\Lnorm2} \lesssim \int_0^1 \hnorm{\nablah \psi^h}{\Hnorm{2}}\,dz + \varepsilon^{-2}\bigl(\hnorm{\zeta_{hh}}{\Lnorm2}\hnorm{\zeta}{\Lnorm\infty} + \hnorm{\zeta_h}{\Lnorm4}^2 \bigr) \lesssim \norm{\nablah \psi^h}{\Hnorm{2}}\\
	& ~~~~ + \varepsilon^2 \norm{\rho_1}{\Hnorm{2}}^2 +  \norm{\varepsilon^{-1} \xi}{\Hnorm{2}}\norm{\varepsilon^{-1} \nablah \xi}{\Hnorm{1}}.
\end{align*}
On the other hand, after substituting \eqref{id:vertical_velocity-200} and \eqref{id:vertical_velocity-110}, below, we have
\begin{align*}
	& R_3 = \int_0^1 \bblparenthese \bigl\lbrack w_{p,z} - ( \dvh \widetilde{\psi^h} + \widetilde v \cdot\nablah \log\rho) \bigr\rbrack \bigl( \zeta_h \psi^h + \rho \psi^h_h\bigr) + \rho \psi^h \bigl\lbrack w_{p,hz} \\
	& ~~~~ ~~~~ - ( \dvh \widetilde{\psi^h_h} + \widetilde v_h \cdot\nablah \log \rho + \widetilde v \cdot\nablah (\log\rho)_h ) \bigr\rbrack \bbrparenthese \,dz.
\end{align*}
Therefore, one has
\begin{align*}
	& \hnorm{R_3}{\Lnorm2} \lesssim \int_0^1 \bblparenthese  \bigl( \hnorm{w_{p,z}}{\Lnorm4} + \hnorm{\nablah \psi^h}{\Lnorm4} + \hnorm{v}{\Lnorm\infty}\hnorm{\nablah\zeta}{\Lnorm4} \bigr) \\
	& ~~~~ \times \bigl( \hnorm{\zeta_h}{\Lnorm4} \hnorm{\psi^h}{\Lnorm\infty} + \hnorm{\rho}{\Lnorm\infty} \hnorm{\psi^h_h}{\Lnorm4} \bigr) + \hnorm{\rho}{\Lnorm\infty} \hnorm{\psi^h}{\Lnorm\infty} \\
	& ~~~~ \times \bigl( \hnorm{w_{p,hz}}{\Lnorm2} + \hnorm{\nablah \psi^h_h}{\Lnorm2} + \hnorm{v_h}{\Lnorm4} \hnorm{\nablah \zeta}{\Lnorm4} \\
	& ~~~~ + \hnorm{v}{\Lnorm\infty} ( \hnorm{\zeta_{hh}}{\Lnorm2} + \hnorm{\zeta_h}{\Lnorm4}^2) \bigr) \bbrparenthese \,dz \\
	& ~~~~ \lesssim Q(\mathcal E) \bigl( \norm{\nablah \psi^h}{\Hnorm{1}} + \varepsilon \norm{\varepsilon^{-1} \nablah \xi}{\Hnorm{1}} \bigr) + \blparenthese Q(\mathcal E) + 1 + \mathfrak G_p \brparenthese \mathfrak H_p^{1/2}.
\end{align*}
After substituting \eqref{nonlinearities}, we obtain
\begin{align*}
	& R_5 = \int_0^1 \bblparenthese \zeta_h Q_p + \zeta Q_{p,h} - \zeta_h \psi^h \cdot\nablah v_p - \rho\psi^h_h \cdot\nablah v_p - \rho \psi^h \cdot \nablah v_{p,h} \\
	& ~~~~ + \psi^z_z (\zeta_h v_p + \rho v_{p,h} ) + \rho v_p \psi^z_{hz} \bbrparenthese \,dz = \int_0^1 \bblparenthese \zeta_h Q_p + \zeta Q_{p,h} - \zeta_h \psi^h \cdot\nablah v_p \\
	& ~~~~ - \rho\psi^h_h \cdot\nablah v_p
	 - \rho \psi^h \cdot \nablah v_{p,h}  - \bigl( \dvh\widetilde{\psi^h} + \widetilde v \cdot\nablah \log\rho \bigr) \bigl(\zeta_h v_p + \rho v_{p,h} \bigr) \\
	& ~~~~ - \bigl( \dvh\widetilde{\psi^h_h} + \widetilde{v}_h \cdot\nablah \log\rho + \widetilde v \cdot\nablah (\log\rho)_h \bigr)\rho v_p \bbrparenthese \,dz.
\end{align*}
Hence, applying the Minkowski, H\"older and Sobolev embedding  inequalities implies
\begin{align*}
	& \hnorm{R_5}{\Lnorm 2} \lesssim \int_0^1 \bblparenthese \hnorm{\zeta_h}{\Lnorm 4} \hnorm{Q_p}{\Lnorm 4} + \hnorm{\zeta}{\Lnorm\infty} \hnorm{Q_{p,h}}{\Lnorm2} + \hnorm{\zeta_h}{\Lnorm4}\hnorm{\psi^h}{\Lnorm\infty}\hnorm{\nablah v_p}{\Lnorm4} \\
	& ~~~~ + \hnorm{\rho}{\Lnorm\infty}\hnorm{\psi^h_h}{\Lnorm4} \hnorm{\nablah v_p}{\Lnorm4}
	+ \hnorm{\rho}{\Lnorm\infty}\hnorm{\psi^h}{\Lnorm\infty} \hnorm{\nablah v_{p,h}}{\Lnorm2} \\
	& ~~~~ + \bigl( \hnorm{\nablah \psi^h}{\Lnorm4} + \hnorm{v}{\Lnorm\infty}\hnorm{\nablah \zeta}{\Lnorm4} \bigr) \bigl( \hnorm{\zeta_h}{\Lnorm4}\hnorm{v_p}{\Lnorm\infty} + \hnorm{\rho}{\Lnorm\infty} \hnorm{v_{p,h}}{\Lnorm4} \bigr) \\
	& ~~~~ + \hnorm{\rho}{\Lnorm\infty} \hnorm{v_p}{\Lnorm\infty} \bigl\lbrack \hnorm{\nablah \psi^h_h}{\Lnorm2} + \hnorm{v_h}{\Lnorm4} \hnorm{\nablah \zeta}{\Lnorm4} \\
	& ~~~~ ~~~~ + \hnorm{v}{\Lnorm\infty} ( \hnorm{\zeta_{hh}}{\Lnorm2} + \hnorm{\zeta_h}{\Lnorm4}^2) \bigr\rbrack \bbrparenthese \,dz \\
	& ~~~~ \lesssim  \blparenthese Q(\mathcal E) + 1 + \mathfrak G_p \brparenthese \mathfrak H_p^{1/2} \blparenthese \norm{\psi^h}{\Hnorm{2}} + \varepsilon \brparenthese.
\end{align*}
Summing up these estimates, we get the following inequality from \eqref{eq:average-002},
\begin{equation}\label{ee:007}
\begin{aligned}
	&  \norm{\varepsilon^{-1}\xi_{hh}}{\Lnorm2}^2 \lesssim \varepsilon^2 \hnorm{\varepsilon^{-2}\xi_{hh}}{\Lnorm2}^2 \lesssim \sum_{i=1}^5 \varepsilon^2\hnorm{R_i}{\Lnorm2}^2 \lesssim \norm{\varepsilon \nabla \dt \psi^h}{\Lnorm{2}}^2 \\
	& ~~~~ +  \varepsilon^2 \norm{ \nabla \psi^h}{\Hnorm{2}}^2
	 + \varepsilon^2 Q(\mathcal E) \blparenthese \norm{\varepsilon^{-1}\nablah \xi}{\Hnorm{1}}^2 + \norm{\nabla \psi^h}{\Hnorm{1}}^2 \\
	& ~~~~ + \norm{\varepsilon \dt \psi^h}{\Lnorm{2}}^2 \brparenthese
	 +  \varepsilon^2 \blparenthese Q(\mathcal E) + 1 + \mathfrak G_p \brparenthese \mathfrak H_p.
\end{aligned}
\end{equation}
We summarize the result in the following:
\begin{lemma}\label{lm:horizontal-diffusion}
Under the same assumptions as in Lemma \ref{lm:temporal-derivative}, the following holds:
	\begin{equation}\label{ee:dhh-xi}
	\begin{aligned}
	&  \norm{\varepsilon^{-1}\nablah \xi}{\Hnorm 1}^2 \leq C \norm{\varepsilon \nabla \dt \psi^h}{\Lnorm{2}}^2 +  \varepsilon^2 \norm{ \nabla \psi^h}{\Hnorm{2}}^2 \\
	& ~~~~ + \varepsilon^2 C Q(\mathcal E) \blparenthese \norm{\varepsilon^{-1}\nablah \xi}{\Hnorm{1}}^2 + \norm{\nabla \psi^h}{\Hnorm{1}}^2 + \norm{\varepsilon \psi^h_t}{\Lnorm{2}}^2 \brparenthese \\
	& ~~~~ + C \varepsilon^2 \blparenthese Q(\mathcal E) + 1 + \mathfrak G_p \brparenthese \mathfrak H_p,
	\end{aligned}
	\end{equation}
	for some positive constant $ C $ independent of $ \varepsilon $.
\end{lemma}
\begin{proof}
	This is the direct consequence of \eqref{ee:007} and the Poincar\'e inequality.
\end{proof}


\subsection{Vertical derivatives estimates}
Now we turn to the required estimates of vertical derivatives. To do so, we first apply $ \dz  $ to system \eqref{eq:perturbation} and write down the resultant system as follows:
\begin{equation}\label{eq:perturbation-100}
	\begin{cases}
		\rho_0(\dvh \psi^h_z + \dz \psi^h_z) = \mathcal G_{1,z} + \mathcal G_{2,z} & \text{in} ~ \Omega, \\
		\rho \dt \psi^h_z + \rho v\cdot\nablah \psi^h_z + \rho w \dz \psi^h_z = \mu \deltah \psi^h_z + \lambda \nablah \dvh \psi^h_z \\
		~~~~ ~~~~ + \partial_{zz} \psi^h_z
		 + \mathcal F_{1,z} + \mathcal F_{2,z} + \mathcal H_{z} & \text{in} ~ \Omega,
	\end{cases}
\end{equation}
where
\begin{equation}\label{nonlinear-dz}
	\mathcal H_z := - \rho v_z \cdot\nablah \psi^h  - \rho w_z \dz \psi^h.
\end{equation}
Then we apply $ \dz $ to system \eqref{eq:perturbation-100} again and obtain the following system:
\begin{equation}\label{eq:perturbation-200}
	\begin{cases}
		\rho_0( \dvh\psi^h_{zz} + \dz \psi^h_{zz} ) = \mathcal G_{1,zz} + \mathcal G_{2,zz} &\text{in} ~ \Omega, \\
		\rho \dt \psi^h_{zz} + \rho v\cdot\nablah \psi^h_{zz} + \rho w \dz \psi^h_{zz} = \mu \deltah \psi^h_{zz} + \lambda\nablah \dvh \psi^h_{zz} \\
		~~~~ + \partial_{zz} \psi^h_{zz}
		 + \mathcal F_{1,zz} + \mathcal F_{2,zz} + \mathcal H_{zz} & \text{in} ~ \Omega,
	\end{cases}
\end{equation}
where
\begin{equation}\label{nonlinear-dzz}
	\mathcal H_{zz} := - \rho v_{zz} \cdot\nablah \psi^h - 2\rho v_z\cdot\nablah \psi^h_z - \rho w_{zz} \dz \psi^h - 2 \rho w_z \dz \psi^h_z.
\end{equation}
Notice, here we have employed the fact that $ \rho, \xi, \rho_1, \mathcal R $ are independent of the $ z $ variable. Also \subeqref{eq:perturbation-200}{2} is a parabolic equation of $ \psi^h_{zz} $. Now we perform standard $ L^2 $ estimate on system \eqref{eq:perturbation-200}.
\begin{lemma}\label{lm:vertical-derivative} Under the same assumptions as in Lemma \ref{lm:temporal-derivative}, we have
	\begin{equation}\label{ee:vertical-est-dzz-psi}
	\begin{aligned}
		& \dfrac{d}{dt}\norm{\rho^{1/2}\psi^h_{zz}}{\Lnorm 2}^2 +  C_{\mu,\lambda} \norm{\nabla \psi^h_{zz}}{\Lnorm 2}^2 \leq \delta \norm{\nabla\psi^h_{zz}}{\Lnorm{2}}^2 + \delta \norm{\psi^h_{zz}}{\Lnorm{2}}^2 \\
		& ~~~~ + C_\delta Q(\mathcal E) \norm{\nabla\psi^h}{\Hnorm{1}}^2 + C_\delta \blparenthese Q(\mathcal E) + 1 + \mathfrak G_p \brparenthese \mathfrak H_p \blparenthese \norm{\psi^h}{\Hnorm{2}}^2 + \varepsilon^2 \brparenthese,
	\end{aligned}
\end{equation}
for some positive constant $ C_{\mu,\lambda} $ independent of $ \varepsilon $.
\end{lemma}
\begin{proof}
After taking the $ L^2 $-inner product of \subeqref{eq:perturbation-200}{2} with $ \psi^h_{zz} $, we have the following:
\begin{equation}\label{ee:004}
\begin{aligned}
	& \dfrac{d}{dt} \biggl\lbrace \dfrac 1 2 \norm{\rho^{1/2}\psi^h_{zz}}{\Lnorm2}^2 \biggr\rbrace + \mu \norm{\nablah \psi^h_{zz}}{\Lnorm2}^2 + \lambda \norm{\dvh \psi^h_{zz}}{\Lnorm2}^2 \\
	& ~~~~ + \norm{\dz \psi^h_{zz}}{\Lnorm2}^2
	 = \int \mathcal F_{1,zz} \cdot \psi^h_{zz} \idx + \int \mathcal F_{2,zz} \cdot \psi^h_{zz} \idx \\
	& ~~~~ + \int \mathcal H_{zz} \cdot \psi^h_{zz} \idx
	 =: I_{13} + I_{14} + I_{15}.
\end{aligned}
\end{equation}
Again, we shall estimate the terms on the right-hand side of \eqref{ee:004}. We  begin with the term $ I_{15} $. Notice first, after taking $ \partial_z, \partial_{zz} $ to \eqref{id:vertical_velocity}, we have the following identities:
\begin{equation}\label{id:vertical_velocity-200}
	\begin{aligned}
		w_z & = w_{p,z} - \bigl( \dvh \widetilde{\psi^h} + \widetilde v\cdot \nablah \log \rho \bigr) = w_{p,z} + \psi^z_z, \\
		w_{zz} & = w_{p,zz} - \bigl( \dvh \widetilde{\psi^h_z} + \widetilde v_z \cdot\nablah \log \rho \bigr) = w_{p,zz} + \psi^z_{zz}.
	\end{aligned}
\end{equation}
Consequently, after substituting \eqref{id:vertical_velocity-200} in $ I_{15} $, it can be estimated as follows:
\begin{align*}
	& I_{15} = - \int \bblparenthese \rho \bigl( v_{zz} \cdot \nablah \psi^h + 2 v_z \cdot\nablah \psi^h_z + w_{p,zz}\dz \psi^h \\
	& ~~~~ - \dvh \widetilde{\psi^h_z}\dz \psi^h
	 - \widetilde v_z \cdot \nablah \log \rho \dz \psi^h
	 + 2 w_{p,z} \dz \psi^h_z \\
	 & ~~~~ - 2 \dvh\widetilde{\psi^h} \dz \psi^h_z
	  - 2 \widetilde v \cdot \nablah \log \rho \dz \psi^h_z  \bigr) \cdot \psi^h_{zz} \bbrparenthese \idx \\
	 & ~~
	\lesssim \norm{\rho}{\Lnorm\infty} \norm{\psi^h_{zz}}{\Lnorm3} \biggl( \norm{v_{zz}}{\Lnorm2} \norm{\nablah \psi^h}{\Lnorm6}
	 + \norm{v_z}{\Lnorm6} \norm{\nablah \psi^h_z}{\Lnorm2} \\
	& ~~~~ + \norm{w_{p,zz}}{\Lnorm2} \norm{\dz \psi^h}{\Lnorm6}
	 + \norm{\nablah \psi^h_z}{\Lnorm2} \norm{\dz \psi^h}{\Lnorm6} \\
	 & ~~~~ + \norm{\nablah \zeta}{\Lnorm3}\norm{v_z}{\Lnorm6} \norm{\dz \psi^h}{\Lnorm6} + \norm{w_{p,z}}{\Lnorm6} \norm{\dz \psi^h_z}{\Lnorm2} \\
	& ~~~~ + \norm{\nablah \psi^h}{\Lnorm6} \norm{\dz \psi^h_z }{\Lnorm2} + \norm{\nablah \zeta}{\Lnorm6} \norm{v}{\Lnorm\infty} \norm{\dz \psi^h_z}{\Lnorm2} \biggr) \\
	& ~~ \lesssim \norm{\rho}{\Lnorm\infty}\bigl( \norm{\psi^h_{zz}}{\Lnorm2}^{1/2} \norm{\nabla \psi^h_{zz}}{\Lnorm2}^{1/2} + \norm{\psi^h_{zz}}{\Lnorm2} \bigr) \\
	& ~~~~ \times \bigl(\norm{v_p}{\Hnorm{2}} + \norm{\nabla \psi^h}{\Hnorm{1}} + \norm{\dz w_{p}}{\Hnorm{1}}
	 + \norm{ \nablah \zeta}{\Hnorm{1}}\norm{v}{\Hnorm{2}} \bigr)\\
	 & ~~~~ \times \norm{\nabla \psi^h}{\Hnorm{1}} \lesssim \delta \norm{\nabla\psi^h_{zz}}{\Lnorm2}^2 +\delta \norm{\psi^h_{zz}}{\Lnorm{2}}^2 \\
	 & ~~~~ + C_\delta  Q(\mathcal E)  \norm{\nabla \psi^h}{\Hnorm{1}}^2
	  + C_\delta \blparenthese Q(\mathcal E) + 1 + \mathfrak G_p \brparenthese \mathfrak H_p \norm{\psi^h}{\Hnorm{2}}^2 .
\end{align*}
The estimates of the terms $ I_{13}, I_{14} $ are listed below:
\begin{align*}
	& I_{14} = - \int \bblparenthese \rho \bigl( \psi^h_{zz} \cdot
	\nablah v_p + 2 \psi^h_z \cdot\nablah v_{p,z} + \psi^h \cdot\nablah v_{p,zz} + \psi^z_{zz} \dz v_p \\
	& ~~~~ + 2 \psi^z_{z} \dz v_{p,z}
	+ \psi^z \dz v_{p,zz} \bigr) \cdot \psi^h_{zz} \bbrparenthese \idx = - \int \bblparenthese \rho \bigl( \psi^h_{zz} \cdot\nablah v_p \\
	& ~~~~ + 2 \psi^h_z \cdot\nablah v_{p,z}
	 + \psi^h \cdot\nablah v_{p,zz}
	 - \dvh \widetilde{\psi^h_z} \dz v_p - \widetilde v_z \cdot \nablah \log \rho \dz v_p \\
	 & ~~~~ - 2 \dvh \widetilde{\psi^h} \dz v_{p,z}
	  - 2 \widetilde v \cdot\nablah \log\rho \dz v_{p,z} \bigr) \cdot \psi^h_{zz} \bbrparenthese \idx \\
	& ~~~~ + \int \bblbrack \int_0^z \bblparenthese \rho\dvh\widetilde{\psi^h} + \widetilde v \cdot\nablah \rho \bbrparenthese \,dz
	 \times \blparenthese \dz v_{p,zz} \cdot \psi^h_{zz} \brparenthese \bbrbrack  \idx \\
	 & ~~ \lesssim \norm{\rho}{\Lnorm\infty} \norm{\psi^h_{zz}}{\Lnorm3}
	\bigl( \norm{\psi^h_{zz}}{\Lnorm2} \norm{\nablah v_p}{\Lnorm6}
	 + \norm{\psi^h_z}{\Lnorm6}\norm{\nablah v_{p,z}}{\Lnorm2}\\
	 & ~~~~ + \norm{\psi^h}{\Lnorm6} \norm{\nablah v_{p,zz}}{\Lnorm2} + \norm{\nablah \psi^h_z}{\Lnorm2} \norm{\dz v_p}{\Lnorm6} \\
	& ~~~~ + \norm{v_z}{\Lnorm6} \norm{\nablah \zeta}{\Lnorm3} \norm{\dz v_p}{\Lnorm6} + \norm{\nablah \psi^h}{\Lnorm6} \norm{v_{p,zz}}{\Lnorm2} \\
	& ~~~~ + \norm{v}{\Lnorm\infty} \norm{\nablah \zeta}{\Lnorm6} \norm{v_{p,zz}}{\Lnorm2} \bigr)
	 + \int_0^1 \bblparenthese \hnorm{\rho}{\Lnorm\infty} \hnorm{\nablah \psi^h}{\Lnorm4} \\
	 & ~~~~ + \hnorm{v}{\Lnorm\infty} \hnorm{\nablah\zeta}{\Lnorm4} \bbrparenthese \,dz
	  \times \int_0^1 \hnorm{\dz v_{p,zz}}{\Lnorm2} \hnorm{\psi^h_{zz}}{\Lnorm4} \,dz \\
	& ~~ \lesssim  \bigl( \norm{\psi^h_{zz}}{\Lnorm2}^{1/2} \norm{\nabla\psi^h_{zz}}{\Lnorm2}^{1/2} + \norm{\psi^h_{zz}}{\Lnorm2} \bigr) \norm{\rho}{\Lnorm\infty}
	 \bigl( \norm{\nabla \psi^h}{\Hnorm{1}}\norm{v_p}{\Hnorm{2}} \\
	& ~~~~ + \norm{\psi^h}{\Hnorm{1}} \norm{v_p}{\Hnorm{3}} + \norm{\nablah \zeta}{\Hnorm{1}} \norm{v}{\Hnorm{2}} \norm{v_p}{\Hnorm{2}} \bigr)\\
	& ~~~~ + \bigl( \norm{\psi^h_{zz}}{\Lnorm2}^{1/2} \norm{\nabla\psi^h_{zz}}{\Lnorm2}^{1/2} + \norm{\psi^h_{zz}}{\Lnorm2} \bigr) \bigl( \norm{\rho}{\Lnorm\infty} \norm{\nablah \psi^h}{\Hnorm{1}} \\
	& ~~~~ + \norm{v}{\Hnorm{2}} \norm{\nablah \zeta}{\Hnorm{1}} \bigr) \norm{v_p}{\Hnorm{3}}\lesssim \delta \norm{\nabla \psi^h_{zz}}{\Lnorm2}^2 + \delta \norm{\psi^h_{zz}}{\Lnorm{2}}^2 \\
	& ~~~~ + C_\delta \blparenthese Q(\mathcal E)+1 + \mathfrak G_p \brparenthese \mathfrak H_p \blparenthese \norm{\psi^h}{\Hnorm{2}}^2 + \varepsilon^2 \norm{\varepsilon^{-1}\xi}{\Hnorm{2}}^2 + \varepsilon^4 \brparenthese, \\
	& I_{13} = \int \zeta Q_{p,zz} \cdot \psi^h_{zz} \idx = - \int \zeta Q_{p,z} \cdot \psi^h_{zzz} \idx \lesssim  \delta \norm{\psi^h_{zzz}}{\Lnorm{2}}^2 \\
	& ~~~~ + C_\delta \blparenthese Q(\mathcal E) + \mathfrak G_p \brparenthese \mathfrak H_p \varepsilon^2.
\end{align*}
After summing the estimates for  $ I_{13}, I_{14}, I_{15} $, above, and \eqref{ee:004}, we conclude \eqref{ee:vertical-est-dzz-psi}.
\end{proof}

\subsection{Mixed horizontal and vertical derivatives estimates}

What is left is to estimate the $ L^2 $ norm of $ \partial_{hz} \psi^h $. We apply $ \partial_h $ to \eqref{eq:perturbation-100} and write down the resultant system:
\begin{equation}\label{eq:perturbation-110}
	\begin{cases}
		\rho_0 ( \dvh \psi^h_{hz} + \dz \psi^h_{hz} ) = \mathcal G_{1,hz} + \mathcal G_{2,hz} & \text{in} ~ \Omega, \\
		\rho\dt \psi^h_{hz} + \rho v\cdot\nablah \psi^h_{hz} + \rho w \dz \psi^h_{hz} = \mu \deltah \psi^h_{hz} \\
		~~~~ + \lambda \nablah \dvh \psi^h_{hz}
		 + \partial_{zz} \psi^h_{hz}
		+ \mathcal F_{1,hz} + \mathcal F_{2,hz} + \mathcal H_{hz} &\text{in}~ \Omega,
	\end{cases}
\end{equation}
where
\begin{equation}\label{nonlinear-dhz}
\begin{aligned}
	& \mathcal H_{hz} : = -\zeta_h \dt \psi^h_z - \zeta_h v_z \cdot\nablah \psi^h - \rho (v_{z} \cdot \nablah \psi^h)_h - \zeta_h w_z \dz \psi^h \\
	& ~~~~ - \rho (w_z \dz \psi^h)_h
	 - \zeta_h v\cdot
	\nablah \psi^h_z - \rho v_h \cdot\nablah \psi^h_z - \zeta_h w\dz \psi^h_z \\
	& ~~~~ - \rho w_h \dz \psi^h_z .
\end{aligned}
\end{equation}

\begin{lemma}\label{lm:horizontal-vertical-derivative} Under the same assumptions as in Lemma \ref{lm:temporal-derivative}, we have
	\begin{equation}\label{ee:vertical-est-dhz-psi}
	\begin{aligned}
		& \dfrac{d}{dt} \norm{\rho^{1/2} \psi^h_{hz}}{\Lnorm 2}^2  + C_{\mu,\lambda} \norm{\nabla \psi^h_{hz}}{\Lnorm 2}^2 \leq \delta \blparenthese \norm{\nabla^3 \psi^h}{\Lnorm2}^2 + \norm{\psi^h_{hz}}{\Lnorm{2}}^2 \\
		& ~~~~ +  \norm{\varepsilon\nabla \dt \psi^h}{\Lnorm2}^2 \brparenthese
		 + C_\delta Q(\mathcal E) \norm{\nabla \psi^h}{\Hnorm{1}}^2 \\
		& ~~~~ + C_\delta \blparenthese Q(\mathcal E) + 1 + \mathfrak G_p \brparenthese \mathfrak H_p \bigl( \varepsilon^2 + \norm{\psi^h}{\Hnorm{2}}^2 \bigr),
	\end{aligned}
\end{equation}
for some positive constant $ C_{\mu,\lambda} $, which is independent of $ \varepsilon $.
\end{lemma}

\begin{proof}	
Take the $ L^2 $ inner produce of \subeqref{eq:perturbation-110}{2} with $ \psi^h_{hz} $. It follows,
\begin{equation}\label{ee:005}
	\begin{aligned}
		& \dfrac{d}{dt} \biggl\lbrace \dfrac 1 2 \norm{\rho^{1/2} \psi^h_{hz}}{\Lnorm2}^2  \biggr\rbrace + \mu \norm{\nablah \psi^h_{hz}}{\Lnorm2}^2 + \lambda \norm{\dvh \psi^h_{hz}}{\Lnorm2}^2 \\
		& ~~~~ + \norm{\dz \psi^h_{hz}}{\Lnorm2}^2
		= \int \mathcal F_{1,hz} \cdot \psi^h_{hz} \idx + \int \mathcal F_{2,hz} \cdot \psi^h_{hz} \idx \\
		& ~~~~ + \int \mathcal H_{hz} \cdot \psi^h_{hz} \idx =: I_{16} + I_{17} + I_{18}.
	\end{aligned}
\end{equation}
As before, we shall estimate the terms on the right-hand side of \eqref{ee:005}. Notice that, we have the following identity after taking $ \partial_h $ to \subeqref{id:vertical_velocity-200}{1}:
\begin{equation}\label{id:vertical_velocity-110}
	w_{hz} = w_{p,hz} - \bigl( \dvh \widetilde{\psi^h_h} + \widetilde v_h \cdot \nablah\log \rho + \widetilde v \cdot\nablah (\log \rho)_h \bigr) = w_{p,hz} + \psi^z_{hz}.
\end{equation}
Therefore, after substituting \eqref{nonlinear-dhz} and \eqref{id:vertical_velocity-110} into $ I_{18} $, we have
\begin{align*}
	& I_{18} = - \int \biggl( \zeta_h \dt \psi^h_z + \zeta_h v_z \cdot \nablah \psi^h + \rho v_{hz} \cdot \nablah \psi^h + \rho v_z \cdot \nablah \psi^h_h \\
	& ~~~~ + \zeta_h w_{p,z} \dz \psi^h + \rho w_{p,hz} \dz \psi^h
	 + \rho w_{p,z} \dz \psi^h_h - ( \dvh \widetilde{\psi^h}\\
	& ~~~~ + \widetilde v \cdot\nablah \log \rho) ( \zeta_h \dz \psi^h + \rho \dz \psi^h_h)
	 - \rho (\dvh\widetilde{\psi^h_h} + \widetilde v_h \cdot\nablah \log\rho \\
	& ~~~~ + \widetilde v \cdot\nablah (\log\rho)_h) \dz \psi^h
	 + \zeta_h v\cdot\nablah \psi^h_z + \rho v_h \cdot\nablah \psi^h_z + \zeta_h w_p \dz \psi^h_z \\
	& ~~~~ + \rho w_{p,h} \dz \psi^h_z    \biggr) \cdot \psi^h_{hz} \idx
	 + \int \bblbrack \int_0^z \zeta_h \bigl(\dvh \widetilde{\psi^h} + \widetilde v \cdot\nablah \log\rho\bigr) \,dz \\
	& ~~~~ \times \blparenthese \dz \psi^h_z \cdot\psi^h_{hz} \brparenthese
	 + \int_0^z \rho \bigl( \dvh \widetilde{\psi^h_{h}} + \widetilde v_h \cdot\nablah \log\rho + \widetilde v \cdot\nablah (\log\rho)_h\bigr) \,dz \\
	& ~~~~ \times \blparenthese \dz \psi^h_z \cdot\psi^h_{hz} \brparenthese \bbrbrack \idx
	 \lesssim \norm{\psi^h_{hz}}{\Lnorm3} \times \biggl( \norm{\varepsilon^{-1}\zeta_h}{\Lnorm6} \norm{\varepsilon\dt \psi^h_z}{\Lnorm2} \\
	& ~~~~ + \norm{\zeta_h}{\Lnorm6} \norm{v_z}{\Lnorm4} \norm{\nablah\psi^h}{\Lnorm4}
	+ \norm{\rho}{\Lnorm\infty} \norm{v_{hz}}{\Lnorm2} \norm{\nablah \psi^h}{\Lnorm6} \\
	& ~~~~
	+ \norm{\rho}{\Lnorm{\infty}} \norm{v_z}{\Lnorm6} \norm{\nablah \psi^h_h}{\Lnorm2}
	 + \norm{\zeta_h}{\Lnorm6} \norm{w_{p,z}}{\Lnorm3} \norm{\dz \psi^h}{\Lnorm6} \\
	 & ~~~~ + \norm{\rho}{\Lnorm\infty} \norm{w_{p,hz}}{\Lnorm2} \norm{\dz \psi^h}{\Lnorm6}
	 + \norm{\rho}{\Lnorm\infty} \norm{w_{p,z}}{\Lnorm6} \norm{\dz \psi^h_h}{\Lnorm2} \\
	 & ~~~~ + (\norm{\nablah \psi^h}{\Lnorm6} + \norm{v}{\Lnorm\infty} \norm{\nablah \zeta}{\Lnorm6} )
	 \times (\norm{\zeta_h}{\Lnorm6} \norm{\dz\psi^h}{\Lnorm3} \\
	 & ~~~~
	 + \norm{\rho}{\Lnorm\infty}\norm{\dz\psi^h_h}{\Lnorm2}) + \norm{\rho}{\Lnorm\infty} \norm{\dz \psi^h}{\Lnorm6}
	  \times \bigl( \norm{\nablah \psi^h_h}{\Lnorm2} \\
	  & ~~~~ + \norm{v_h}{\Lnorm3} \norm{\nablah\zeta}{\Lnorm6}
	 + \norm{v}{\Lnorm\infty} (\norm{\zeta_{hh}}{\Lnorm2} + \norm{\zeta_h}{\Lnorm4}^2 ) \bigr) \\
	 & ~~~~+  \norm{\zeta_h}{\Lnorm6} \norm{v}{\Lnorm\infty} \norm{\nablah \psi^h_z}{\Lnorm2} + \norm{\rho}{\Lnorm\infty} \norm{v_h}{\Lnorm6} \norm{\nablah \psi^h_z}{\Lnorm2} \\
	& ~~~~ + \norm{\zeta_h}{\Lnorm6} \norm{w_p}{\Lnorm3} \norm{\dz \psi^h_z}{\Lnorm6} + \norm{\rho}{\Lnorm\infty} \norm{w_{p,h}}{\Lnorm2} \norm{\dz \psi^h_z}{\Lnorm6} \biggr) \\
	& ~~~~ + \int_0^1 \bblparenthese \hnorm{\zeta_h}{\Lnorm4} \bigl( \hnorm{\nablah \psi^h}{\Lnorm4} + \hnorm{v}{\Lnorm\infty} \hnorm{\nablah \zeta}{\Lnorm4} \bigr) + \hnorm{\rho}{\Lnorm\infty} \bigl(\hnorm{\nablah \psi^h_h}{\Lnorm2} \\
	& ~~~~ + \hnorm{v_h}{\Lnorm4} \hnorm{\nablah \zeta}{\Lnorm4}
	 + \hnorm{v}{\Lnorm\infty} ( \hnorm{\zeta_{hh}}{\Lnorm2} + \hnorm{\zeta_h}{\Lnorm4}^2 ) \bigr) \bbrparenthese \,dz \\
	 & ~~~~ \times \int_0^1 \hnorm{\dz \psi^h_z}{\Lnorm4} \hnorm{\psi^h_{hz}}{\Lnorm4} \,dz
	  \lesssim \biggl( \norm{\psi^h_{hz}}{\Lnorm2}^{1/2} \norm{\nabla\psi^h_{hz}}{\Lnorm2}^{1/2} + \norm{\psi^h_{hz}}{\Lnorm2}\biggr)  \\
	 & ~~~~ \times \biggl( \norm{\varepsilon^{-1}\nablah \zeta}{\Hnorm{1}} \norm{\varepsilon \dt \psi^h_z }{\Lnorm2}
	 + (\norm{\zeta}{\Hnorm{2}} + \norm{\zeta}{\Hnorm{2}}^2 + \norm{\rho}{\Lnorm\infty} \\
	 & ~~~~ + \norm{\rho}{\Lnorm\infty} \norm{\zeta}{\Hnorm{2}}^2 + \norm{\rho}{\Lnorm\infty} \norm{\zeta}{\Hnorm{2}} )
	 ( \norm{\psi^h}{\Hnorm{2}} +\norm{v}{\Hnorm{2}} \\
	& ~~~~ + \norm{w_{p,z}}{\Hnorm{1}})\norm{ \nabla \psi^h}{\Hnorm{1}}
	 + ( \norm{\rho}{\Lnorm\infty} \norm{\nablah\zeta}{\Hnorm{1}}^2
	 + \norm{\rho}{\Lnorm\infty} \norm{\nablah \zeta}{\Hnorm{1}} \\
	 & ~~~~ + \norm{\nablah \zeta}{\Hnorm{1}}^2 + \norm{\nablah \zeta}{\Hnorm{1}} + \norm{\rho}{\Lnorm\infty}  ) ( \norm{w_p}{\Hnorm{1}}
	+ \norm{\nabla \psi^h}{\Hnorm{1}} \\
	& ~~~~ + \norm{v}{\Hnorm{2}}) ( \norm{\nabla^2 \psi^h}{\Lnorm2} + \norm{\nabla^3 \psi^h}{\Lnorm2} ) \biggr) \\
	&  \lesssim \delta \norm{\nabla^3 \psi^h}{\Lnorm2}^2 + \delta \norm{\varepsilon \nabla \dt \psi^h}{\Lnorm2}^2 + C_\delta  Q(\mathcal E) \norm{\nabla \psi^h}{\Hnorm{1}}^2 \\
	& ~~~~ + C_\delta \blparenthese Q(\mathcal E) + 1 + \mathfrak G_p \brparenthese \mathfrak H_p \norm{\psi^h}{\Hnorm{2}}^2 .
\end{align*}
Similarly, after substituting \eqref{nonlinearities}, \eqref{id:vertical-velocity-020} and \eqref{id:vertical_velocity-200}, we have
\begin{align*}
	& I_{17} = - \int \bblbrack \biggl( \zeta_h \psi^h_z \cdot \nablah v_p + \zeta_h \psi^h \cdot \nablah v_{p,z} - \zeta_h (\dvh\widetilde{\psi^h} \\
	& ~~~~  + \widetilde v \cdot\nablah \log\rho ) \dz v_p
	 + \rho \psi^h_{hz} \cdot\nablah v_p + \rho \psi^h_h \cdot\nablah v_{p,z} \\
	& ~~~~ - \rho (\dvh \widetilde{\psi^h_h} + \tilde v_h \cdot\nablah \log\rho
	 + \widetilde v \cdot\nablah (\log\rho)_h) \dz v_p \\
	& ~~~~ + \rho \psi^h_z \cdot \nablah v_{p,h} + \rho \psi^h \cdot\nablah v_{p,hz}
	 - \rho (\dvh\widetilde{\psi^h} \\
	 & ~~~~ + \widetilde v \cdot\nablah \log\rho ) \dz v_{p,h} \biggr) \cdot \psi^h_{hz} \bbrbrack \idx
	 + \int \bblbrack \int_0^z \bblparenthese \dvh\widetilde{\psi^h} \\
	& ~~~~ + \widetilde v \cdot\nablah \log\rho \bbrparenthese  \,dz \times  \bigl( \zeta_h \dz v_{p,z} \cdot\psi^h_{hz} + \rho \dz v_{p,hz} \cdot\psi^h_{hz}\bigr) \\
	& ~~~~
	+ \int_0^z \bblparenthese \dvh \widetilde{\psi^h_h}
	 + \widetilde v_h \cdot\nablah \log\rho + \widetilde v \cdot\nablah (\log\rho)_h \bbrparenthese \,dz \\
	& ~~~~ \times \blparenthese \rho \dz v_{p,z} \cdot\psi^h_{hz}\brparenthese  \bbrbrack \idx
	\lesssim \biggl( \norm{\nablah \zeta}{\Hnorm{1}} ( \norm{\psi^h}{H^2} \\
	& ~~~~ + \norm{v}{H^2} \norm{\nablah \zeta}{\Hnorm{1}} + \norm{\nablah \zeta}{\Hnorm{1}} )
	 + \norm{\rho}{\Lnorm{\infty}} \bigl( \norm{\psi^h}{\Hnorm{2}} \\
	 & ~~~~ + \norm{v}{\Hnorm{2}} ( \norm{\nablah \zeta}{\Hnorm{1}}+ \norm{\nablah \zeta}{\Hnorm{1}}^2 ) \bigr) \biggr) \norm{v_p}{H^3} \\
	& ~~~~ \times \biggl( \norm{\psi^h_{hz}}{\Lnorm2}^{1/2} \norm{\nabla\psi^h_{hz}}{\Lnorm2}^{1/2}
	 + \norm{\psi^h_{hz}}{\Lnorm2}\biggr)
	 + \int_0^1 \bblparenthese \hnorm{\nablah \psi^h}{\Lnorm4} \\
	 & ~~~~ + \varepsilon \hnorm{v}{\Lnorm\infty} \hnorm{\varepsilon^{-1}\nablah\zeta}{\Lnorm4} \bbrparenthese \,dz \times \int_0^1 \bblparenthese \hnorm{\zeta_h}{\Lnorm4} \hnorm{\dz v_{p,z}}{\Lnorm4} \hnorm{\psi^h_{hz}}{\Lnorm4} \\
	 & ~~~~ + \hnorm{\rho}{\Lnorm\infty}
	 \hnorm{\dz v_{p,hz}}{\Lnorm2} \hnorm{\psi^h_{hz}}{\Lnorm4} \bbrparenthese  \,dz
	 + \int_0^1 \bblparenthese \hnorm{\nablah \psi^h_h}{\Lnorm2} + \hnorm{v_h}{\Lnorm4} \hnorm{\nablah \zeta}{\Lnorm4} \\
	 & ~~~~+ \hnorm{v}{\Lnorm\infty}( \hnorm{\zeta_{hh}}{\Lnorm2}+ \hnorm{\zeta_h}{\Lnorm4}^2) \bbrparenthese \,dz \cdot\int_0^1\hnorm{\rho}{\Lnorm\infty} \hnorm{\dz v_{p,z}}{\Lnorm4} \hnorm{\psi^h_{hz}}{\Lnorm4}  \,dz   \\
	& ~~ \lesssim \delta \norm{\nabla \psi^h_{hz}}{\Lnorm2}^2 + \delta \norm{\psi^h_{hz}}{\Lnorm{2}}^2 + C_\delta \blparenthese Q(\mathcal E) + 1 + \mathfrak G_p \brparenthese \mathfrak H_p \\
	& ~~~~ \times \blparenthese \norm{\psi^h}{\Hnorm{2}}^2 + \varepsilon^2 \brparenthese.\\
	& I_{16} = \int \blparenthese \zeta_h Q_{p,z} \cdot\psi^h_{hz} + \zeta Q_{p,hz} \cdot\psi^h_{hz} \brparenthese \idx = - \int \blparenthese \zeta_h Q_{p} \cdot \psi^h_{hzz} \\
	& ~~~~ + \zeta Q_{p,h} \cdot \psi^h_{hzz} \brparenthese \idx
	 \lesssim \delta \norm{\psi^h_{hzz}}{\Lnorm{2}}^2  + C_\delta \blparenthese Q(\mathcal E)  + \mathfrak G_p \brparenthese \mathfrak H_p \varepsilon^2 .
\end{align*}
Therefore after summing up the estimates of $ I_{16}, I_{17},I_{18} $, above, and \eqref{ee:005}, we conclude the proof of \eqref{ee:vertical-est-dhz-psi}.
\end{proof}

\subsection{Proof of Proposition \ref{prop:a-prior-estimates}}

After applying similar arguments as in the proofs of Lemma \ref{lm:horizontal-derivative}, Lemma \ref{lm:vertical-derivative} and Lemma \ref{lm:horizontal-vertical-derivative}, one can easily check the following inequalities hold:
\begin{align*}
	& \dfrac{d}{dt} \biggl\lbrace \dfrac{1}{2} \norm{\rho^{1/2} \psi^h}{\Lnorm{2}}^2 + \dfrac{c^2_s}{2\rho_0} \norm{\varepsilon^{-1}\xi}{\Lnorm{2}}^2  \biggr\rbrace + C_{\mu,\delta} \norm{\nabla \psi^h}{\Lnorm{2}}^2 \\
	& ~~~~ \leq \delta \blparenthese \norm{\nabla \psi^h}{\Lnorm{2}}^2 + \norm{\varepsilon^{-1} \nablah \xi }{\Hnorm{1}}^2 \brparenthese + C_\delta Q(\mathcal E) \blparenthese  \norm{\nabla\psi^h}{\Hnorm{1}}^2 \\
	& ~~~~ + \norm{\varepsilon^{-1} \nablah \xi }{\Hnorm{1}}^2 \brparenthese
	 + C_\delta \blparenthese Q(\mathcal E) + 1 + \mathfrak G_p \brparenthese \mathfrak H_p \blparenthese \varepsilon^2 + \norm{\psi^h}{\Hnorm{2}}^2 \brparenthese;    \\
	& \dfrac{d}{dt} \biggl\lbrace \dfrac{1}{2} \norm{\rho^{1/2}\psi^h_h}{\Lnorm{2}}^2 + \dfrac{c^2_s}{2\rho_0} \norm{\varepsilon^{-1} \xi_h}{\Lnorm{2}}^2 \biggr\rbrace + C_{\mu,\lambda} \norm{\nabla \psi^h_h}{\Lnorm{2}}^2 \\
	& ~~~~ \leq \delta \blparenthese \norm{\nabla^2 \psi^h}{\Lnorm{2}}^2 + \norm{\varepsilon^{-1} \nablah \xi}{\Hnorm{1}}^2 + \norm{\varepsilon \psi^h_t}{\Lnorm{2}}^2 \brparenthese \\
	& ~~~~ + C_\delta Q(\mathcal E) \blparenthese \norm{\nabla\psi^h}{\Hnorm{1}}^2 + \norm{\varepsilon^{-1} \nablah \xi }{\Hnorm{1}}^2 \brparenthese \\
	& ~~~~ + C_\delta  \blparenthese Q(\mathcal E) + 1 + \mathfrak G_p \brparenthese \mathfrak H_p \blparenthese \varepsilon^2 + \norm{\psi^h}{\Hnorm{2}}^2 \brparenthese;\\
	& \dfrac{d}{dt} \norm{\rho^{1/2} \psi^h_z}{\Lnorm{2}}^2 + C_{\mu,\lambda} \norm{\nabla \psi^h_z}{\Lnorm{2}}^2 \leq \delta \blparenthese \norm{\nabla^2 \psi^h}{\Lnorm{2}}^2 + \norm{\psi^h_z}{\Lnorm{2}}^2 \\
	& ~~~~ + \norm{\varepsilon \dt \psi^h}{\Lnorm{2}}^2 \brparenthese + C_\delta Q(\mathcal E) \norm{\nabla\psi^h}{\Hnorm{1}}^2 \\
	& ~~~~ + C_\delta \blparenthese Q(\mathcal E) + 1 + \mathfrak G_p \brparenthese \mathfrak H_p \blparenthese \varepsilon^2 + \norm{\psi^h}{\Hnorm{2}}^2 \brparenthese.
\end{align*}
Therefore, the above inequalities, together with \eqref{ene:temporal-derivative}, \eqref{ee:006}, \eqref{ee:008}, \eqref{ee:horizontal-derivative}, \eqref{ee:dhh-xi}, \eqref{ee:vertical-est-dzz-psi} and \eqref{ee:vertical-est-dhz-psi}, imply that there exist positive constants
$ c_i, ~ i \in \lbrace 1,2 \cdots 10 \rbrace $, such that
\begin{equation}\label{estimates-total-01}
	\begin{aligned}
		& \dfrac{d}{dt} \mathcal E_{LM} + \mathcal D_{LM} \leq \blparenthese \delta + \varepsilon^2 C_\delta (Q(\mathcal E) + 1 + \mathfrak G_p) + (C_\delta + \varepsilon^2) Q(\mathcal E) \brparenthese \mathcal D  \\
		& ~~~~ + C_\delta \blparenthese Q(\mathcal E) + 1 + \mathfrak G_p \brparenthese \mathfrak H_p \blparenthese \varepsilon^2 + \mathcal E\brparenthese,
	\end{aligned}
\end{equation}
where we denote by
\begin{align}
	& \mathcal E_{LM} = \mathcal E_{LM}(t) := \dfrac{c_1}{2} \norm{\rho^{1/2} \varepsilon \psi^h_t}{\Lnorm{2}}^2 + \dfrac{c_1 c^2_s}{2\rho_0} \norm{\xi_t}{\Lnorm{2}}^2 + \dfrac{c_2}{2} \norm{\rho^{1/2} \nablah^2 \psi^h}{\Lnorm{2}}^2 \nonumber \\
	& ~~~~ + \dfrac{c_2 c^2_s}{2\rho_0} \norm{\varepsilon^{-1} \nablah^2 \xi}{\Lnorm{2}}^2 + c_3 \norm{\rho^{1/2} \psi^h_{zz}}{\Lnorm{2}}^2 + c_4 \norm{\rho^{1/2}\nablah \psi^h_z}{\Lnorm{2}}^2 \nonumber\\
	& ~~~~ + \dfrac{c_5}{2} \norm{\rho^{1/2} \nablah \psi^h}{\Lnorm{2}}^2 + \dfrac{c_5 c^2_s}{2\rho_0} \norm{\varepsilon^{-1}\nablah \xi}{\Lnorm{2}}^2 + c_6 \norm{\rho^{1/2} \psi^h_z}{\Lnorm{2}}^2 \nonumber\\
	& ~~~~ + \dfrac{c_7}{2} \norm{\rho^{1/2} \psi^h}{\Lnorm{2}}^2 + \dfrac{c_7c^2_s}{2\rho_0} \norm{\varepsilon^{-1}\xi}{\Lnorm{2}}^2, \label{functional:instant-energy}  \\
	& \mathcal D_{LM} = \mathcal D_{LM}(t) := c_1 \norm{\varepsilon \nabla \psi^h_t}{\Lnorm{2}}^2 + c_2 \norm{\nabla \nablah^2 \psi^h}{\Lnorm{2}}^2 + c_3 \norm{\nabla \psi^h_{zz}}{\Lnorm{2}}^2 \nonumber \\
	& ~~~~ + c_4 \norm{\nabla \nablah \psi^h_z}{\Lnorm{2}}^2 + c_5 \norm{\nabla \nablah \psi^h}{\Lnorm{2}}^2 + c_6 \norm{\nabla \psi^h_z}{\Lnorm{2}}^2 \nonumber\\
	& ~~~~ + c_7 \norm{\nabla \psi^h}{\Lnorm{2}}^2 + c_8 \norm{\xi_t}{\Lnorm{2}}^2 + c_9 \norm{\varepsilon^{-1} \nablah \xi }{\Hnorm{1}}^2 + C_{10} \norm{\varepsilon \rho
	\psi^h_t}{\Lnorm{2}}^2. \label{functional:instant-dissipation}
\end{align}
Under the assumption \eqref{boundness-density}, it is easy to check that
\begin{equation}\label{equivalent-functionals}
\begin{gathered}
	\mathcal E \lesssim \mathcal E_{LM} \lesssim \mathcal E, ~
	\mathcal D \lesssim \mathcal D_{LM} \lesssim \mathcal D,
\end{gathered}
\end{equation}
where $ \mathcal E = \mathcal E(t) $ and $ \mathcal D = \mathcal D(t) $ are defined in \eqref{functional:energy} and \eqref{functional:dissipation}, respectively.
Therefore \eqref{estimates-total-01} can be written, after choosing $ \delta $ small enough, as
\begin{equation}\label{estimates-total-02}
	\begin{aligned}
		& \dfrac{d}{dt} \mathcal E_{LM} + \mathcal D_{LM} \leq \blparenthese \varepsilon^2 C(Q(\mathcal E) + 1 + \mathfrak G_p) + (1 + \varepsilon^2) Q(\mathcal E) \brparenthese \mathcal D_{LM}  \\
		& ~~~~ + C \blparenthese Q(\mathcal E) + 1 + \mathfrak G_p \brparenthese \mathfrak H_p \blparenthese \varepsilon^2 + \mathcal E_{LM} \brparenthese.
	\end{aligned}
\end{equation}
Then after applying Gr\"onwall's inequality to \eqref{estimates-total-02}, one concludes that,
\begin{align*}
	& \sup_{0\leq t \leq T} \mathcal E_{LM}(t) + \int_0^T \mathcal D_{LM}(t) \,dt \lesssim e^{C \int_0^T \blparenthese Q(\mathcal E) + 1 + \mathfrak G_p\brparenthese \mathfrak H_p \,dt} \\
	& ~~ \times \biggl\lbrace  \varepsilon^2 + \mathcal E_{LM}(0) + \int_0^T \bblbrack \bigl\lbrack \varepsilon^2 C (Q(\mathcal E) + 1 + \mathfrak G_p) \\
	& ~~~~ + (1 + \varepsilon^2) Q(\mathcal E) \bigr\rbrack \mathcal D_{LM} \bbrbrack \,dt  \biggr\rbrace.
\end{align*}
Under the assumptions of Proposition \ref{prop:a-prior-estimates}, this completes the proof of \eqref{prop:stability-CPE}.


\section{Low Mach number limit}\label{sec:low-mach-number-limit}
In this section, we will establish the asymptotic behavior of $ (\xi, \psi^h) = (\xi^\varepsilon, \psi^{\varepsilon,h}) $ as $ \varepsilon \rightarrow 0^+ $. In particular, we prove Theorem \ref{thm:low-mach} in this section.

 First, as a consequence of Theorem \ref{thm:global-pe}, we have the following:
\begin{corollary}\label{cor:choosing-s}
Under the same assumptions of Theorem \ref{thm:global-pe}, consider any integer $ s \geq 3 $. Then \eqref{boundness-and-decay} holds true for $ v_{p,in} \in H^s(\Omega_h \times 2\mathbb T) $, with the compatibility conditions as in \eqref{compatible-condition-pe}.
\end{corollary}
\begin{proof}
	Directly from the conclusion of Theorem \ref{thm:global-pe}, for $ s \geq 3 $ and $ v_{p,in} \in H^s(\Omega) $ as stated in the theorem, one has
	\begin{align*}
		& \sup_{0\leq t < \infty} \blparenthese \norm{v_p(t)}{\Hnorm{3}} + \norm{v_{p,t}(t)}{\Lnorm{2}} \brparenthese +  \int_0^\infty \bblparenthese \norm{v_{p}(t)}{\Hnorm{3}}^2 + \norm{v_{p,t}(t)}{\Hnorm{1}}^2 \\
		& ~~~~ + \norm{v_p(t)}{\Hnorm{2}} \bbrparenthese \,dt \leq C_{p,in,3} + 1,
	\end{align*}
	where we have used the Gagliardo-Nirenberg interpolation inequality:
	\begin{align*}
		& \int_0^\infty \norm{v_{p}(t)}{\Hnorm{2}} \,dt  \lesssim \int_0^\infty \norm{v_{p}(t)}{\Lnorm{2}}^{1/3} \norm{v_{p}(t)}{\Hnorm{3}}^{2/3} \,dt \\
		& ~~~~ \lesssim \int_0^\infty e^{-\frac{c}{3}t} \,dt \times C_{p,in,3}^{1/2} \leq C_{p,in,3} + 1 .
	\end{align*}
	Moreover, applying the Minkowski and H\"older inequalities to the expression of $ w_p $ as in \eqref{id:vertical-vel-PE} yields
	\begin{gather*}
		\sup_{0\leq t < \infty}\norm{w_p(t)}{\Hnorm{1}} \leq \sup_{0\leq t < \infty} \norm{v_{p}(t)}{\Hnorm{2}} \leq C_{p,in,2} + 1, \\
		\int_0^\infty\norm{w_p(t)}{\Hnorm{2}}^2 \,dt \leq \int_0^\infty \norm{v_p(t)}{\Hnorm{3}}^2 \,dt \leq C_{p,in,2}.
	\end{gather*}

	What is left is to estimate $$ \norm{\rho_1}{\Hnorm{2}}, \norm{\rho_{1,t}}{\Lnorm{2}}, \int_0^\infty\blparenthese \norm{\rho_{1,tt}}{\Lnorm{2}}^2 + \norm{\rho_{1,t}}{\Hnorm{1}}^2 + \norm{\rho_1}{\Hnorm{2}}^2 \brparenthese \,dt. $$ To do so, we write down the elliptic problem for $ \rho_1 $, which is obtained by taking average over the $ z $-variable and then taking $ \dvh $ in \subeqref{PE}{2}, below:
	\begin{equation}\label{eq:pressure-pe}
		- c^2_s \deltah \rho_1 = \rho_0 \int_0^1 \dvh \blparenthese \dvh (v_p\otimes v_p ) \brparenthese \,dz ~~ \text{in} ~ \Omega_h,  ~ \int_{\Omega_h} \rho_1 \idxh = 0.
	\end{equation}
	Then the $ L^p $ estimate of the Riesz transform implies, together with the Minkowski, H\"older and Sobolev embedding inequalities,
	\begin{equation*}
		\norm{\rho_1}{\Hnorm{2}} \leq \hnorm{\rho_1}{\Hnorm{2}} \lesssim \norm{ \abs{v_p}{2}}{\Hnorm{2}} \lesssim \norm{v_p}{\Hnorm{2}}^2.
	\end{equation*}
	Consequently, for $ s \geq 2 $,
	\begin{equation}\label{19Jun2018-pressure}
	\begin{gathered}
		\sup_{0\leq t < \infty}\norm{\rho_1(t)}{\Hnorm{2}} + \int_0^\infty \norm{\rho_1(t)}{\Hnorm{2}}^2 \,dt \lesssim \sup_{0\leq t <\infty} \norm{v_p(t)}{\Hnorm{2}}^2 \\
		+ \sup_{0\leq t < \infty}  \norm{v_p(t)}{\Hnorm{2}}^2 \times \int_0^\infty \norm{v_p(t)}{\Hnorm{2}}^2 \,dt \leq ( 1 + C_{p,in,1} ) C_{p,in,2} .
	\end{gathered}
	\end{equation}
	Furthermore, after taking time derivatives of \eqref{eq:pressure-pe}, we have the following elliptic problems:
	\begin{align}
	& \label{eq:pressure-pe-1} \begin{cases} - c^2_s \deltah \rho_{1,t} = 2 \rho_0 \int_0^1 \dvh \blparenthese \dvh (v_p\otimes v_{p,t} ) \brparenthese \,dz ~~ \text{in} ~ \Omega_h,  \\
	\int_{\Omega_h} \rho_{1,t} \idxh = 0; \end{cases} \\
	&\label{eq:pressure-pe-2}  \begin{cases} - c^2_s \deltah \rho_{1,tt} = 2 \rho_0 \int_0^1 \dvh \blparenthese \dvh (v_p\otimes v_{p,tt} + v_{p,t} \otimes v_{p,t} ) \brparenthese \,dz \\ ~~ \text{in} ~ \Omega_h,  ~~
	\int_{\Omega_h} \rho_{1,tt} \idxh = 0. \end{cases}
	\end{align}
	Thus similarly, one has, for $ s \geq 3 $,
	\begin{equation}\label{cor-est:001}
	\begin{aligned}
		& \norm{\rho_{1,t}}{\Lnorm{2}}  \lesssim \norm{\abs{v_p}{} \abs{v_{p,t}}{}}{\Lnorm{2}} \lesssim \norm{v_{p}}{\Hnorm{2}} \norm{v_{p,t}}{\Lnorm{2}} \leq C_{p,in,2} , \\
		& \int_0^\infty \norm{\rho_{1,t}}{\Hnorm{1}}^2 \,dt  \lesssim \int_0^\infty \blparenthese \norm{\abs{\nablah v_p}{} \abs{v_{p,t}}{}}{\Lnorm{2}}^2 + \norm{\abs{v_p}{} \abs{\nablah v_{p,t}}{}}{\Lnorm{2}}^2 \brparenthese \,dt \\
		& ~~ \lesssim \int_0^\infty \norm{v_p}{\Hnorm{2}}^2 \norm{v_{p,t}}{\Hnorm{1}}^2 \,dt \lesssim \sup_{0\leq t <\infty} \norm{v_p}{\Hnorm{2}}^2 \\
		& ~~~~  \times \int_0^\infty  \norm{v_{p,t}}{\Hnorm{1}}^2 \,dt \leq C_{p,in,2}^2 , \\
		& \int_0^\infty \norm{\rho_{1,tt}}{\Lnorm{2}}^2  \lesssim \int_0^\infty \blparenthese \norm{\abs{v_p}{} \abs{v_{p,tt}}{}}{\Lnorm{2}}^2 + \norm{\abs{v_{p,t}}{2}}{\Lnorm{2}}^2 \brparenthese \,dt  \\
		& ~~ \lesssim \int_0^\infty \blparenthese  \norm{v_p}{\Hnorm{2}}^2 \norm{v_{p,tt}}{\Lnorm{2}}^2 + \norm{v_{p,t}}{\Hnorm{1}}^4 \brparenthese \,dt\\
		& ~~ \lesssim \sup_{0\leq t < \infty} \norm{v_{p}}{\Hnorm{2}}^2 \times \int_0^\infty \norm{v_{p,tt}}{\Lnorm{2}}^2 \,dt + \sup_{0\leq t < \infty } \norm{v_{p,t}}{\Hnorm{1}}^2 \\
		& ~~~~ \times  \int_0^\infty \norm{v_{p,t}}{\Hnorm{1}}^2 \leq C_{p,in,2}\int_0^\infty \norm{v_{p,tt}}{\Lnorm{2}}^2 \,dt  + C_{p,in,3}^2 .
	\end{aligned}	
	\end{equation}
	On the other hand, after taking time derivative of \subeqref{PE}{2}, we have the identity,
	\begin{align*}
		& \rho_0 v_{p,tt} = - \rho_0 \blparenthese v_{p} \cdot \nablah v_p + w_p \dz v_p \brparenthese_t - c^2_s \nablah \rho_{1,t} \\
		& ~~~~ + \mu \deltah v_{p,t} + \lambda\nablah \dvh v_{p,t} + \partial_{zz} v_{p,t}.
	\end{align*}
	Therefore, directly one has,
	\begin{align*}
		& \norm{v_{p,tt}}{\Lnorm{2}} \lesssim \norm{\rho_{1,t}}{\Hnorm{1}} + \norm{v_{p,t}}{\Hnorm{2}} + \norm{v_p \cdot \nablah v_{p,t}}{\Lnorm{2}} + \norm{v_{p,t} \cdot \nablah v_{p}}{\Lnorm{2}} \\
		& ~~~~ + \norm{w_{p} \dz v_{p,t}}{\Lnorm{2}} + \norm{w_{p,t} \dz v_{p}}{\Lnorm{2}} \lesssim \norm{\rho_{1,t}}{\Hnorm{1}} + \norm{v_{p,t}}{\Hnorm{2}} \\
		& ~~~~ + \norm{v_{p}}{\Hnorm{2}} \norm{v_{p,t}}{\Hnorm{1}} + \norm{w_p}{\Hnorm{2}} \norm{v_{p,t}}{\Hnorm{1}} + \norm{w_{p,t}}{\Lnorm{2}} \norm{v_{p}}{\Hnorm{3}}\\
		& ~~ \lesssim \norm{\rho_{1,t}}{\Hnorm{1}} + \norm{v_{p,t}}{\Hnorm{2}} + \norm{v_{p}}{\Hnorm{3}} \norm{v_{p,t}}{\Hnorm{1}},
	\end{align*}
	where we have applied the Minkowski, Sobolev embedding and H\"older inequalities, and the following inequalities as the consequence of \eqref{id:vertical-vel-PE},
	\begin{equation*}
		\norm{w_p}{\Hnorm{2}} \leq \norm{v_{p}}{\Hnorm{3}}, ~ \norm{w_{p,t}}{\Lnorm{2}} \leq \norm{v_{p,t}}{\Hnorm{1}}.
	\end{equation*}
	Consequently, one concludes that, for $ s \geq 3 $,
	\begin{align*}
		& \int_{0}^\infty\norm{v_{p,tt}}{\Lnorm{2}}^2 \,dt \lesssim \int_0^\infty \norm{\rho_{1,t}}{\Hnorm{1}}^2\,dt + \blparenthese 1 + \sup_{0\leq t < \infty} \norm{v_{p}}{\Hnorm{3}}^2 \brparenthese \\
		& ~~~~ \times \int_0^\infty \norm{v_{p,t}}{\Hnorm{2}}^2 \,dt \lesssim C_{p,in,2}^2 + \bigl( 1 + C_{p,in,3} \bigr) C_{p,in,3},
	\end{align*}
	where we have substituted inequality \subeqref{cor-est:001}{2}. Thus \subeqref{cor-est:001}{3} yields
	\begin{equation*}
		 \int_0^\infty \norm{\rho_{1,tt}}{\Lnorm{2}}^2 \lesssim C_{p,in,3}^2 + C_{p,in,2} \blparenthese C_{p,in,2}^2 + ( 1 + C_{p,in,3} ) C_{p,in,3} \brparenthese < \infty.
	\end{equation*}
	This completes the proof.
\end{proof}

Now, given $ v_{p,in} \in H^s(\Omega_h \times 2\mathbb T) $, for any integer $ s \geq 3 $, which is even in the $ z $-variable and satisfies the compatibility conditions \eqref{compatible-condition-pe}, one can apply the conclusion of Proposition \ref{prop:a-prior-estimates} to establish the global bound of the perturbation energy $ \mathcal E $, provided it is initially small. This is done through a continuity argument. We state first the proposition concerning the local well-posedness of solutions to system \eqref{eq:perturbation} with $ \mathcal E_{in} $ small enough:
\begin{proposition}\label{prop:local-well-perturbation}
	Let $ v_p $ be the solution to system \eqref{PE}, as stated in Theorem \ref{thm:global-pe} with initial data $ v_{p,in} \in H^s(\Omega_h \times 2\mathbb T) $, for an integer $ s \geq 3 $. Consider the initial data $ (\xi_{in}, \psi^h_{in} ) \in H^2(\Omega_h\times 2 \mathbb T) \times H^2(\Omega_h \times 2\mathbb T)  $ as in \eqref{initial-data-well-prepared} and satisfying the compatibility condition \eqref{cmpt-cds-perturbed}.
	There is a positive constant $ \bar\varepsilon \in (0,1) $, small enough, and a positive time $ T_{\bar\epsilon}\in (0,\infty) $,  such that if $ \varepsilon \in (0, \bar\varepsilon) $ and
	$ \mathcal E_{in} 
	\leq \bar\varepsilon $, there exists a unique strong solution $ (\xi^\varepsilon, \psi^{\varepsilon,h}) \in L^\infty(0,T_{\bar\varepsilon};H^2(\Omega_h \times 2\mathbb T)) $, with $ \psi^{\varepsilon,z} $ as in \eqref{id:vertical_perturbation}, to system \eqref{eq:perturbation} in the time interval $ [0,T_{\bar\varepsilon}] $. The existence time $ T_{\bar{\varepsilon}} $ depends only on $ \bar{\varepsilon} $ and $ \norm{v_{p,in}}{\Hnorm{3}} $ and is independent of $ \varepsilon $. Here $ \mathcal E_{in} $ is as in \eqref{functional-initial-energy}. Moreover, $ (\dt \xi^\varepsilon, \dt \psi^{\varepsilon,h}) \in L^\infty(0,T_{\bar\varepsilon};L^2(\Omega_h \times 2\mathbb T)) $, $ \rho = \rho_0 + \varepsilon^2 \rho_1 + \xi^\varepsilon \in (\frac{1}{2}\rho_0, 2\rho_0) $ in $ \Omega \times [0,T_{\bar\varepsilon}]  $, and
	 there is a constant $ C'' > 0 $, independent of $ \varepsilon $, such that,
	 $$ \sup_{0\leq t \leq T_{\bar\varepsilon}} \mathcal E(t) \leq C'' \mathcal E_{in}, $$
	where $ \mathcal E(t) $ is as in \eqref{functional:energy}.
\end{proposition}
The proof of Proposition \ref{prop:local-well-perturbation} can be done via a fixed point argument similar to that in our previous work \cite{LT2018a} and it is omitted here.

Now we are ready to establish the proof of Theorem \ref{thm:low-mach}
\begin{proof}[Proof of Theorem \ref{thm:low-mach}]
Consider $ s \geq 3 $ and $ \varepsilon \in (0,\bar{\varepsilon}) $ with $ \bar\varepsilon \in (0,1) $ as given in Proposition \ref{prop:local-well-perturbation}. Let the initial data $ (\xi_{in}, \psi^{h}_{in} ) \in H^2(\Omega_h) \times H^2(\Omega_h \times 2\mathbb T) $ satisfy \eqref{initial-data-well-prepared}, the compatibility conditions \eqref{cmpt-cds-perturbed}, and $ \mathcal E_{in} \leq \varepsilon^2 $, where $ \mathcal E_{in} $ is as in \eqref{functional-initial-energy}. Then $ \mathcal E_{in} \leq \bar\varepsilon $, and there is a strong solution to system \eqref{eq:perturbation} as stated by Proposition \ref{prop:local-well-perturbation} in the time interval $ [0,T_{\bar\varepsilon}] $, for some $ T_{\bar\varepsilon} \in (0,\infty) $, independent of $ \varepsilon $. The strong solution satisfies
\begin{gather*}
	\rho \in (\frac{1}{2}\rho_0, 2\rho_0) ~ \text{in}~ (\Omega_h \times 2\mathbb T) \times [0,T_{\bar\varepsilon}], \\
	\text{and} ~~~~ \sup_{0\leq t \leq T_{\bar\varepsilon}} \mathcal E(t) \leq C'' \mathcal E_{in} \leq C'' \varepsilon^2.
\end{gather*}

Such estimates, together with Theorem \ref{thm:global-pe} and Corollary \ref{cor:choosing-s}, imply that the assumptions in Proposition \ref{prop:a-prior-estimates} hold true in the time interval $ [0,T_{\bar\varepsilon}] $.
Therefore applying \eqref{prop:stability-CPE} yields
\begin{equation*}
	\begin{aligned}
		& \sup_{0 \leq t \leq T_{\bar\varepsilon}} \mathcal E(t) + \int_0^{T_{\bar\varepsilon}} \mathcal D(t) \,dt  \leq C' e^{C' + Q(C''\varepsilon^2)} \biggl\lbrace \varepsilon^2 + \varepsilon^2 \\
		& ~~~~ + \bblparenthese \varepsilon^2 + (\varepsilon^2  + 1 ) Q(C''\varepsilon^2) \bbrparenthese \int_0^{T_{\bar\varepsilon}} \mathcal D(t) \,dt \biggr\rbrace\\
		& ~~ \leq 2 C' e^{2C'} \varepsilon^2 + \dfrac{1}{2} \int_0^{T_{\bar\varepsilon}} \mathcal D(t)\,dt \biggr\rbrace,
		\end{aligned}
\end{equation*}
provided $ \varepsilon \in (0,\varepsilon_1) \subset (0,\bar\varepsilon) $, where $ \varepsilon_1 $ is small enough such that $ Q(C''\varepsilon_1^2) \leq C' $ and $ C' e^{2C'} (\varepsilon_1^2 + (\varepsilon_1^2 + 1) Q(C''\varepsilon_1^2) ) \leq 1/2 $.
This inequality yields that,
\begin{equation}\label{loot-initial}
	\sup_{0 \leq t \leq T_{\bar\varepsilon}} \mathcal E(t) + \int_0^{T_{\bar\varepsilon}} \mathcal D(t) \,dt  \leq C''' \varepsilon^2 < \bar\varepsilon,
\end{equation}
where $ C''' = 4 C' e^{2C'} $, and provided $ \varepsilon_1 $ is small such that $ C''' \varepsilon_1^2 < \bar\varepsilon $.
In particular, $ \mathcal E(T_{\bar\varepsilon}) \leq \bar\varepsilon $. We apply Proposition \ref{prop:local-well-perturbation} again in the time interval $ [T_{\bar\varepsilon}, 2T_{\bar\varepsilon}] $, which states that  there exists a strong solution satisfying
\begin{gather*}
	\rho \in (\frac{1}{2} \rho_0, 2 \rho_0) ~ \text{in} ~ (\Omega_h \times 2\mathbb T) \times [T_{\bar\varepsilon}, 2T_{\bar\varepsilon}], ~~~~\\
	 \text{and}  ~~~~\sup_{T_{\bar\varepsilon} \leq t \leq 2T_{\bar\varepsilon}} \mathcal E(t) \leq C'' \mathcal E(T_{\bar\varepsilon}) \leq C'' C''' \varepsilon^2\end{gather*}
Together with \eqref{loot-initial}, this implies
\begin{equation*}
	\sup_{0\leq t \leq 2T_{\bar\varepsilon}} \mathcal E(t) \leq C'''' \varepsilon^2 ~ \text{with} ~ C'''' = \max\lbrace C'', C''C''' \rbrace,
\end{equation*}
	Consequently Proposition \ref{prop:a-prior-estimates} applies. In particular, \eqref{prop:stability-CPE} yields
\begin{equation*}
	\begin{aligned}
		& \sup_{0 \leq t \leq 2T_{\bar\varepsilon}} \mathcal E(t) + \int_0^{2T_{\bar\varepsilon}} \mathcal D(t) \,dt  \leq C' e^{C' + Q(C''''\varepsilon^2)} \biggl\lbrace \varepsilon^2 + \varepsilon^2 \\
		& ~~~~ + \bblparenthese \varepsilon^2 + (\varepsilon^2  + 1 ) Q(C''''\varepsilon^2) \bbrparenthese \int_0^{2T_{\bar\varepsilon}} \mathcal D(t) \,dt \biggr\rbrace.
		\end{aligned}
\end{equation*}
As above, this implies 
\begin{equation}\label{global-stability-est}
	\sup_{0 \leq t \leq 2T_{\bar\varepsilon}} \mathcal E(t) + \int_0^{2T_{\bar\varepsilon}} \mathcal D(t) \,dt  \leq C''' \varepsilon^2 < \bar\varepsilon,
\end{equation}
provided $ \varepsilon \in (0,\varepsilon_2) \subset (0,\varepsilon_1) \subset (0,\bar\varepsilon) $, for $ \varepsilon_2 $ small enough such that $ Q(C'''' \varepsilon_2^2) \leq Q(C''\varepsilon_1^2) $.
Then inductively, without needing to determine the smallness of $ \varepsilon $ again, the arguments from \eqref{loot-initial} to \eqref{global-stability-est} hold true for $ T_{\bar\varepsilon}, 2T_{\bar\varepsilon} $ replaced by $ n T_{\bar\varepsilon}, (n+1)T_{\bar\varepsilon} $, $ n \geq 2 $, respectively.
In particular, \eqref{global-stability-est} holds true for $ 2 T_{\bar\varepsilon} $ replaced by $ (n+1) T_{\bar\varepsilon} $. Recall that $ T_{\bar\varepsilon} $ is independent of $ \varepsilon $. This concludes the proof of \eqref{global-stability-thm}. \eqref{low-mach-limit-thm} is a direct consequence of \eqref{global-stability-thm}, \eqref{ansatz}, \eqref{id:vertical_velocity}, \eqref{id:vertical-vel-PE} and the fact that $ \norm{\rho_1}{\Hnorm{2}} < \infty $ as in \eqref{19Jun2018-pressure}. Therefore let $ \varepsilon_0 = \varepsilon_2 $, and we complete the proof of Theorem \ref{thm:low-mach}.
\end{proof}

\section{Global regularity estimates of the solution to the primitive equations}\label{sec:global-pe}

Let us recall the primitive equations  \eqref{PE} first. We shorten the notations $ v= v_p , w = w_p  $ in this section. Recall that
\begin{equation*}\tag{\ref{PE}}
\begin{cases}
\dvh v + \dz w = 0 &\text{in} ~ \Omega_h \times 2\mathbb T, \\
\rho_0 (\dt v + v \cdot \nablah v + w \dz v ) + \nablah (c^2_s \rho_1) = \mu \deltah v \\
~~~~ ~~~~ + \lambda \nablah \dvh v + \partial_{zz}v &\text{in} ~ \Omega_h \times 2\mathbb T,\\
\dz (c^2_s \rho_1) = 0 &\text{in} ~ \Omega_h \times 2\mathbb T,
\end{cases}
\end{equation*}
where $ \mu > 0, \lambda > \frac{1}{3} \mu > 0 $.
Here the symmetry \eqref{SYM-PE} and the side condition \eqref{zero-average-pressure} are imposed.
We will make further assumptions on the viscous coefficients later.

In this section, we will study the global regularity of the solution $ (\rho_1, v, w) $ to \eqref{PE} with initial data $ v_{p,in} \in H^s(\Omega_h\times 2\mathbb T) $ for arbitrary integer $ s \in \lbrace 1, 2,3 \cdots \rbrace $, with $ v_{p,in} $ being even in the $ z $-variable and satisfying the compatible conditions:
\begin{gather*}
	\int_{\Omega_h\times 2\mathbb T}  v_{p,in} \idx = 0, ~ \dvh  \overline{v}_{p,in} = 0.
\end{gather*}

We will show the following proposition:
\begin{proposition}\label{prop:PE-Hs-data}
For $ 0 < \lambda < 4\mu < 12 \lambda  $, suppose \eqref{PE} is complemented with initial data $ v_{p,in} \in H^1(\Omega_h\times 2\mathbb T) $ as above. Then the unique solution $ v $ 
to the primitive equations \eqref{PE}, satisfies
\begin{equation*}
	\sup_{0\leq t< \infty} \norm{v(t)}{\Hnorm{1}}^2 + \int_0^\infty \biggl( \norm{\nabla v(t)}{\Hnorm{1}}^2 + \norm{\dt v(t)}{\Lnorm 2}^2 \biggr) \,dt  \leq C_{p,in},
\end{equation*}
for some positive constant $ C_{p,in} $ depending only $ \norm{v_{p,in}}{\Hnorm{1}} $. Furthermore,
\begin{equation*}
	\norm{v(t)}{\Lnorm2}^2 \leq C e^{-ct} \norm{v_{p,in}}{\Lnorm2}^2 ~~~~ t \in [0,\infty),
\end{equation*}
for some positive constants $ c, C $.
Moreover, for any integer $ s \geq 2 $, if $ v_{p,in} \in H^s(\Omega_h\times 2\mathbb T) $, there is a constant $ C_{p,in,s} $ depending only on $ \norm{v_{p,in}}{\Hnorm{s}} $  such that
\begin{equation*}
\begin{aligned}
	& \sup_{0\leq t< \infty} \blparenthese \norm{v(t)}{H^s}^2 + \norm{v_t(t)}{\Hnorm{s-2}}^2\brparenthese + \int_0^\infty \biggl( \norm{v(t)}{H^{s+1}}^2  \\
	& ~~~~ + \norm{\dt v(t)}{H^{s-1}}^2 \biggr) \,dt \leq C_{p,in,s}.
\end{aligned}
\end{equation*}
\end{proposition}

As in \cite{Li2017}, we focus on the {\it a prior} estimates below. In fact, the local-in-time regularity in Sobolev and analytic function spaces has been studied in \cite{Petcu2005}, and therefore following a continuity argument, one can check the validity of the proof below.



\subsection*{Basic energy estimate}
Take the $ L^2 $-inner product of \subeqref{PE}{2} with $ v $. We have
\begin{equation}\label{pe:basic-energy-diff}
	\dfrac{d}{dt}\biggl\lbrace \dfrac{\rho_0}{2} \norm{v}{\Lnorm2}^2 \biggr\rbrace + \mu \norm{\nablah v}{\Lnorm2}^2 + \lambda \norm{\dvh v}{\Lnorm2}^2 + \norm{\dz v}{\Lnorm2}^2 = 0.
\end{equation}
Integrating the above equation in the time variable yields
\begin{equation}\label{pe:basic-energy}
	\begin{aligned}
	& \sup_{0\leq t < \infty} \biggl\lbrace \dfrac{\rho_0}{2} \norm{v(t)}{\Lnorm2}^2 \biggr\rbrace + \int_0^\infty \bblparenthese \mu \norm{\nablah v(t)}{\Lnorm2}^2 + \lambda \norm{\dvh v(t)}{\Lnorm2}^2 \\
	& ~~~~ + \norm{\dz v(t)}{\Lnorm2}^2 \bbrparenthese \,dt = \dfrac{\rho_0}{2} \norm{v_{p,in}}{\Lnorm2}^2.
	\end{aligned}
\end{equation}
Moreover, under the assumption \eqref{conservation}, after applying the Poincar\'e inequality in \eqref{pe:basic-energy-diff}, we have the inequality
\begin{equation*}
	\dfrac{d}{dt} \norm{v}{\Lnorm{2}}^2 + c \norm{v}{\Lnorm{2}}^2 \leq 0,
\end{equation*}
for some positive constant $ c $. Thus
one can derive from above that 
\begin{equation}\label{pe:exponential-decay}
	\norm{v(t)}{\Lnorm2}^2 \lesssim e^{-ct} \norm{v_{p,in}}{\Lnorm2}^2,
\end{equation}
for all $ t \in [0,\infty) $.

\subsection*{$ H^1 $ estimate}


After applying $ \dz $ to \subeqref{PE}{2}, we write down the following equation:
\begin{equation}\label{pe:vertical-derivative}
	\begin{aligned}
	& \rho_0 ( \dt v_{z} + v \cdot\nablah v_{z} + w \dz v_{z} ) = \mu \deltah v_{z} + \lambda \nablah \dvh v_{z} + \partial_{zz} v_{z} \\
	& ~~~~ ~~~~ - \rho_0 ( v_{z} \cdot \nablah v + w_{z} \dz v ).
	\end{aligned}
\end{equation}
Then take the $ L^2 $ inner product of \eqref{pe:vertical-derivative} with $ v_{z} $. It follows, after substituting \eqref{id:vertical-vel-PE}, that
\begin{equation}\label{pe:001}
\begin{aligned}
	& \dfrac{d}{dt} \biggl\lbrace \dfrac{\rho_0}{2} \norm{v_{z}}{\Lnorm2}^2 \biggr\rbrace + \mu \norm{\nablah v_{z}}{\Lnorm2}^2 + \lambda \norm{\dvh v_{z}}{\Lnorm2}^2 + \norm{\dz v_{z}}{\Lnorm2}^2\\
	& ~~~~ = - \rho_0 \int \blparenthese v_{z}\cdot\nablah v \brparenthese \cdot v_{z} \idx + \rho_0 \int \dvh v \blparenthese \dz v \cdot v_{z} \brparenthese \idx =: L_{1} + L_{2}.
\end{aligned}
\end{equation}
After applying integration by parts, one will have
\begin{align*}
	& L_1 = \rho_0 \int \blparenthese v_{z} \cdot \nablah v_{z} \brparenthese \cdot v \idx + \rho_0 \int \blparenthese v_{}\cdot v_{z} \brparenthese \dvh{v_{z}} \idx, \\
	& L_2 = - 2 \rho_0 \int \blparenthese v \cdot\nablah v_{z} \brparenthese \cdot v_{z} \idx.
\end{align*}
Therefore, let $ p, q $ be some positive constants, to be determined later,  satisfying
\begin{equation}\label{constraint-001}
	\dfrac{1}{p} + \dfrac{1}{q} = \dfrac{1}{2}, ~~~ 2 < p < 6 ~~ (\text{equivalently} ~ q > 3 ).
\end{equation}
After applying H\"older's, the Gagliardo-Nirenberg interpolation and Young's inequalities, one has
\begin{align*}
	& L_1 + L_2 \lesssim \norm{v_{z}}{\Lnorm p}\norm{\nablah v_{z}}{\Lnorm2}\norm{v_{}}{\Lnorm q} \lesssim \norm{v_{z}}{\Lnorm2}^{3/p-1/2} \norm{\nabla v_{z}}{\Lnorm2}^{3/2-3/p} \\
	& ~~~~ \times \norm{\nablah v_{z}}{\Lnorm2} \norm{v}{\Lnorm q}
	 \lesssim \delta \norm{\nabla v_{z}}{\Lnorm2}^2 + C_\delta \norm{v_{z}}{\Lnorm2}^{2} \norm{v}{\Lnorm q}^{4p/(6-p)} \\
	& ~~ = \delta \norm{\nabla v_{z}}{\Lnorm2}^2 + C_\delta \norm{v_{z}}{\Lnorm2}^{2} \norm{v}{\Lnorm q}^{2q/(q-3)}.
\end{align*}
Hence, after choosing $ \delta > 0 $, sufficiently small, in the above estimate, we have
\begin{align*}
	& \dfrac{d}{dt} \biggl\lbrace \dfrac{\rho_0}{2} \norm{v_{z}}{\Lnorm2}^2 \biggr\rbrace + \dfrac{\mu}{2} \norm{\nablah v_{z}}{\Lnorm2}^2 + \lambda \norm{\dvh v_{z}}{\Lnorm2}^2 + \dfrac{1}{2} \norm{\dz v_{z}}{\Lnorm2}^2\\
	& ~~~~ \lesssim \norm{v_{z}}{\Lnorm2}^{2} \norm{v}{\Lnorm q}^{\frac{2q}{q-3}}.
\end{align*}
Integrating the above inequality in the time variable yields
\begin{equation}\label{pe:008}
	\begin{aligned}
	& \sup_{0\leq t < \infty} \norm{v_z(t)}{\Lnorm2}^2 + \int_0^\infty \norm{\nabla v_z(t)}{\Lnorm2}^2 \,dt \lesssim \sup_{0\leq t <\infty} \norm{v(t)}{\Lnorm q}^{\frac{2q}{q-3}} \\
	& ~~~~ \times \int_0^\infty\norm{\nabla v(t)}{\Lnorm2}^2 \,dt
	 + \norm{v_{p,in,z}}{\Lnorm2}^2
	 \lesssim  \sup_{0\leq t <\infty} \norm{v(t)}{\Lnorm q}^{\frac{2q}{q-3}} \\
	 & ~~~~ \times \norm{v_{p,in}}{\Lnorm2}^2 + \norm{v_{p,in,z}}{\Lnorm 2}^2,
	\end{aligned}
\end{equation}
for $ q \in (3,\infty) $, where we have substituted \eqref{pe:basic-energy}.



On the other hand,
after applying $ \partial_h $ to \subeqref{PE}{2}, one gets the equation:
\begin{equation}\label{pe:horizontal-derivative}
	\begin{aligned}
		& \rho_0 (\dt v_h + v \cdot\nablah v_h + w \dz v_h) + \nablah (c^2_s \rho_{1,h}) = \mu\deltah v_h \\
		& ~~~~ + \lambda\nablah \dvh v_h + \partial_{zz}v_h
		 - \rho_0 (v_h \cdot \nablah v + w_h \dz v ).
	\end{aligned}
\end{equation}
Then after taking $ L^2 $-inner product of \eqref{pe:horizontal-derivative} with $ v_h $, we have
\begin{equation}\label{pe:002}
	\begin{aligned}
		& \dfrac{d}{dt} \biggl\lbrace \dfrac{\rho_0}{2} \norm{v_h}{\Lnorm2}^2 \biggr\rbrace + \mu \norm{\nablah v_h}{\Lnorm2}^2 + \lambda\norm{\dvh v_h}{\Lnorm2}^2 + \norm{\dz v_h}{\Lnorm2}^2 \\
		& ~~~~ = - \rho_0 \int \blparenthese v_h \cdot \nablah v \brparenthese \cdot v_h \idx - \rho_0 \int w_h \blparenthese \dz v \cdot v_h \brparenthese \idx =: L_3 + L_4.
	\end{aligned}
\end{equation}
Now we estimate the terms on the right-hand side of the above equality. As before, let $ q > 3, \frac{1}{p} + \frac{1}{q} = \frac{1}{2} $. After  applying integration by parts and the H\"older, Gagliardo-Nirenberg interpolation and Young inequalities, one has
\begin{align*}
	& L_3 = \rho_0 \int \bblparenthese \blparenthese v_h \cdot\nablah v_h \brparenthese \cdot v +  ( v \cdot  v_h )  \dvh v_h \bbrparenthese \idx \lesssim \int \abs{v}{} \abs{\nablah v}{} \abs{\nablah^2 v}{} \idx\\
	& ~~~~ \lesssim \norm{\nablah v}{\Lnorm p}\norm{\nablah^2 v}{\Lnorm2}\norm{v}{\Lnorm q} \lesssim \norm{\nablah v}{\Lnorm2}^{3/p-1/2} \norm{\nabla \nablah v}{\Lnorm2}^{3/2-3/p} \\
	& ~~~~ \times  \norm{\nablah^2 v}{\Lnorm2} \norm{v}{\Lnorm q}
	 \lesssim \delta \norm{\nabla \nablah v}{\Lnorm 2}^2 + C_\delta \norm{\nablah v}{\Lnorm 2}^2 \norm{v}{\Lnorm q}^{2q/(q-3)}.
\end{align*}
On the other hand, after substituting \eqref{id:vertical-vel-PE} in the term $ L_4 $ and applying the Minkowski, H\"older's, the Gagliardo-Nirenberg interpolation and Young's inequalities, one has
\begin{align*}
	& L_4 = \rho_0 \int \bblparenthese \int_0^z \dvh v_h \,dz \bbrparenthese \times \blparenthese  \dz v \cdot v_h \brparenthese \idx \lesssim \int_0^1 \hnorm{\nablah^2 v}{\Lnorm2} \,dz \\
	& ~~~~  \times \int_0^1 \hnorm{v_z}{\Lnorm4} \hnorm{v_h}{\Lnorm 4} \,dz
	 \lesssim \norm{\nablah^2 v}{\Lnorm 2} \norm{v_z}{\Lnorm 2}^{1/2}\norm{\nablah v_z}{\Lnorm 2}^{1/2} \\
	& ~~~~ \times \norm{v_h}{\Lnorm 2}^{1/2}\norm{\nablah v_h}{\Lnorm 2}^{1/2}
	 \lesssim \delta \norm{\nabla \nablah v}{\Lnorm 2}^2 + C_\delta \norm{v_z}{\Lnorm 2}^2 \\
	 & ~~~~ \times \norm{\nablah v_z}{\Lnorm 2}^2 \norm{v_h}{\Lnorm2}^2.
\end{align*}
Similarly, take the $ L^2 $ inner product of \subeqref{PE}{2} with $ v_t $. One has,
\begin{equation}\label{pe:003}
	\begin{aligned}
		& \dfrac{d}{dt} \biggl\lbrace \dfrac{\mu}{2} \norm{\nablah v}{\Lnorm 2}^2 + \dfrac{\lambda}{2} \norm{\dvh v}{\Lnorm 2}^2 + \dfrac{1}{2} \norm{\dz v}{\Lnorm2}^2 \biggr\rbrace + \rho_0 \norm{\dt v}{\Lnorm 2}^2 \\
		& ~~~~ = - \rho_0 \int  \blparenthese v \cdot\nablah v \brparenthese  \cdot v_t \idx - \rho_0 \int w \blparenthese \dz v \cdot v_t  \brparenthese \idx =: L_{5} + L_6.
	\end{aligned}
\end{equation}
As before, applying the H\"older, Minkowski, Gagliardo-Nirenberg interpolation and Young inequalities yield, for $ q > 3, ~ \frac{1}{p} + \frac{1}{q} = \frac{1}{2} $,
\begin{align*}
	& L_5 \lesssim \norm{\nablah v}{\Lnorm p} \norm{v_t}{\Lnorm 2}\norm{v}{\Lnorm q} \lesssim  \norm{\nablah v}{\Lnorm 2}^{3/p-1/2} \norm{\nabla\nablah v}{\Lnorm 2}^{3/2-3/p} \\
	& ~~~~ \times \norm{v_t}{\Lnorm 2} \norm{v}{\Lnorm q}
	\lesssim \delta \bigl( \norm{v_t}{\Lnorm 2}^2 + \norm{\nabla\nablah v}{\Lnorm 2}^2 \bigr) \\
	& ~~~~ + C_\delta  \norm{\nablah v}{\Lnorm 2}^2 \norm{v}{\Lnorm q}^{2q/(q-3)}, \\
	& L_6 = \rho_0 \int \bblparenthese \int_0^z \dvh v \,dz \bbrparenthese \times \blparenthese \dz v \cdot v_t \brparenthese \idx \lesssim \int_0^1 \hnorm{\nablah v}{\Lnorm 4} \,dz \\
	& ~~~~ \times  \int_0^1 \hnorm{\dz v}{\Lnorm 4} \hnorm{v_t}{\Lnorm 2} \,dz
	\lesssim \norm{\nablah v}{\Lnorm 2}^{1/2}\norm{\nablah^2 v}{\Lnorm 2}^{1/2} \norm{v_z}{\Lnorm 2}^{1/2} \norm{\nablah v_z}{\Lnorm 2}^{1/2}\\
	& ~~~~ \times \norm{v_t}{\Lnorm 2}\lesssim \delta \bigl( \norm{v_t}{\Lnorm 2}^2 + \norm{\nablah^2 v}{\Lnorm 2}^2 \bigr)
	 + C_\delta \norm{v_z}{\Lnorm 2}^2 \norm{\nablah v_z}{\Lnorm 2}^2 \\
	& ~~~~ \times \norm{\nablah v}{\Lnorm 2}^2.
\end{align*}
After summing \eqref{pe:002}, \eqref{pe:003} and the estimates of $ L_3, L_4, L_5, L_6 $ above with sufficiently small $ \delta $, one has
\begin{align*}
	& \dfrac{d}{dt} \biggl\lbrace \dfrac{\rho_0}{2} \norm{v_h}{\Lnorm 2}^2 + \dfrac{\mu}{2} \norm{\nablah v}{\Lnorm 2}^2 + \dfrac{\lambda}{2} \norm{\dvh v}{\Lnorm 2}^2 + \dfrac{1}{2} \norm{\dz v}{\Lnorm 2}^2 \biggr\rbrace \\
	& ~~~~ ~~~~ + \dfrac{\mu}{2} \norm{\nablah v_h}{\Lnorm 2}^2 + \lambda\norm{\dvh v_h}{\Lnorm2}^2 + \dfrac 1 2 \norm{\dz v_h}{\Lnorm2}^2 + \dfrac{\rho_0}{2}\norm{\dt v}{\Lnorm2}^2 \\
	& ~~~~\lesssim \norm{\nablah v}{\Lnorm2}^2 \norm{v}{\Lnorm q}^{2q/(q-3)} + \norm{v_z}{\Lnorm2}^2 \norm{\nablah v_z}{\Lnorm2}^2 \norm{\nablah v}{\Lnorm2}^2.
\end{align*}
Then after applying the Gr\"onwall's inequality, it follows
\begin{equation}\label{pe:009}
\begin{aligned}
	& \sup_{0\leq t < \infty} \norm{\nabla v(t)}{\Lnorm2}^2 + \int_0^\infty \bblparenthese \norm{\nabla v_h(t)}{\Lnorm2}^2 + \norm{\dt v(t)}{\Lnorm2}^2 \bbrparenthese \,dt \\
	& ~~~~ \lesssim e^{C \int_0^\infty \norm{v_z(t)}{\Lnorm2}^2 \norm{\nablah v_z(t)}{\Lnorm2}^2 \,dt } \bblparenthese \norm{\nabla v_{p,in}}{\Lnorm2}^2 \\
	& ~~~~ + \int_0^\infty\norm{\nablah v(t)}{\Lnorm2}^2 \,dt \times \sup_{0\leq t < \infty} \norm{v(t)}{\Lnorm q}^{\frac{2q}{q-3}} \bbrparenthese,\\
	& ~~~~ \lesssim e^{C ( \sup_{0\leq t < \infty } \norm{v(t)}{\Lnorm{q}}^{\frac{4q}{q-3}} \times \norm{v_{p,in}}{\Lnorm{2}}^4 + \norm{v_{p,in,z}}{\Lnorm{2}}^4  ) } \biggl( \norm{\nabla v_{p,in}}{\Lnorm{2}}^2 \\
	& ~~~~ + \norm{v_{p,in}}{\Lnorm{2}}^2  \times \sup_{0\leq t< \infty} \norm{v(t)}{\Lnorm{q}}^{\frac{2q}{q-3}} \biggr),
\end{aligned}
\end{equation}
for some positive constant $ C $ and $ q \in (3,\infty) $, where we have substituted \eqref{pe:basic-energy} and \eqref{pe:008}.

\subsection*{$ L^q $ estimate }


We take the $ L^2 $-inner product of \subeqref{PE}{2} with $ \abs{v}{q-2} v $. It follows that,
\begin{equation}\label{pe:004}
	\begin{aligned}
		& \dfrac{d}{dt} \biggl\lbrace \dfrac{\rho_0}{q} \norm{v}{\Lnorm q}^q \biggr\rbrace + \mu \int \bblparenthese \abs{v}{q-2}\abs{\nablah v}{2} + (q-2) \abs{v}{q-2}\abs{\nablah \abs{v}{}}{2} \bbrparenthese \idx \\
		& ~~~~  + \lambda \int \abs{v}{q-2} \abs{\dvh v}{2} \idx
		 + \int \bblparenthese \abs{v}{q-2}\abs{\dz v}{2} \\
		& ~~~~ ~~~~ + (q-2)\abs{v}{q-2}\abs{\dz\abs{v}{}}{2} \bbrparenthese \idx = \int c^2_s \rho_1  \dvh (\abs{v}{q-2} v) \idx \\
		& ~~~~ - \lambda (q-2) \int \abs{v}{q-3} \blparenthese v \cdot \nablah \abs{v}{} \brparenthese \dvh v \idx =: L_{7} + L_8.
	\end{aligned}
\end{equation}
By using the Cauchy-Schwarz inequality, it holds,
\begin{align*}
	& L_8 \leq \lambda (q-2) \int \abs{v}{q-2} \abs{\nablah \abs{v}{}}{} \abs{\dvh v}{} \idx \leq \lambda \int \abs{v}{q-2}\abs{\dvh v}{2} \idx \\
	& ~~~~ + \dfrac{\lambda (q-2)^2}{4} \int \abs{v}{q-2} \abs{\nablah \abs{v}{}}{2} \idx.
\end{align*}
Therefore, \eqref{pe:004} implies
\begin{equation}\label{pe:005}
	\begin{aligned}
		& \dfrac{d}{dt} \biggl\lbrace \dfrac{\rho_0}{q} \norm{v}{\Lnorm q}^q \biggr\rbrace + \mu \norm{\abs{v}{\frac{q}{2}-1}\nablah v}{\Lnorm 2}^2 + \norm{\abs{v}{\frac{q}{2}-1}\dz v}{\Lnorm 2}^2 \\
		& ~~  \lesssim \dfrac{d}{dt} \biggl\lbrace \dfrac{\rho_0}{q} \norm{v}{\Lnorm q}^q \biggr\rbrace
		 + \int \biggl( \mu  \abs{v}{q-2}\abs{\nablah v}{2} \\
		 & ~~~~ + \bigl( \mu (q-2) - \dfrac{\lambda(q-2)^2}{4} \bigr) \ \abs{v}{q-2}\abs{\nablah \abs{v}{}}{2} \biggr) \idx  \\
		& ~~~~ + \int \biggl( \abs{v}{q-2}\abs{\dz v}{2} + (q-2)\abs{v}{q-2}\abs{\dz\abs{v}{}}{2} \biggr) \idx \leq
		L_{7},
	\end{aligned}
\end{equation}
provided
\begin{equation}\label{constraint-002}
	q-2 \geq 0 ~ \text{and} ~ \mu - \dfrac{\lambda(q-2)}{4} \geq 0,~~ \text{or equivalently}, ~~ 2 \leq q  \leq \dfrac{4\mu}{\lambda} + 2.
\end{equation}
In order to estimate $ L_7 $, we first derive an estimate for the ``pressure'' $ \rho_1 $. Recall the elliptic problem \eqref{eq:pressure-pe},
\begin{equation*}
		 - c^2_s \deltah \rho_1 = \rho_0 \int_0^1 \dvh \bigl( \dvh ( v \otimes v ) \bigr) \,dz ~~ \text{in} ~ \Omega_h, ~~ \text{with} ~ \int_{\Omega_h} \rho_1 \idxh = 0.
\end{equation*}
Now we consider the $ L^{p_1} $ estimate of $ \rho_1 $. In fact, as the consequence of the $ L^p $ estimate of the Riesz transform, one has
\begin{equation}\label{pe:006}
	\begin{aligned}
		& \norm{\rho_1}{\Lnorm{p_1}} = \hnorm{\rho_1}{\Lnorm{p_1}} \lesssim \int_0^1 \hnorm{\abs{v}{2}}{\Lnorm{p_1}} \,dz = \int_0^1 \hnorm{v}{\Lnorm{2p_1}}^2 \,dz \lesssim \int_0^1 \hnorm{v}{\Lnorm 4} \hnorm{v}{\Lnorm{4p_1/(4-p_1)}} \,dz \\
		& ~~~~ = \int_0^1 \hnorm{v}{\Lnorm 4} \hnorm{\abs{v}{\frac q 2}}{\Lnorm{\frac{8p_1}{q(4-p_1)}}}^{\frac 2 q} \,dz \lesssim \int_0^1 \hnorm{v}{\Lnorm 2}^{\frac 1 2} \hnorm{\nablah v}{\Lnorm 2}^{\frac 1 2} \hnorm{\abs{v}{\frac{q}{2}}}{\Lnorm 2}^{\frac{2}{p_1} - \frac 1 2} \hnorm{\nablah \abs{v}{\frac q 2}}{\Lnorm 2}^{\frac 1 2 + \frac 2 q - \frac{2}{p_1}} \,dz \\
		& ~~~~ \lesssim \norm{v}{\Lnorm 2}^{\frac 1 2} \norm{\nablah v}{\Lnorm 2}^{\frac 1 2} \norm{\abs{v}{\frac{q}{2}}}{\Lnorm 2}^{\frac{2}{p_1} - \frac 1 2} \norm{\nablah \abs{v}{\frac q 2}}{\Lnorm 2}^{\frac 1 2 + \frac 2 q - \frac{2}{p_1}} \\
		& ~~~~ \lesssim \norm{v}{\Lnorm 2}^{\frac 1 2} \norm{\nablah v}{\Lnorm 2}^{\frac 1 2} \norm{v}{\Lnorm q}^{\frac{q}{p_1} - \frac q 4} \norm{\abs{v}{\frac{q}{2} - 1}\nablah v}{\Lnorm 2}^{\frac 1 2 + \frac 2 q - \frac{2}{p_1}} ,
	\end{aligned}
\end{equation}
provided
\begin{equation}\label{constraint-003}
		p_1 > 2, ~~ \dfrac{1}{4} < \dfrac{1}{p_1} < \dfrac{1}{q} + \dfrac{1}{4}, ~~ q \geq 2,
\end{equation}
where we have applied the Minkowski,  H\"older and Gagliardo-Nirenberg interpolation inequalities.
Let $ q_1 > 0 $ be such that
\begin{equation*}
	\dfrac{1}{p_1} + \dfrac{1}{q_1} = \dfrac 1 2.
\end{equation*}
Then we have, after applying the Minkowski, H\"older and Gagliardo-Nirenberg interpolation inequalities,
\begin{equation}\label{Jun14}
\begin{aligned}
	& L_7 \lesssim \int  \abs{\rho_1}{} \abs{\abs{v}{\frac{q}{2}-1}\abs{\nablah v}{}}{} \abs{v}{\frac{q}{2}-1} \idx  \lesssim \hnorm{\rho_1}{\Lnorm{p_1}} \\
	& ~~~~ \times \int_0^1 \hnorm{\abs{v}{\frac q 2 - 1} \nablah v}{\Lnorm 2} \hnorm{\abs{v}{\frac{q}{2} - 1}}{\Lnorm{q_1}} \,dz
	 = \norm{\rho_1}{\Lnorm{p_1}} \\
	& ~~~~ \times  \int_0^1  \hnorm{\abs{v}{\frac q 2 - 1} \nablah v}{\Lnorm 2} \hnorm{\abs{v}{\frac{q}{2}}}{\Lnorm{\frac{q_1(q-2)}{q}}}^{\frac{q-2}{q}} \,dz \lesssim \norm{\rho_1}{\Lnorm{p_1}} \\
	& ~~~~ \times \int_0^1 \hnorm{\abs{v}{\frac q 2 - 1} \nablah v}{\Lnorm 2} \hnorm{\abs{v}{\frac q 2}}{\Lnorm 2}^{\frac{2}{q_1}} \hnorm{\nablah\abs{v}{\frac q 2}}{\Lnorm 2}^{1 - \frac{2}{q} - \frac{2}{q_1}} \,dz \lesssim \norm{\rho_1}{\Lnorm{p_1}} \\
	& ~~~~ \times \norm{\abs{v}{\frac q 2 - 1} \nablah v}{\Lnorm 2} \norm{\abs{v}{\frac q 2}}{\Lnorm 2}^{1- \frac{2}{p_1}} \norm{\nablah\abs{v}{\frac q 2}}{\Lnorm 2}^{\frac{2}{p_1} - \frac{2}{q}}\\
	& ~~ \lesssim \norm{\rho_1}{\Lnorm{p_1}} \norm{\abs{v}{\frac q 2 - 1} \nablah v}{\Lnorm 2} \norm{v}{\Lnorm q}^{\frac{q}{2}- \frac{q}{p_1}} \norm{\abs{v}{\frac{q}{2}-1}\nablah v}{\Lnorm 2}^{\frac{2}{p_1} - \frac{2}{q}} \\
	& ~~ \lesssim \norm{v}{\Lnorm 2}^{1/2}\norm{\nablah v}{\Lnorm 2}^{1/2} \norm{\abs{v}{\frac q 2 - 1} \nablah v}{\Lnorm 2}^{\frac{3}{2}} \norm{v}{\Lnorm q}^{\frac{q}{4}}\\
	& ~~ \lesssim \delta \norm{\abs{v}{\frac q 2 - 1} \nablah v}{\Lnorm 2}^2 + C_\delta \norm{v}{\Lnorm 2}^2 \norm{\nablah v}{\Lnorm 2}^2 \norm{v}{\Lnorm q}^q,
\end{aligned}
\end{equation}
provided
\begin{equation}\label{constraint-004}
	0 < \dfrac{1}{q_1} < \dfrac{1}{2} - \dfrac{1}{q}, ~~ \text{or equivalently}, ~~ \dfrac{1}{q} < \dfrac{1}{p_1}~ \text{and} ~ q > 2,
\end{equation}
where we have substituted \eqref{pe:006} in the second but last inequality.
Therefore after combining \eqref{constraint-002}, \eqref{constraint-003}, \eqref{constraint-004}, for $ q $ satisfying
\begin{equation}\label{constraint-005}
	2 < q  \leq \dfrac{4\mu}{\lambda} + 2,
\end{equation}
we conclude from \eqref{pe:005} and \eqref{Jun14},
\begin{equation*}
	\begin{aligned}
		& \dfrac{d}{dt}\biggl\lbrace \dfrac{\rho_0}{q} \norm{v}{\Lnorm q}^q \biggr\rbrace + \dfrac{\mu}{2} \norm{\abs{v}{\frac{q}{2}-1}\nablah v}{\Lnorm 2}^2 + \norm{\abs{v}{\frac{q}{2}-1}\dz v}{\Lnorm 2}^2 \\
		& ~~~~ \lesssim \norm{v}{\Lnorm 2}^2 \norm{\nablah v}{\Lnorm 2}^2 \norm{v}{\Lnorm q}^q,
	\end{aligned}
\end{equation*}
after choosing $ \delta $ sufficiently small above.
Applying Gr\"onwall's inequality to the above inequality implies, for $ 2 < q  \leq \frac{4\mu}{\lambda} + 2 $,
\begin{equation}\label{pe:007}
\begin{aligned}
	& \sup_{0\leq t < T}\norm{v(t)}{\Lnorm q}^q + \int_0^\infty \norm{\abs{v}{\frac{q}{2}-1}\nabla v(t)}{\Lnorm 2}^2 \,dt \\
	& ~~~~ \lesssim C_q e^{C \int_0^\infty \norm{v(t)}{\Lnorm 2}^2 \norm{\nablah v(t)}{\Lnorm 2}^2 \,dt } \norm{v_{p,in}}{\Lnorm q}^q \\
	& ~~~~\lesssim C_q e^{C \norm{v_{p,in}}{\Lnorm 2}^4}\norm{v_{p,in}}{\Lnorm q}^q,
\end{aligned}
\end{equation}
for some positive constant $ C $ and $ C_q $ depending on $ q $.

Therefore, after summing up the inequalities \eqref{pe:basic-energy}, \eqref{pe:008}, \eqref{pe:009}, \eqref{pe:007}, for $ \lambda < 4\mu < 12 \lambda $, 
one will have
\begin{equation}\label{pe:H^1}
\begin{aligned}
	& \sup_{0\leq t< \infty} \norm{v(t)}{\Hnorm{1}}^2 + \int_0^\infty \biggl( \norm{\nabla v(t)}{\Hnorm{1}}^2 + \norm{\dt v(t)}{\Lnorm 2}^2 \biggr) \,dt \\
	& ~~~~ \leq C_{p,in}(\norm{v_{p,in}}{\Hnorm{1}}, \norm{v_{p,in}}{\Lnorm{q}} ),
\end{aligned}	
\end{equation}
for some positive constant $ C_{p,in} $ depending on $ \norm{v_{p,in}}{\Hnorm{1}}, \norm{v_{p,in}}{\Lnorm{q}} $ with $$ 3 < q  \leq \dfrac{4\mu}{\lambda} + 2, ~~~~ \dfrac{4\mu}{\lambda} \in (1,12).$$
In particular, it suffices to take
\begin{equation*}
	q = \begin{cases}
		4 & \text{if} ~~ \dfrac{4\mu}{\lambda} \in [4 ,12) ,\\
		\dfrac{4\mu}{\lambda} + 2 & \text{if} ~~ \dfrac{4\mu}{\lambda} \in (1,4),
	\end{cases}
	~~~~ \text{such that} ~ q \in [2,6],
\end{equation*}
and therefore, \begin{equation}\label{pe:H^1-001} C_{p,in}(\norm{v_{p,in}}{\Hnorm{1}},\norm{v_{p,in}}{\Lnorm{q}}) = C_{p,in}(\norm{v_{p,in}}{\Hnorm{1}}) \end{equation} depends only on $ \norm{v_{p,in}}{\Hnorm{1}} $ by noticing that $ \norm{v_{p,in}}{\Lnorm{q}} \lesssim \norm{v_{p,in}}{\Hnorm{1}} $ in this case.

\subsection*{$ H^s $ estimates}
Next, we will show the global regularity of the solution $ v $ to system \eqref{PE} with more regular initial data $ v_{p,in} $. That is, we complement \eqref{PE} with the initial data $ v_{p,in} \in H^s(\Omega) $, with $ s \geq 2 $. In fact, we will use the mathematical induction principle to show that for any $ s \geq 1, s \in \mathbb Z^+ $,
\begin{equation}\label{pe:mi-001}
\begin{cases}
 	\sup_{0\leq t< \infty} \norm{v(t)}{H^1}^2  + \int_0^\infty \biggl( \norm{v(t)}{H^{2}}^2
	 + \norm{\dt v(t)}{\Lnorm{2}}^2 \biggr) \,dt \\
	~~~~ \leq C_{p,in,1} ~~~ \text{if} ~ s = 1,\\
	 \sup_{0\leq t< \infty} \blparenthese \norm{v(t)}{H^s}^2 + \norm{v_t(t)}{\Hnorm{s-2}}^2\brparenthese + \int_0^\infty \biggl( \norm{v(t)}{H^{s+1}}^2  \\
	~~~~ + \norm{\dt v(t)}{H^{s-1}}^2 \biggr) \,dt \leq C_{p,in,s}~~~ \text{if} ~ s \geq 2,
\end{cases}
\end{equation}
for some positive constant $ C_{p,in,s} $ depending on $ \norm{v_{p,in}}{\Hnorm{s}} $. Notice, the case when $ s = 1 $ has been shown in \eqref{pe:H^1}.

First, for any integer $ s \geq 1 $, it is assumed that \eqref{pe:mi-001} holds true. Our goal is to show that the same estimate is also true for $ s $ replaced by $ s + 1 $.
In order to do so, we apply $ \partial^{s+1} $ to \subeqref{PE}{2} with $ \partial \in \lbrace \partial_{x}, \partial_y, \partial_z \rbrace $ and denote the $k $-order derivative by $ \cdot_{k} := \partial^k \cdot $ for any $ k \in \lbrace 0,1,2 \cdots s+1 \rbrace $. Then we have the following equation:
\begin{equation}\label{pe:mi-002}
	\begin{aligned}
		& \rho_0 ( \dt v_{s+1} + v\cdot\nablah v_{s+1} + w \dz v_{s+1}) + \nablah (c^2_s \rho_{1,s+1}) \\
		& ~~~~ = \mu \deltah v_{s+1} + \lambda \nablah \dvh v_{s+1}
		 + \partial_{zz} v_{s+1} \\
		& ~~~~ - \rho_0 \bigl( \partial^{s+1} (v\cdot\nablah v) - v \cdot\nablah v_{s+1} + \partial^{s+1}(w\dz v) - w\dz v_{s+1} \bigr).
	\end{aligned}
\end{equation}
Take the $ L^2 $-inner product of \eqref{pe:mi-002} with $ v_{s+1} $. It follows,
\begin{equation}\label{pe:mi-003}
	\begin{aligned}
		& \dfrac{d}{dt} \biggl\lbrace \dfrac{\rho_0}{2} \norm{v_{s+1}}{\Lnorm 2}^2 \biggr\rbrace + \mu \norm{\nablah v_{s+1}}{\Lnorm 2}^2 + \lambda\norm{\dvh v_{s+1}}{\Lnorm 2}^2 \\
		& ~~~~ + \norm{\dz v_{s+1}}{\Lnorm 2}^2
		 = \rho_0 \int \bigl( v \cdot\nablah v_{s+1}  - \partial^{s+1} (v\cdot\nablah v)  \bigr) \cdot v_{s+1} \idx \\
		& ~~~~ + \rho_0 \int \bigl(w\dz v_{s+1} -  \partial^{s+1}(w\dz v)   \bigr)\cdot v_{s+1} \idx =: K_1 + K_2.
	\end{aligned}
\end{equation}
We estimate $ K_1, K_2 $ on the right-hand side of \eqref{pe:mi-003} below. First, notice that $ K_1, K_2 $ can be written as
\begin{equation*}
	\begin{aligned}
		& K_1 = - \rho_0 \sum_{i=0}^{s} \biggl( \begin{array}{c}
		s+1\\i
		\end{array} \biggr) 
		\int \blparenthese v_{s+1-i} \cdot\nablah v_{i} \brparenthese \cdot v_{s+1} \idx=: \sum_{i=0}^{s}K_{1,i},\\
		& K_2 = - \rho_0 \sum_{i=0}^{s} \biggl( \begin{array}{c}
		s+1\\i
		\end{array} \biggr) 
		 \int w_{s+1-i} \dz v_{i} \cdot v_{s+1} \idx =: \sum_{i=0}^{s}K_{2,i}.
	\end{aligned}
\end{equation*}
We consider the estimates of $ K_{j,i} $, for $ j\in \lbrace 1,2 \rbrace $ and $ i \in \lbrace 0 ,1 ,2 \cdots s \rbrace $ in three cases:
\begin{equation*}
	\biggl\lbrace \begin{array}{l}
		2 \leq i \leq s, \\
		i = 1, \\
		i = 0.
	\end{array}
\end{equation*}

In the case when $ i \geq 2 $,  we have
$$
s+1 - i \leq s-1, ~ 3 \leq 1 + i \leq s + 1.
$$
Therefore, applying the H\"older, Sobolev embedding and Young inequalities implies
\begin{align*}
	& K_{1,i} \lesssim \norm{v_{s+1-i}}{\Lnorm 3} \norm{\nablah v_i}{\Lnorm 6} \norm{v_{s+1}}{\Lnorm 2} \lesssim \norm{v}{H^{s+2-i}} \norm{v}{H^{i+2}}\norm{v_{s+1}}{\Lnorm 2}\\
	& ~~~~ \lesssim \norm{v}{H^s} \norm{v}{H^{s+2}} \norm{v_{s+1}}{\Lnorm 2}\lesssim \delta \norm{v}{H^{s+2}}^2  + C_\delta \norm{v}{H^s}^2\norm{v_{s+1}}{\Lnorm 2}^2.
\end{align*}
Similarly, we apply the Minkowski, H\"older, Sobolev embedding and Young inequalities to estimate $ K_{2,i}  $.
On the one hand, if $ w_{s+1-i} = \partial^{s+1-i} w = \partial_h^{s+1-i} w $, we have, thanks to \eqref{id:vertical-vel-PE},
\begin{align*}
	& K_{2,i}=  \rho_0 \biggl( \begin{array}{c}
		s+1\\i
		\end{array} \biggr)  \int \bblparenthese \partial_h^{s+1-i} \bigl(\int_0^z \dvh v \,dz'\bigr) \times \blparenthese \dz v_i \cdot v_{s+1} \brparenthese \bbrparenthese \idx \\
	& ~~ \lesssim \int_0^1 \hnorm{v_{s+2-i}}{\Lnorm 4} \,dz' \times \int_0^1 \hnorm{\dz v_i}{\Lnorm 4} \hnorm{v_{s+1}}{\Lnorm 2} \,dz \\
	& ~~ \lesssim \int_0^1 \hnorm{v_{s+2-i}}{\Lnorm 2}^{1/2}\hnorm{v_{s+3-i}}{\Lnorm 2}^{1/2} \,dz'
	\times \int_0^1 \hnorm{\dz v_i}{\Lnorm 2}^{1/2}\hnorm{\dz v_{i+1}}{\Lnorm 2}^{1/2} \hnorm{v_{s+1}}{\Lnorm 2} \,dz \\
	& ~~ \lesssim \norm{v}{H^{s+2-i}}^{1/2}\norm{v}{H^{s+3-i}}^{1/2}\norm{v}{H^{i+1}}^{1/2} \norm{v}{H^{i+2}}^{1/2} \norm{v_{s+1}}{2}\\
	& ~~ \lesssim \norm{v}{H^{s}}^{1/2}\norm{v}{H^{s+1}}\norm{v}{H^{s+2}}^{1/2}\norm{v_{s+1}}{2} \lesssim \delta \norm{v}{H^{s+2}}^2  \\
	& ~~~~ + C_\delta  \norm{v}{H^{s+1}}^2 \norm{v_{s+1}}{\Lnorm 2}^2.
\end{align*}
On the other hand, if $ w_{s+1-i} = \dz w_{s-i} $, we have
\begin{align*}
	& K_{2,i} = \rho_0 \biggl( \begin{array}{c}
		s+1\\i
		\end{array} \biggr) \int \dvh v_{s-i} \dz v_i \cdot v_{s+1}\idx \lesssim \norm{v_{s+1-i}}{\Lnorm 3} \norm{\dz v_{i}}{\Lnorm 6} \norm{v_{s+1}}{\Lnorm 2} \\
	& ~~~~ \lesssim \delta\norm{v}{H^{s+2}}^2 + C_\delta\norm{v}{H^s}^2 \norm{v_{s+1}}{\Lnorm 2}^2.
\end{align*}

In the case when $ i = 1 $, direct application of the H\"older, Sobolev embedding and Young inequalities yields,
\begin{align*}
	& K_{1,1} \lesssim \norm{v_{s}}{\Lnorm 3} \norm{v_{2}}{\Lnorm 2} \norm{v_{s+1}}{\Lnorm 6} \lesssim \norm{v}{H^{s+1}}\norm{v}{H^2} \norm{v}{H^{s+2}} \\
	& ~~~~ \lesssim \delta\norm{v}{H^{s+2}}^2 + C_\delta \norm{v}{H^2}^2 \norm{v}{H^{s+1}}^2.
\end{align*}
Meanwhile, to estimate $ K_{2,1} $, we will again apply the Minkowski, H\"older, Sobolev embedding and Young inequalities.
If $ w_s = \partial_h^s w $, we have, after substituting \eqref{id:vertical-vel-PE},
\begin{align*}
	& K_{2,1} = \rho_0 \biggl( \begin{array}{c}
		s+1\\ 1
		\end{array} \biggr) \int \bblparenthese \blparenthese \int_0^z \dvh \partial_h^s v \,dz' \brparenthese \times  \blparenthese  \dz v_1 \cdot v_{s+1} \brparenthese\bbrparenthese \idx \\
	& ~~ \lesssim \int_0^1 \hnorm{\dvh v_s}{\Lnorm 4} \,dz' \times \int_0^1 \hnorm{\dz v_1}{\Lnorm 2} \hnorm{v_{s+1}}{\Lnorm 4} \,dz \\
	& ~~ \lesssim \int_0^1 \hnorm{\dvh v_s}{\Lnorm 2}^{1/2} \hnorm{\nablah \dvh v_s}{\Lnorm 2}^{1/2} \,dz' \times \int_0^1 \hnorm{\dz v_1}{\Lnorm 2} \hnorm{v_{s+1}}{\Lnorm 2}^{1/2} \hnorm{\nablah v_{s+1}}{\Lnorm 2}^{1/2} \,dz\\
	& ~~ \lesssim \norm{v}{H^{s+1}} \norm{v}{H^{s+2}} \norm{v}{H^2} \lesssim \delta \norm{v}{H^{s+2}}^2 + C_\delta \norm{v}{H^2}^2 \norm{v}{H^{s+1}}^2.
\end{align*}
If $ w_s = \dz w_{s-1} $, we have
\begin{align*}
	& K_{2,1} = \rho_0 \biggl( \begin{array}{c}
		s+1\\ 1
		\end{array} \biggr) \int \dvh v_{s-1} \dz v_1 \cdot v_{s+1} \idx \lesssim  \norm{v_{s}}{\Lnorm3} \norm{v_{2}}{\Lnorm2} \norm{v_{s+1}}{\Lnorm 6}\\
	& ~~~~ \lesssim \delta\norm{v}{H^{s+2}}^2 + C_\delta \norm{v}{H^2}^2 \norm{v}{H^{s+1}}^2.
\end{align*}

Finally, in the case when $ i = 0 $, we apply the H\"older, Sobolev embedding and Young inequalities to get
\begin{align*}
	& K_{1,0} \lesssim \norm{v_{s+1}}{\Lnorm 6} \norm{\nablah v}{\Lnorm 3} \norm{v_{s+1}}{\Lnorm 2} \lesssim \norm{v}{H^{s+2}} \norm{v}{H^2} \norm{v_{s+1}}{\Lnorm 2} \\
	& ~~~~ \lesssim \delta \norm{v}{H^{s+2}}^2 + C_\delta \norm{v}{H^2}^2 \norm{v_{s+1}}{\Lnorm 2}^2.
\end{align*}
When $ w_{s+1} = \partial_h^{s+1} w $, applying the Minkowski, H\"older, Sobolev embedding and Young inequalities yields, after substituting \eqref{id:vertical-vel-PE},
\begin{align*}
	& K_{2,0} = \rho_0 \biggl( \begin{array}{c}
		s+1\\ 0
		\end{array} \biggr) \int \bblparenthese \partial_h^{s+1} \blparenthese \int_0^z \dvh v \,dz' \brparenthese \times \blparenthese \dz v \cdot v_{s+1} \brparenthese \bbrparenthese \idx\\
	& ~~ \lesssim \int_0^1 \hnorm{v_{s+2}}{\Lnorm 2}\,dz' \times  \int_0^1 \hnorm{\dz v}{\Lnorm 4} \hnorm{v_{s+1}}{\Lnorm 4} \,dz
	\lesssim  \int_0^1 \hnorm{v_{s+2}}{\Lnorm 2}\,dz' \\
	& ~~~~ \times \int_0^1 \hnorm{\dz v}{\Lnorm 2}^{1/2} \hnorm{\nablah \dz v}{\Lnorm 2}^{1/2} \hnorm{v_{s+1}}{\Lnorm 2}^{1/2}\hnorm{\nablah v_{s+1}}{\Lnorm 2}^{1/2} \,dz\\
	& ~~ \lesssim \norm{v}{H^{s+2}}^{3/2} \norm{v_{s+1}}{\Lnorm 2}^{1/2} \norm{v}{H^1}^{1/2}\norm{v}{H^2}^{1/2} \\
	& ~~ \lesssim \delta \norm{v}{H^{s+2}}^2 + C_\delta \norm{v}{H^1}^2 \norm{v}{H^2}^2 \norm{v_{s+1}}{\Lnorm 2}^2.
\end{align*}
When $ w_{s+1} = \dz w_{s} $, we have, after substituting \eqref{id:vertical-vel-PE},
\begin{align*}
	& K_{2,0} = \rho_0 \biggl( \begin{array}{c}
		s+1\\ 0
		\end{array} \biggr)\int \dvh v_{s} \dz v \cdot v_{s+1} \idx \lesssim \norm{v_{s+1}}{\Lnorm 6} \norm{\nabla v}{\Lnorm 3} \norm{v_{s+1}}{\Lnorm 2} \\
	& ~~~~ \lesssim \delta \norm{v}{H^{s+2}}^2 + C_\delta \norm{v}{H^2}^2 \norm{v_{s+1}}{\Lnorm 2}^2.
\end{align*}

From the above estimates, one can conclude from \eqref{pe:mi-003} that for any integer $ s \geq 1 $,
\begin{equation}\label{pe:Hs18Jun001}
\begin{aligned}
	& \dfrac{d}{dt} \biggl\lbrace \dfrac{\rho_0}{2} \norm{v_{s+1}}{\Lnorm 2}^2 \biggr\rbrace + \mu \norm{\nablah v_{s+1}}{\Lnorm 2}^2 + \lambda\norm{\dvh v_{s+1}}{\Lnorm 2}^2 \\
	& ~~~~ + \norm{\dz v_{s+1}}{\Lnorm 2}^2
	\lesssim \delta \norm{v}{H^{s+2}}^2 + C_\delta\norm{v}{H^{s+1}}^2 \norm{v}{H^{s+1}}^2 \\
	& ~~~~ + C_\delta \norm{v}{H^1}^2 \norm{v}{H^2}^2\norm{v_{s+1}}{2}^2
	 \lesssim \delta \norm{\nabla v_{s+1}}{2}^2 \\
	 & ~~~~ + \bigl(\delta + C_\delta \norm{v}{H^s}^2 \bigr) \norm{v}{H^{s+1}}^2 + C_\delta \bigl( \norm{v}{H^{s+1}}^2 + \norm{v}{H^1}^2 \norm{v}{H^2}^2 \bigr) \\
	 & ~~~~ \times \norm{v_{s+1}}{\Lnorm 2}^2.
\end{aligned}
\end{equation}
Here we have used the notation $ \norm{v_{s+1}}{\Lnorm 2}^2 $ to denote $ \sum_{\partial\in \lbrace \partial_x,\partial_y,\partial_z \rbrace} \norm{\partial^{s+1} v}{\Lnorm 2}^2 $. Then after taking $ \delta  > 0 $ small enough and applying Gr\"onwall's inequality, together with the inequalities \eqref{pe:mi-001} and \eqref{pe:H^1}, we have
\begin{equation}\label{pe:mi-004}
\begin{aligned}
	& \sup_{0\leq t < \infty} \norm{v(t)}{H^{s+1}}^2 + \int_0^\infty \norm{v(t)}{H^{s+2}}^2 \,dt \\
	& ~~ \lesssim e^{C \int_0^\infty \blparenthese \norm{v(t)}{H^{s+1}}^2 + \norm{v(t)}{H^1}^2 \norm{v(t)}{H^2}^2 \brparenthese \,dt}\bigl( \norm{v_{p,in}}{H^{s+1}}^2 \\
	& ~~~~  + \int_0^\infty (1 + \norm{v(t)}{H^s}^2) \norm{v(t)}{H^{s+1}}^2 \,dt \bigr) \lesssim e^{C_{p,in,s} + C_{p,in,1}^2}\\
	& \times \bigl( \norm{v_{p,in}}{\Hnorm{s+1}}^2
	 + C_{p,in,s} + C_{p,in,s}^2 \bigr),
\end{aligned}
\end{equation}
where $ C_{p,in,1} = C_{p,in} (\norm{v_{p,in}}{\Hnorm{1}}) $ is as in \eqref{pe:H^1-001}.

On the other hand, after replacing $ s+1 $ by $ s $ in \eqref{pe:mi-002}, we have the identity:
\begin{equation}\label{pe:mi-005}
	\begin{aligned}
		& \rho_0 \dt v_{s} = - \nablah (c^2_s \rho_{1,s}) +  \mu \deltah v_{s} + \lambda \nablah \dvh v_{s}  + \partial_{zz} v_{s} \\
		& ~~~~ - \rho_0  \partial^{s} (v\cdot\nablah v) - \rho_0 \partial^{s}(w\dz v)  .
	\end{aligned}
\end{equation}
After taking the $ L^2 $-inner product of \eqref{pe:mi-005} with $ \dt v_{s} $ and noticing the fact that
\begin{equation*}
	- \int \nablah (c^2_s \rho_{1,s}) \cdot \dt v_s \idx = \int c^2_s \rho_{1,s} \dvh \dt v_{s} \idx  = 0,
\end{equation*}
this implies
\begin{equation}\label{pe:mi-006}
\begin{aligned}
	& \norm{\dt v_s}{\Lnorm 2}^2 \lesssim \norm{\nabla v_{s+1}}{\Lnorm 2}^2 + \norm{\partial^s (v \cdot\nablah v)}{\Lnorm 2}^2 + \norm{\partial^s (w \dz v)}{\Lnorm 2}^2 \\
	& ~~~~ \lesssim  \norm{\nabla v_{s+1}}{\Lnorm 2}^2
	 + \sum_{i=0}^s \bigl( \norm{v_i\cdot\nablah v_{s-i}}{\Lnorm 2}^2 + \norm{w_i\dz v_{s-i}}{\Lnorm 2}^2 \bigr) \\
	 & ~~~~  \lesssim  (1 + \norm{v}{H^{s+1}}^2)  \norm{v}{H^{s+2}}^2,
\end{aligned}
\end{equation}
where we have applied the inequalities below: since $ s \geq 1 $,
\begin{align*}
	& \norm{v_i\cdot\nablah v_{s-i}}{\Lnorm 2}^2 \lesssim \norm{v_i}{\Lnorm 3}^2 \norm{\nablah v_{s-i}}{\Lnorm 6}^2 \lesssim \norm{v}{H^{i+1}}^2 \norm{v}{H^{s+2-i}}^2\\
	& ~~~~ \lesssim \norm{v}{H^{s+1}}^2 \norm{v}{H^{s+2}}^2,\\
	& \norm{w_i\dz v_{s-i}}{\Lnorm 2}^2 
	\lesssim \norm{w}{H^{i+1}}^2 \norm{v}{H^{s+2-i}}^2 \lesssim \norm{v}{H^{i+2}}^2 \norm{v}{H^{s+2-i}}^2 \\
	& ~~~~ \lesssim \norm{v}{H^{s+1}}^2 \norm{v}{H^{s+2}}^2,
\end{align*}
due to the fact that from \eqref{id:vertical-vel-PE},
\begin{equation*}
	\norm{w}{H^{i+1}} = \norm{\int_0^z \dvh v\,dz'}{H^{i+1}} \lesssim \norm{v}{H^{i+2}}.
\end{equation*}
Here we have applied the Minkowski, H\"older and Sobolev embedding inequalities. Similarly, taking $ s = s-1 $ in \eqref{pe:mi-006} yields
\begin{equation}\label{pe:mi-007}
	\norm{\dt v_{s-1}}{\Lnorm{2}}^2 \lesssim \blparenthese 1 + \norm{v}{\Hnorm{s}}^2 \brparenthese \norm{v}{\Hnorm{s+1}}^2.
\end{equation}
Integrating \eqref{pe:mi-004} in the time variable, together with \eqref{pe:mi-006}, \eqref{pe:mi-007}, implies  \eqref{pe:mi-001} with $ s $ replaced by $ s + 1 $. This finishes the mathematical induction.
Hence, this concludes the proof of Proposition \ref{prop:PE-Hs-data}.

\paragraph{Acknowledgements}
	This work was supported in part by the Einstein Stiftung/Foundation - Berlin, through the Einstein Visiting Fellow Program, and by the John Simon Guggenheim Memorial Foundation


\begin{thebibliography}{10}
	
	\bibitem{AlazardReview}
	Thomas Alazard.
	\newblock {Incompressible limit of the non-isentropic Euler equations}.
	\newblock {Available at \href{http://talazard.perso.math.cnrs.fr/Ade.pdf}{http://talazard.perso.math.cnrs.fr/Ade.pdf}}.
	
	\bibitem{Alazard2005}
	Thomas Alazard.
	\newblock {Low Mach number flows and combustion}.
	\newblock {\em SIAM J. Math. Anal.}, 38(4):1186--1213, 2006.
	
	\bibitem{Alazard2006}
	Thomas Alazard.
	\newblock {Low Mach number limit of the full Navier--Stokes equations}.
	\newblock {\em Arch. Ration. Mech. Anal.}, 180(1):1--73, 2006.
	
	\bibitem{Azerad2001}
	Pascal Az{\'{e}}rad and Francisco Guill{\'{e}}n.
	\newblock {Mathematical justification of the hydrostatic approximation in the
		primitive equations of geophysical fluid dynamics}.
	\newblock {\em SIAM J. Math. Anal.}, 33(4):847--859, 2001.
	
	\bibitem{Brenier1999}
	Yann Brenier.
	\newblock {Homogeneous hydrostatic flows with convex velocity profiles}.
	\newblock {\em Nonlinearity}, 12(3):495--512, 1999.
	
	\bibitem{Cao2015}
	Chongsheng Cao, Slim Ibrahim, Kenji Nakanishi, and Edriss~S. Titi.
	\newblock {Finite-time blowup for the inviscid primitive equations of oceanic
		and atmospheric dynamics}.
	\newblock {\em Commun. Math. Phys.}, 337(2):473--482, 2015.
	
	\bibitem{Cao2014}
	Chongsheng Cao, Jinkai Li, and Edriss~S. Titi.
	\newblock {Global well-posedness of strong solutions to the 3D primitive
		equations with horizontal eddy diffusivity}.
	\newblock {\em J. Differential Equations}, 257(11):4108--4132, 2014.
	
	\bibitem{Cao2014b}
	Chongsheng Cao, Jinkai Li, and Edriss~S. Titi.
	\newblock {Local and global well-posedness of strong solutions to the 3D
		primitive equations with vertical eddy diffusivity}.
	\newblock {\em Arch. Ration. Mech. Anal.}, 214(1):35--76, 2014.
	
	\bibitem{Cao2016}
	Chongsheng Cao, Jinkai Li, and Edriss~S. Titi.
	\newblock {Global well-posedness of the three-dimensional primitive equations
		with only horizontal viscosity and diffusion}.
	\newblock {\em Commun. Pure Appl. Math.}, 69(8):1492--1531, 2016.
	
	\bibitem{Cao2016a}
	Chongsheng Cao, Jinkai Li, and Edriss~S. Titi.
	\newblock {Strong solutions to the 3D primitive equations with only horizontal
		dissipation: Near $H^1$ initial data}.
	\newblock {\em Journal of Functional Analysis}, 272(11):4606--4641, 2017.
	
	\bibitem{Cao2017}
	Chongsheng Cao, Jinkai Li, and Edriss~S. Titi.
	\newblock {Global well-posedness of the 3D primitive equations with horizontal
		viscosity and vertical diffusivity}.
	\newblock Available at \href{https://arxiv.org/abs/1703.02512}{arXiv:1703.02512}, 2017.
	
	\bibitem{Cao2007}
	Chongsheng Cao and Edriss S. Titi.
	\newblock {Global well-posedness of the three-dimensional viscous primitive
		equations of large scale ocean and atmosphere dynamics}.
	\newblock {\em Ann. Math.}, 166(1):245--267, 2007.
	
	\bibitem{Cao2003}
	Chongsheng Cao and Edriss~S. Titi.
	\newblock {Global well-posedness and finite-dimensional global attractor for a
		3-D planetary geostrophic viscous model}.
	\newblock {\em Commun. Pure Appl. Math.}, 56(2):198--233, 2003.
	
	\bibitem{Cao2012}
	Chongsheng Cao and Edriss~S. Titi.
	\newblock {Global well-posedness of the 3D primitive equations with partial
		vertical turbulence mixing heat diffusion}.
	\newblock {\em Commun. Math. Phys.}, 310(2):537--568, 2012.
	
	\bibitem{Zelati2015}
	Michele~Coti Zelati, Aimin Huang, Igor Kukavica, Roger Temam, and Mohammed
	Ziane.
	\newblock {The primitive equations of the atmosphere in presence of vapour
		saturation}.
	\newblock {\em Nonlinearity}, 28(3):625--668, 2015.
	
	\bibitem{Danchin2002per}
	Rapha\"{e}l~Danchin.
	\newblock {Zero Mach number limit for compressible flows with periodic boundary
		conditions}.
	\newblock {\em Am. J. Math.}, 124(6):1153--1219, 2002.
	
	\bibitem{Danchin2002}
	Rapha\"{e}l~Danchin.
	\newblock {Zero Mach number limit in critical spaces for compressible
		Navier--Stokes equations}.
	\newblock {\em Ann. Scient. {\'{E}}c. Norm. Sup.}, 35(1):27--75,
	2002.
	
	\bibitem{Danchin2005}
	Rapha\"{e}l Danchin.
	\newblock {Low Mach number limit for viscous compressible flows}.
	\newblock {\em ESAIM Math. Model. Numer. Anal.}, 39(3):459--475, 2005.
	
	\bibitem{DesjardinsGrenier1999}
	Benoit~Desjardins and Emmanuel~Grenier.
	\newblock {Low Mach number limit of viscous compressible flows in the whole
		space}.
	\newblock {\em Proc. R. Soc. A Math. Phys. Eng. Sci.}, 455(1986):2271--2279,
	1999.
	
	\bibitem{Ersoy2012}
	Mehmet Ersoy and Timack Ngom.
	\newblock {Existence of a global weak solution to compressible primitive
		equations}.
	\newblock {\em C. R. Acad. Sci. Paris, Ser. I}, 350(7-8):379--382, 2012.
	
	\bibitem{Ersoy2011a}
	Mehmet Ersoy, Timack Ngom, and Mamadou Sy.
	\newblock {Compressible primitive equations: Formal derivation and stability of
		weak solutions}.
	\newblock {\em Nonlinearity}, 24(1):79--96, 2011.
	
	\bibitem{Feireisl2004}
	Eduard Feireisl.
	\newblock {\em {Dynamics of Viscous Compressible Fluids}}.
	\newblock Oxford Lecture Series in Mathematics and its Applications, 26. Oxford
	University Press, 2004.
	
	\bibitem{Feireisl2011}
	Eduard Feireisl.
	\newblock {Flows of viscous compressible fluids under strong stratification: Incompressible limits for long-range potential forces}.
	\newblock {\em Math. Model. Methods Appl. Sci.}, 21(01):7--27, 2011.
	
	\bibitem{Feireisl2015a}
	Eduard Feireisl, Rupert Klein, Anton{\'{i}}n Novotn{\'{y}}, and Ewelina
	Zatorska.
	\newblock {On singular limits arising in the scale analysis of stratified fluid
		flows}.
	\newblock {\em Math. Model. Methods Appl. Sci.}, 26(03):419--443, 2016.
	
	\bibitem{Feireisl2008}
	Eduard Feireisl, Josef M{\'{a}}lek, Anton{\'{i}}n Novotn{\'{y}}, and Ivan
	Stra{\v{s}}kraba.
	\newblock {Anelastic approximation as a singular limit of the compressible
		Navier--Stokes system}.
	\newblock {\em Communications in Partial Differential Equations}, 33(1):157--176, 2008.
	
	\bibitem{FeireislSingularLimits}
	Eduard Feireisl and Anton{\'{i}}n Novotn{\'{y}}.
	\newblock {\em {Singular Limits in Thermodynamics of Viscous Fluids}}.
	\newblock Advances in Mathematical Fluid Mechanics. Springer International
	Publishing, Cham, 2017.

	
	\bibitem{Efeireisl2012}
	Eduard Feireisl and Maria~E. Schonbek.
	\newblock {On the Oberbeck-Boussinesq approximation on unbounded domains}.
	\newblock In {\em Nonlinear Partial Differ. Equations}, pages 131--168.
	Springer Berlin Heidelberg, Berlin, Heidelberg, 2012.
	
	\bibitem{Gallagher1998}
	Isabelle~Gallagher.
	\newblock {Applications of Schochet's methods to parabolic equations}.
	\newblock {\em J. Math. Pures Appl.}, 77(10):989--1054, 1998.
	
	\bibitem{Gatapov2005}
	Bair~V. Gatapov and Aleksandr~V. Kazhikhov.
	\newblock {Existence of a global solution to one model problem of atmosphere
		dynamics}.
	\newblock {\em Sib. Math. J.}, 46(5):805--812, 2005.
	
	\bibitem{Gerard-Varet2018}
	David Gerard-Varet, Nader Masmoudi, and Vlad Vicol.
	\newblock {Well-posedness of the hydrostatic Navier--Stokes equations}.
	\newblock {\em J. Math. Fluid Mech.}, 14(2):355--361, 2018.
	
	\bibitem{Ginibre1995}
	Jean~Ginibre and Giorgio~Velo.
	\newblock {Generalized Strichartz inequalities for the wave equation}.
	\newblock {\em Journal of Functional Analysis}, 133(1):50--68, 1995.
	
	\bibitem{GuillenGonzalez2001}
	Francisco M. Guill\'{e}n-Gonz\'{a}lez, Nader~Masmoudi, and Mar\'{i}a \'{A}. Rodr\'{i}guez-Bellido.
	\newblock {Anisotropic estimates and strong solutions of the primitive equations}.
	\newblock {\em Differential and Integral Equations}, 14(11):1381--1408, 2001.
	
	\bibitem{Hieber2016}
	Matthias Hieber and Takahito Kashiwabara.
	\newblock {Global strong well-posedness of the three dimensional primitive
		equations in ${L^p}$-Spaces}.
	\newblock {\em Arch. Ration. Mech. Anal.}, 221(3):1077--1115, 2016.
	
	\bibitem{hittmeir2017}
	{Sabine Hittmeir}, {Rupert Klein}, {Jinkai Li}, and {Edriss~S. Titi}.
	\newblock {Global well-posedness for passively transported nonlinear moisture
		dynamics with phase changes}.
	\newblock {\em Nonlinearity}, 30:3676--3718, 2017.
	
	\bibitem{Hoff1998}
	David Hoff.
	\newblock {The zero-mach limit of compressible flows}.
	\newblock {\em Commun. Math. Phys.}, 192:543--554, 1998.
	
	\bibitem{HuTemamZiane2003}
	Changbing Hu, Roger Temam, and Mohammed Ziane.
	\newblock {The primitive equations on the large scale ocean under the small
		depth hypothesis}.
	\newblock {\em Discret. Contin. Dyn. Syst.}, 9(1):97--131, 2002.
	
	\bibitem{Ignatova2012}
	Mihaela Ignatova, Igor Kukavica, and Mohammed Ziane.
	\newblock {Local existence of solutions to the free boundary value problem for
		the primitive equations of the ocean}.
	\newblock {\em J. Math. Phys.}, 53(10):103101, 2012.
	
	\bibitem{Jiang2011}
	Song Jiang and Yaobin Ou.
	\newblock {Incompressible limit of the non-isentropic Navier--Stokes equations
		with well-prepared initial data in three-dimensional bounded domains}.
	\newblock {\em J. Math. Pures Appl.}, 96(1):1--28, 2011.
	
	\bibitem{Keel1998}
	Markus Keel and Terence Tao.
	\newblock {Endpoint Strichartz estimates}.
	\newblock {\em Am. J. Math.}, 120(5):955--980, 1998.
	
	\bibitem{Klainerman1981}
	Sergiu Klainerman and Andrew Majda.
	\newblock {Singular limits of quasilinear hyperbolic systems with large
		parameters and the incompressible limit of compressible fluids}.
	\newblock {\em Commun. Pure Appl. Math.}, 34(4):481--524, 1981.
	
	\bibitem{Klainerman1982}
	Sergiu Klainerman and Andrew Majda.
	\newblock {Compressible and incompressible fluids}.
	\newblock {\em Commun. Pure Appl. Math.}, 35(5):629--651, 1982.
	
	\bibitem{Klein2001}
	Ruper~Klein, Nicola~Botta, Thomas~Schneider, Claus-Dieter Munz, Sabine~Roller, Andreas~Meister,
	L.~Hoffmann, and Thomas~Sonar.
	\newblock {Asymptotic adaptive methods for multi-scale problems in fluid
		mechanics}.
	\newblock {\em J. Eng. Math.}, 39(1/4):261--343, 2001.
	
	\bibitem{Klein2000}
	Rupert Klein.
	\newblock {Asymptotic analyses for atmospheric flows and the construction of
		asymptotically adaptive numerical methods}.
	\newblock {\em ZAMM}, 80(11-12):765--777, 2000.
	
	\bibitem{Klein2005}
	Rupert Klein.
	\newblock {Multiple spatial scales in engineering and atmospheric low Mach
		number flows}.
	\newblock {\em ESAIM Math. Model. Numer. Anal.}, 39(3):537--559, 2005.
	
	\bibitem{RKlein2010}
	Rupert Klein.
	\newblock {Scale-dependent models for atmospheric flows}.
	\newblock {\em Annu. Rev. Fluid Mech.}, 42(1):249--274, 2010.
	
	\bibitem{Kobelkov2006}
	Georgij~M. Kobelkov.
	\newblock {Existence of a solution `in the large' for the 3D large-scale ocean
		dynamics equations}.
	\newblock {\em C. R. Acad. Sci. Paris, Ser. I}, 343(4):283--286, 2006.
	
	\bibitem{Kukavica2014}
	Igor Kukavica, Yuan Pei, Walter Rusin, and Mohammed Ziane.
	\newblock {Primitive equations with continuous initial data}.
	\newblock {\em Nonlinearity}, 27(6):1135--1155, 2014.
	
	\bibitem{Kukavica2011}
	Igor Kukavica, Roger Temam, Vlad~C. Vicol, and Mohammed Ziane.
	\newblock {Local existence and uniqueness for the hydrostatic Euler equations
		on a bounded domain}.
	\newblock {\em J. Differential  Equations}, 250(3):1719--1746, 2011.
	
	\bibitem{Kukavica2007}
	Igor Kukavica and Mohammed Ziane.
	\newblock {On the regularity of the primitive equations of the ocean}.
	\newblock {\em Nonlinearity}, 20(12):2739--2753, 2007.
	
	\bibitem{Kukavica2007a}
	Igor Kukavica and Mohammed Ziane.
	\newblock {The regularity of solutions of the primitive equations of the ocean
		in space dimension three}.
	\newblock {\em C. R. Acad. Sci. Paris, Ser. I}, 345(5):257--260, 2007.
	
	\bibitem{Li2017a}
	Jinkai Li and Edriss~S Titi.
	\newblock {Existence and uniqueness of weak solutions to viscous primitive
		equations for a certain class of discontinuous initial data}.
	\newblock {\em SIAM J. Math. Anal.}, 49(1):1--28, 2017.
	
	\bibitem{Li2017}
	Jinkai Li and Edriss~S. Titi.
	\newblock {The primitive equations as the small aspect ratio limit of the
		Navier--Stokes equations: Rigorous justification of the hydrostatic
		approximation}.
	\newblock {\em J. Math. Pures Appl.}, pages 1--31, 2018.
	
	\bibitem{Lindblad1995}
	Hans~Lindblad and Christopher D. Sogge.
	\newblock {On existence and scattering with minimal regularity for semilinear
		wave equations}.
	\newblock {\em Journal of Functional Analysis}, 130(2):357--426, 1995.
	
	\bibitem{Lions1993}
	Jacques-Louis Lions, Roger~Temam, and Shouhong~Wang.
	\newblock {Models of the Coupled Atmosphere and Ocean (CAO I) Part I}.
	\newblock {\em Comput. Mech. Adv.}, 1(1):5 -- 54, 1993.
	
	\bibitem{JLLions1994}
	Jacques-Louis Lions, Roger~Temam, and Shouhong~Wang.
	\newblock {Geostrophic asymptotics of the primitive equations of the
		atmosphere}.
	\newblock {\em Topol. Methods Nonlinear Anal.}, 4:253 -- 287, 1994.
	
	\bibitem{JLLions1992}
	Jacques-Louis Lions, Roger~Temam, and Shouhong Wang.
	\newblock {On the equations of the large-scale ocean}.
	\newblock {\em Nonlinearity}, 5(5):1007--1053, 1992.
	
	\bibitem{Lions1992}
	Jacques-Louis Lions, Roger Temam, and Shouhong Wang.
	\newblock {New formulations of the primitive equations of atmosphere and
		applications}.
	\newblock {\em Nonlinearity}, 5(2):237--288, 1992.
	
	\bibitem{Lions1998a}
	Pierre-Louis Lions and Nader~Masmoudi.
	\newblock {Incompressible limit for a viscous compressible fluid}.
	\newblock {\em J. Math. Pures Appl.}, 77(6):585--627, 1998.
	
	\bibitem{Lions1996}
	Pierre-Louis Lions.
	\newblock {\em {Mathematical Topics in Fluid Mechanics. Volume 1.
			Incompressible Models}}.
	\newblock Oxford Lecture Series in Mathematics and Its Applications, 3. Oxford
	University Press, 1996.
	
	\bibitem{Lions1998}
	Pierre-Louis Lions.
	\newblock {\em {Mathematical Topics in Fluid Mechanics. Volume 2. Compressible
			Models}}.
	\newblock Oxford Lecture Series in Mathematics and Its Applications , Vol 2, No
	10. Oxford University Press, 1998.
	
	\bibitem{LT2018b}
	Xin Liu and Edriss~S. Titi.
	\newblock {Global existence of weak solutions to compressible primitive
		equations with degenerate viscosities}.
	\newblock Available at \href{https://arxiv.org/abs/1808.03975}{arXiv:1808.03975}, accepted by {\em SIAM J. Math. Anal.}, 2019.
	
	\bibitem{LT2018a}
	Xin Liu and Edriss~S. Titi.
	\newblock {Local well-posedness of strong solutions to the three-dimensional
		compressible primitive equations}.
	\newblock Available at \href{https://arxiv.org/abs/1806.09868}{arXiv:1806.09868}, 2018.
	
	\bibitem{CPE2PE2}
	Xin Liu and Edriss~S. Titi.
	\newblock{Zero Mach number limit of the compressible primitive equations Part II: Ill-prepared initial data}.
	\newblock In preperation.
	
	\bibitem{MajdaAtmosphereOcean}
	Andrew Majda.
	\newblock {\em {Introduction to PDEs and Waves for the Atmosphere and Ocean}}.
	\newblock Courant Lecture Notes in Mathematics 9. American Mathematical
	society, 2003.
	
	\bibitem{Majda2002vorticity}
	Andrew~J. Majda and Andrea~L. Bertozzi.
	\newblock {\em {Vorticity and Incompressible Flow}}.
	\newblock Number~27. Cambridge University Press, 2002.
	
	\bibitem{Masmoudi2001}
	Nader Masmoudi.
	\newblock {Incompressible, inviscid limit of the compressible Navier--Stokes
		system}.
	\newblock {\em Ann. Inst. Henri Poincar\'e, Anal. non lin\'eaire}
	18(2):199--224, 2001.
	
	\bibitem{Masmoudi2007}
	Nader Masmoudi.
	\newblock {Rigorous derivation of the anelastic approximation}.
	\newblock {\em J. Math. Pures Appl.}, 88(3):230--240, 2007.
	
	\bibitem{Wong2012}
	Nader Masmoudi and Tak~Kwong Wong.
	\newblock {On the $H^s$ theory of hydrostatic euler equations}.
	\newblock {\em Arch. Ration. Mech. Anal.}, 204(1):231--271, 2012.
	
	\bibitem{Metivier2001}
	Guy~M{\'{e}}tivier and Steve~Schochet.
	\newblock {The incompressible limit of the non-isentropic Euler equations}.
	\newblock {\em Arch. Ration. Mech. Anal.}, 158(1):61--90, 2001.
	
	\bibitem{Novotny2011}
	Anton\'in Novotn\'y, Michael R\r{u}\v{z}i\v{c}ka, and Gudrun Th\"{a}ter.
	\newblock {Rigorous derivation of the anelastic approximation to the
		Oberbeck--Boussinesq equations}.
	\newblock {\em Asymptot. Anslysis}, 75:93--123, 2011.
	
	\bibitem{Petcu2005}
	Madalina~Petcu and Djoko~Wirosoetisno.
	\newblock {Sobolev and Gevrey regularity results for the primitive equations in
		three space dimensions}.
	\newblock {\em Appl. Anal.}, 84(8):769--788, 2005.
	
	\bibitem{Rajagopal1996}
	K.R. Rajagopal, M.~Ruzicka, and A.R. Srinivasa.
	\newblock {On the Oberbeck-Boussinesq approximation}.
	\newblock {\em Math. Model. Methods Appl. Sci.}, 06(08):1157--1167, 1996.
	
	\bibitem{Renardy2009}
	Michael Renardy.
	\newblock {Ill-posedness of the hydrostatic Euler and Navier--Stokes equations}.
	\newblock {\em Arch. Ration. Mech. Anal.}, 194(3):877--886, 2009.
	
	\bibitem{Richardson1965}
	Lewis~F. Richardson.
	\newblock {\em {Weather Prediction by Numerical Process}}.
	\newblock Cambridge University Press, 2007.
	
	\bibitem{Schochet1994}
	Steve~Schochet.
	\newblock {Fast singular limits of hyperbolic PDEs}.
	\newblock {\em J. Differential Equations}, 114(2):476--512, 1994.
	
	\bibitem{SchochetCMP1986}
	Steve Schochet.
	\newblock {The compressible Euler equations in a bounded domain: Existence of
		solutions and the incompressible limit}.
	\newblock {\em Commun. Math. Phys.}, 104(1):49--75, 1986.
	
	\bibitem{Schochet1988}
	Steven Schochet.
	\newblock {Asymptotics for symmetric hyperbolic systems with a large
		parameter}.
	\newblock {\em J. Differential Equations}, 75(1):1--27, 1988.
	
	\bibitem{Tang2015}
	Tong Tang and Hongjun Gao.
	\newblock {On the stability of weak solution for compressible primitive
		equations}.
	\newblock {\em Acta Appl. Math.}, 140(1):133--145, 2015.
	
	\bibitem{Ukai1986}
	Seiji Ukai.
	\newblock {The incompressible limit and the initial layer of the compressible
		Euler equation}.
	\newblock {\em J. Math. Kyoto Univ.}, 26(2):323--331, 1986.
	
	\bibitem{wang2017global}
	{Fengchao Wang, Changsheng Dou, Quansen Jiu}.
		\newblock {Global weak solutions to 3D compressible primitive equations with density-dependent viscosity}.
		\newblock Available at \href{https://arxiv.org/abs/1712.04180}{arXiv:1712.04180},
		{2017}.

	
	\bibitem{Washington2005}
	Warren~M. Washington and Claire~L. Parkinson.
	\newblock {\em {An Introduction to Three-Dimensional Climate Modeling}}.
	\newblock University Science Books, 2005.
	
	\bibitem{Wong2014}
	Tak~Kwong Wong.
	\newblock {Blowup of solutions of the hydrostatic Euler equations}.
	\newblock {\em Proc. Am. Math. Soc.}, 143(3):1119--1125, 2014.
	
	\bibitem{wk}
	Aneta Wr{\'{o}}blewska-Kami{\'{n}}ska.
	\newblock {The asymptotic analysis of the complete fluid system on a varying
		domain: from the compressible to the incompressible flow}.
	\newblock {\em SIAM J. Math. Anal.}, 49(5):3299--3334, 2017.


	
\end{thebibliography}

\end{document}